\UseRawInputEncoding
\documentclass[12pt]{amsart}
\usepackage{amssymb}
\usepackage[all]{xy}
\usepackage[toc,page,title,titletoc,header]{appendix}
\usepackage{xcolor}

\theoremstyle{plain}
\newtheorem{thm}{Theorem}[section]

\newtheorem*{thm1.1}{Theorem 1.1}

\newtheorem{lem}[thm]{Lemma}
\newtheorem{cor}[thm]{Corollary}
\newtheorem{pro}[thm]{Proposition}

\theoremstyle{definition}

\newtheorem{rem}[thm]{Remark}

\newtheorem{defi}[thm]{Definition}
\newtheorem{exe}[thm]{Example}
\newtheorem{facts}[thm]{Facts}

\newtheorem{que}[thm]{Question}

\numberwithin{equation}{section}
\newcounter{elno}                

\newcommand{\Tan}{\rm Tan}

\newcommand{\la}{\lambda}

\newcommand{\an}{{\rm an}}

\newcommand{\Spec}{{\rm Spec \,}}

\newcommand{\Char}{{\rm char}}

\newcommand{\ord}{{\rm ord}}
\newcommand{\id}{{\rm id}}

\newcommand{\Proj}{{\rm Proj}}
\newcommand{\Gal}{{\rm Gal}}

\newcommand{\End}{{\rm End}}
\newcommand{\Per}{{\rm Per}\,}
\newcommand{\Rad}{{\rm Rad}}

\renewcommand{\d}{\stackrel{\mbox{\scriptsize{$\bullet$}}}{}}

\newcommand{\boxtensor}{{\Box\kern-9.03pt\raise1.42pt\hbox{$\times$}}}

\newcommand{\propsubset}
{\mbox{$\textstyle{
\subseteq_{\kern-5pt\raise-1pt\hbox{\mbox{\tiny{$/$}}}}}$}}

\newcommand{\sC}{{\mathcal C}}
\newcommand{\sD}{{\mathcal D}}

\newcommand{\sG}{{\mathcal G}}
\newcommand{\sH}{{\mathcal H}}
\newcommand{\sI}{{\mathcal I}}

\newcommand{\sM}{{\mathcal M}}

\newcommand{\sO}{{\mathcal O}}
\newcommand{\sP}{{\mathcal P}}

\newcommand{\sS}{{\mathcal S}}
\newcommand{\sT}{{\mathcal T}}
\newcommand{\sU}{{\mathcal U}}
\newcommand{\sV}{{\mathcal V}}
\newcommand{\sW}{{\mathcal W}}

\newcommand{\sY}{{\mathcal Y}}

\newcommand{\A}{{\mathbb A}}

\newcommand{\C}{{\mathbb C}}

\newcommand{\F}{{\mathbb F}}

\renewcommand{\P}{{\mathbb P}}
\newcommand{\Q}{{\mathbb Q}}
\newcommand{\R}{{\mathbb R}}

\newcommand{\Z}{{\mathbb Z}}
\newcommand{\bk}{{\mathbf{k}}}

\newcommand{\Fix}{\mathrm{Fix}}

\begin{document}
\title[]{Algebraicity criteria, invariant subvarities and transcendence problems from arithmetic dynamics}

\author{Junyi Xie}

\address{BICMR, Peking University, Haidian District, Beijing 100871, China}
\email{xiejunyi@bicmr.pku.edu.cn}


\date{\today}

\bibliographystyle{plain}


\maketitle

\begin{abstract}
We introduce an algebraicity criteria. It has the following form: under certain conditions, an analytic subvariety of some algebriac variety over a global field $K$, if it contains many $K$-points, then it is algebraic over $K.$ This gives a way to show the transcendence of points via the transcendence of analytic subvarieties.
Such a situation often appears when we have a dynamical system, because we can often produce infinitely many points from one point via iterates.

Combing this  criteria and the study of invariant subvarieties, we get some results on the transcendence in arithmetic dynamics. 
We get a characterization of products of B\"ottcher coordinates or products of multiplicative canonical heights for polynomial dynamical pairs to be algebraic. For this, we studied the invariant subvarieties for product of endomorphisms. In particular, we partially generate Medvedev-Scanlon's classification of invariant subvarieties of split polynomial maps to separable endomorphisms on $(\P^1)^N$ in any characteristic. We also get some hight dimensional partial generalization via introducing of a notion of 
independence.
We then study dominant endomorphisms $f$ on $\A^N$ over a number field of algebraic degree $d\geq 2$. We show that in most of the cases (e.g. when such an endomorphism extends to an endomorphism on $\P^N$), there are many analytic curves centered at infinity which are periodic. We show that for most of them, it is algebraic if and only if it contains one algebraic point. We also study the periodic curves. We show that for most of $f$, all periodic curves has degree at most $2$. When $N=2$, we get a more precise classification result. We show that under a condition which is satisfied for a general $f$, if $f$ has infinitely many periodic curves, then $f$ is homogenous up to a changing of origin.
\end{abstract}

\tableofcontents

\section{Introduction}
Let us consider the following naive question:  
\begin{que}\label{quenaive}Let $C$ be an analytic curve in $\C^2$ passing through the origin $o=(0,0).$
Let $p_n, n\geq 0$ be a sequence of points in $K^2\subseteq \C^2$ for some number field $K$ satisfying $p_n\in \C^2$ and $\lim\limits_{n\geq 0}p_n=o$. Is $C$ contained in an algebraic curve?
\end{que}

Without any assumption, the answer is negative. 
\begin{exe}\label{exenonalg}For $n\geq 2$, set $h_n(x):=x\prod_{i=2}^n(1-ix).$
For $|x|<r<1$, there is $c(r)>0$ such that $$h_n(x)/n!< x(\prod_{i=2}^n(1/i+x))\leq c(r)(\frac{1+r}{2})^n.$$
Hence $$f(x):=\sum_{n\geq 2}h_n(x)/n!$$ converges for $|x|<1.$
For every $n=2,3,\dots$, $f(1/n)=\sum_{i= 2}^nh_i(1/n)/i!\in \Q.$
Set $p_{1,n}:=(1/n, f(1/n))\in \Q^2, n\geq 2.$
Let $C_1$ be the analytic curve $\{y=f(x), |x|<1\}$ in $\C^2.$ 
It is clear that $p_{1,n}\to 0$ for $n\to 0$ and $p_{1,n}\in C_1$ for $n\geq 2.$
But $C_1$ is not contained in any algebraic curve. Otherwise $f'(0)$ should be an algebraic number. But 
$$f'(0)=\sum_{n\geq 2}h_n'(0)/n!=\sum_{n\geq 2}1/n!=e-2,$$
which is transcendental.
Even we assume that the formal curve induced by $C$ at $o$ is defined over $K$, Question \ref{quenaive} still has a negative answer. 
Set $$g(x):=\sum_{n\geq 2}\frac{(-x)^{n-1}h_n(x)}{((n-1)!)^2}.$$
It is clear that $g$ converges in $\C$.
Let $C_2$ be that analytic curve in $\C^2$ defined by $\{y=g(x), x\in \C\}.$
Set $p_{2,n}:=(1/n, g(1/n))\in \Q^2$. Then $p_{2,n}\to 0$ for $n\to 0$ and $p_{2,n}\in C_2$ for $n\geq 2.$
Moreover, $$g(x)=\sum_{2}^n\frac{(-x)^{i-1}h_i(x)}{((n-1)!)^2} \mod (x^{n+1})$$ for $n\geq 2.$
Hence the formal curve induced by $C_2$ at $o$ is defined over $\Q.$
On the other hand, $C_2$ is not contained in any algebraic curve. Otherwise, $g(1)$ should be algebraic.
But $$g(1)=\sum_{n\geq 2}\frac{(-1)^{n-1}h_n(1)}{((n-1)!)^2}=\sum_{n\geq 2}\frac{(n-1)!}{((n-1)!)^2}=\sum_{n\geq 1}\frac{1}{n!}=e-1,$$
which is transcendental.
\end{exe}

In arithmetic dynamics, we often meet some situations similar to Question \ref{quenaive}. So it is interesting to find some natural condition, under which, Question \ref{quenaive} (or its high dimensional/nonarchimedean generalizations) has positive answer.  In \cite[Theorem 1.5]{Xie2015ring}, the author gave an criteria for the algebraicity of adelic curves in $\A^2.$
Such an criteria plays impotent role in the recent breakthrough \cite{FavreGauthier} by Favre and Gauthier, which proved the dynamical Andr\'e-Oort conjecture for curves in the moduli spaces of polynomials. In this paper, we prove a more general  algebraicity criteria.

\medskip

Let $K$ be a finite field extension of $\F$, where $\F$ is $\Q$ or $\bk(t)$ for a field $\bk.$
Let $M_K$ be the set of places of $K$. For $v\in M_K$, denote by $K_v$ the completion of $K$ w.r.t the place $v.$
For each $v\in M_K$, we normalize the absolute value $|\cdot|_v$ to be the unique extension of the usual absolute value $|\cdot|_p$ on $\F$ where $p$ is the restriction of $v$ on $\F.$
More precisely, when $\F=\Q$ and $p$ is a prime, $|\cdot|_p$ is the usual $p$-adic norm with $|p|_p=1/p.$ When $\F=\Q$ and $p=\infty$, $|\cdot|_p$ is the usual archimedean absolute value on $\R$ or  $\C.$ When $\F=\bk(t)$, then $p$ is a closed point in $\P^1(\bk)$ with residue field $\kappa(p)$, then for every $h\in \F$, $|h|_p=e^{-[\kappa(p):\bk]\ord_p(h)}.$

\medskip

Let $X$ be a projective variety over $K.$ 
For $v\in M_K$, denote by $X_v$ the analytification of $X$ w.r.t. $v.$
More precisely, when $v$ is archimedean, $X_v=X(K_v)$ is the associated real or complex variety;
when $v$ is non-archimedean, we use Berkovich's analytification \cite{Berkovich1990}. 
In particular, $X(K_v)\subseteq X_v(K_v)$ is naturally endowed with the $v$-adic topology. 
Using an embedding of $X(K_v)$ to a projective space $\P^N(K_v)$, one may get a distance function $d_v$ on $X(K_v)$. The equivalence class of such distance function does not depend on the choice of embedding and it induces the $v$-adic topology on $X(K_v)$.

For a point $o\in X_v(K_v)$, a \emph{local analytic subspace at $o$} is a closed analytic subspace $\sV$ of some open subset $\sW$ of $X_v$ satisfying $o\in \sV.$ 
Assume that $o\in X(K)\subseteq X_v(K_v)$. 
We say that  \emph{$\sV$ is defined over $K$} if the closed formal subscheme 
$\widehat{\sV_{o}}\subseteq \widehat{(X_{K_v})_{o}}$
is induced by some closed formal subscheme of $\widehat{X_{o}}.$
We say that $\sV$ is \emph{$K$-algberaic at $o$} if there is a subvariety $Y$ of $X$ such that $\dim Y=\dim \sV$ and $Y_v\cap \sV$ contains
a neighborhood of $o$ in $\sV.$
It is clear that if $\sV$ is $K$-algberaic then it  is defined over $\overline{K}$.

\begin{thm}\label{thmgeneralalgebraicity}
Denote by $h: X(K)\to \R$ any Weil height. 
Let $v\in M_K$ and let $\sV$ be a local analytic subspace of $X_v$ at $o$ which is defined over $K$.
Let $x_n, n\geq 0$ be a set of $K$-points in $X(K)$ such that:
\begin{points}
\item $x_n\in \sV$ for $n\geq 0;$
\item $x_n\to o$ in $X(K_v)$ as $n\to \infty;$
\item $\liminf\limits_{n\to \infty} \frac{-\log d_v(o, x_n)}{h(x_n)}>0.$
\end{points}
Set $H:=\cap_{m\geq 0}\overline{\{x_n, n\geq m\}}\subseteq X.$ Then $H_v \cap \sV$ is $K$-algebraic at $o$.
\end{thm}
We indeed proved a stronger statement Theorem \ref{thmgengeneralalgebraicity} which allow us to replace the point $o$ by a closed subscheme $Z$ of $X$ over $K.$ 
\begin{rem}
One may think the assumption (iii) as follows. The assumptions (i) and (ii) only concern the information of $x_n$ at the place $v$.
If one think that the height $h(x_n)$ measures the arithmetic information of $x_n$, then 
the assumption (iii) means that the information of $x_n$ at the place $v$ is not negligible asymptotically. 
\end{rem}

\subsection{A general strategy in the dynamical setting}\label{subsectiongenstra}
Theorem \ref{thmgeneralalgebraicity} provides a general strategy to show the  transcendence of numbers when we have a dynamics. Let us consider the following simple situation.
Assume that $K$ is a number field, $X=\A^2$ and the place $v$ is induced by an embedding $\tau: K\hookrightarrow \C.$
Via $\tau$, we view $K$ as a subfield of $\C.$ Let $f$ be an endomorphism of $X$ defined over $K,$ such that $f(o)=o$ where $o=(0,0).$
Let $h(x)$ be a power series with coefficients in $K$ which convergences when $|x|<r$ for some $r>0.$
Assume that the analytic curve $\sC:=\{y=h(x), |x|<r\}$ is $f$ invariant and $(f|_{\sC})^n\to o$ as $n\to \infty.$
Let $a\in K$ with $|a|<r$ and we want to show the transcendence of $h(a)$. 
We note that if $h(a)$ is algebraic, we may replace $K$ by a suitable finite extension to assume that $(a,h(a))\in K^2.$
Then for every $n\geq 0$, $p_n:=f^n((a,h(a)))\in \sC\cap K^2.$ Moreover $p_n\to o$ as $n\to \infty.$
Hence the assumptions (i) and (ii) of Theorem \ref{thmgeneralalgebraicity} are satisfied.
Then we only need the following two steps:
\begin{points}
\item[(1)] show that the assumption (iii) holds for the sequence $p_n, n\geq 0$;
\item[(2)] show that $\sC$ is not algebraic. 
\end{points}
In our paper, we get two applications following this strategy.  In these applications,
we consider those $f$ having very strong attraction property at $o$ to get (1).
Step (2) concerns the transcendence of  analytic subvarieties, which is in general a difficult problem. However, in many cases, a geometric transcendence problem is easier than proving the transcendence of a number.  In the dynamical setting, such a geometric transcendence problem relates to the understanding of invariant subvarieties of $f$, which is one of the central problems in algebraic dynamics.  In our applications, we get some classification results of invariant subvarieties. These allow us to check the algebraicity of an invariant analytic subvariety. 

\subsection{Invariant subvarieties of products of endomorphisms}
Let $f_i: X_i\to X_i, i=1,\dots,m$ be endomorphisms of projective varieties over an algebraically closed field $\bk$. We want to understand the invariant (or periodic) subvarieties of $\prod_{i=1}^mf_i: \prod_{i=1}^mX_i\to \prod_{i=1}^mX_i.$
We study this problem in two basic cases.


\subsubsection{Products of independent endomorphisms}
An endomorphism $f: X\to X$ of a projective variety is called \emph{amplified} \cite{Krieger2017}, if it is dominant and  there exists a line bundle $L$ on $X$ such that $f^*L\otimes L^{-1}$ is ample. 
Basic examples of amplified endomorphisms are polarized endomorphisms.
\begin{exe}
A \emph{polarized endomorphism} is an endomorphism $f: X\to X$ of a projective variety such that there is an ample line bundle $L$ on $X$ satisfying 
$f^*L=L^{\otimes \la_1(f)}$ for some $\la_1(f)\geq 2.$ We note that $\la_1(f)$ does not depend on $L$, indeed $\la_1(f)^{\dim X}=\deg f.$
A polarized endomorphism is alway amplified. Dominant endomorphisms of $\P^N$ of degree $>1$ are polarized.
\end{exe}
We then introduce a notion of independence for amplified endomorphisms:
For amplified endomorphisms $f: X\to X$ and $g: Y\to Y$, we say that $(X,f)$ and $(Y,g)$ are \emph{independent} and write $(X,f)\perp (Y,g)$ (or $f\perp g$) if 
for every closed irreducible $(f\times g)$-periodic subvariety $Z\subseteq X\times Y$, $Z$ takes form $Z=Z_1\times Z_2$ where 
$Z_1,Z_2$ are closed periodic subvarieties for $f$ and $g$ respectively. 
In Section \ref{subsectionind}, we show that this notion has good behaviors under iterates. Moreover, we prove the following property.
\begin{cor}\label{corprodindind}Let $f_1,\dots, f_m, g_1,\dots, g_n$ be amplified endomorphisms. Assume that $f_i\perp g_j$ for $i=1,\dots, m, j=1,\dots,n.$ Then we have
$$(f_1\times\dots\times f_m) \perp (g_1\times \dots\times g_n).$$
\end{cor}
This property is directly implied by Proposition \ref{proindprod} in Section \ref{subsectionind}.

\begin{rem}
In model theory, their are two notions  ``almost orthogonal" and ``orthogonal" \cite[Definition 2.8]{Medvdev}.  
In \cite[Proposition 2.18]{Medvdev}, Medvedev and Scanlon proved a similar result for orthogonal endomorphisms or more generally for orthogonal $\sigma$-varieties in terminology of model theory. We note that, for two amplified endomorphisms, they are independent in our terminology if and only if they are almost orthogonal. A priori, the notion ``almost orthogonal" is weaker than ``orthogonal", hence \cite[Proposition 2.18]{Medvdev} does not applies in this case directly. We suspect that for amplified endomorphisms these two notions `almost orthogonal" and ``orthogonal" coincide. 
\end{rem}

\begin{rem}
In the definition of the independence, the amplifiedness assumption is necessary for having a result like Corollary \ref{corprodindind}.

For example, let $a_1,a_2\in \bk^*$ which are multiplicatively independent i.e. for $n,m\in \Z$, $a_1^na_2^m=1$ if and only if $m=n=0.$
For example, when $\bk=\C$, we may take $a_1=2, a_2=3.$ Set $a_3:=a_1a_2$. Then $a_1,a_2,a_3$ are pairwise multiplicatively independent.
Let $f_i: \P^1\to \P^1, i=1,2,3$ be the endomorphism $f_i: x_i\mapsto a_ix_i$. It is easy to see that the only periodic curves of $f_i\times f_j, i\neq j$ are 
$\{0\}\times \P^1, \P^1\times \{0\}, \{\infty\}\times \P^1$ and $\P^1\times \{\infty\}.$  But there is an $(f_1\times f_2)\times f_3$-invariant hypersurface 
$\{x_1x_2x_3=1\}$. It is not a product of invariant subvarieties of $f_1\times f_2$ and $f_3$. 
\end{rem}

\subsubsection{Products of separable endomorphisms of curves}
Now we study the case where $X_i$ are smooth projective curves and $f_i$ are separable endomorphisms of degree at least $2.$
There are three types of such endomorphisms: Latt\`es, monomial and nonexceptional (c.f. Section \ref{sectionsependo}).
Endomorphisms of different types are independent to each other. By Corollary \ref{corprodindind}, we only need to study the invariant subvarieties for product of $f_i$ of the same type. 
Basically, endomorphims of Latt\`es or monomial type come from some group structure. The invariant subvarieties of the product of $f_i$ of Latt\`es or monomial comes from the subgroups of the corresponding algebraic group. So in this paper, we focus in the nonexceptional  endomorphisms.  We note that if $f_i$ is nonexceptional, then $X_i\simeq \P^1.$
We get the following classification result.
\begin{pro}\label{prosplitendocurveintro}
Let $f_1,\dots, f_m$ be separable endomorphisms of degree at least $2$ which are nonexceptional.
Let $V$ be an $\prod_{i=1}^mf_i$-invariant subvariety, then there is a partition 
$\{1,\dots,m\}=J_0\sqcup(\sqcup_{j=1}^l J_j),$ fixed points $o_s$ of $f_s$ for $s\in J_0$, $(\prod_{s\in J_j}f_s)$-invariant curves 
$C_j\subseteq \prod_{s\in J_j} X_s,  j=1,\dots,l,$ such that $V=\prod_{j=1}^lC_j$.
\end{pro}
See Proposition \ref{prosplitendocurve} for a more precise form. 
When $\Char\, \bk=0$, this result was obtained in \cite{Medvdev} using model theory and in \cite[Appendix B]{Xie2019} using purely geometric method.
When $\bk=\overline{\Q}$, it was also obtained in \cite[Theorem 1.2]{Ghioca2018c}, as a consequence of their solution of the Dynamical 
Manin-Mumford Conjecture in this case.  Here we follows the method in \cite[Proposition 9.2]{Xie2019}.

\subsection{Transcendence of B\"ottcher coordinates and multiplicative canonical heights}
This application is strongly inspired by the recent preprint \cite{Nguyen} of Nguyen.   

\subsubsection{Transcendence of B\"ottcher coordinates}
Assume that $d\geq 2$ is an integer which is not divided by $\Char\, K.$
Let $$f(z)=a_dz^d+\dots+ a_0\in K[z]$$ be a polynomial of degree $d\geq 2.$ 
 A B\"ottcher coordinate of $f$ is a Laurent series $\phi_f(z)\in \overline{K}((z^{-1}))$
satisfying $$(\phi_f\circ f)(z)=\phi_f(z)^d$$ and
of order $-1$ at $\infty$ i.e.  
it takes form 
$$\phi_f(z)=b_1z+b_0+b_{-1}/z+b_{-1}/z^2\dots.$$
It exists (c.f. Proposition \ref{probottcherexists}) and is unique up to multiply by a $(d-1)$-root of unity (c.f.  Lemma \ref{lemuniquemultibottcher}).
For $v\in M_K$, there is $B_v>0$, such that $\phi_f(z)$ converges in the neighborhood of infinity $\Omega_v(f):=\{x\in \P^1_v|\,\, |z(x)|_v> B_v\}$ and for every $z\in \Omega_v$, $f^{n}(z)\to \infty$ (c.f. Proposition \ref{proconvergencebottcher}). 

\medskip

Our aim is to understand when a product of B\"ottcher coordinates is algebraic. 
Recall that, for a $f\in K[z]$ of degree $d$,  it is either  nonexceptional or of monomial type.
If $f$ is of monomial type, any B\"ottcher coordinate $\phi_f$ of $f$ is algebraic (c.f. Corollary \ref{cormono}). Hence for $a\in K\cap \Omega_v(f)$, $\phi_f(a)$ is algebraic. 
With this fact, we only need to consider the products of B\"ottcher coordinates of nonexceptional polynomials.
\begin{thm}\label{thmtranbottintro}
Assume that $d\geq 2$ is an integer which is not divided by $\Char\, K.$
Let $f_1,\dots f_r$ be nonexceptional  polynomials of degree $d.$ 
Let $a_i, i=1,\dots,r$ be points in $\A^1(K)$ with $(a_i)_v\in \Omega_v(f_i).$
Then the following holds:
\begin{points}
\item For integers $n_1,\dots, n_r$, $\prod_{i=1}^r\phi^{n_i}_{f_i}(a_i)$ is either transcendence over $K$ or a root of unity.
\item The product $\prod_{i=1}^r\phi^{n_i}_{f_i}(a_i)$ is a root of unity, if and only if, for every $j=1,\dots,l,$
$\sum_{s\in J_j}n_sd_{s/j}=0$.
\end{points}
\end{thm}
Here $\{1,\dots,r\}=\sqcup_{j=1}^l J_j$ is a partition. This partition and the positive integers $d_{s/j}, s\in J_j, j=1,\dots,l$ are purely geometric invariants on the pairs $(f_i,a_i), i=1,\dots, r$. In particular, they do not depend on the place $v.$ See Section \ref{subsectionbottcoortran} or Remark \ref{remgeomdata} for their definitions.

When $\Char\, K=0$, Part (i) of Theorem \ref{thmtranbottintro} was proved by Nguyen in \cite[Theorem 1.4]{Nguyen}.
Part (ii) of Theorem \ref{thmtranbottintro} answers the first question proposed by Nguyen in \cite[Section 4.3]{Nguyen}.
The proof of \cite[Theorem 1.4]{Nguyen} relies on a construction of auxiliary polynomials by hand, which is more elementary.

The proof of Theorem \ref{thmtranbottintro} is a typical realization of the strategy presented in Section \ref{subsectiongenstra}. It is by combing a stronger version of  Theorem \ref{thmgeneralalgebraicity} (c.f. Theorem \ref{thmgengeneralalgebraicity}) and the classification of invariant subvarieties of products of nonexceptional polynomials, which is a  special case of Proposition \ref{prosplitendocurveintro}.

\subsubsection{Transcendence of multiplicative canonical heights}
For any polynomial $f\in \overline{\Q}[z]$ and $a\in \A^1(\overline{\Q})$, we denote by $\hat{H}_f(a)$ the multiplicative canonical height (c.f. Section \ref{subsectioncanheight}).
Our aim is to understand when a product of multiplicative canonical heights is algebraic. 
When $f$ is of monomial type, $\hat{H}_f(a)$ is algebraic. When $a$ is $f$-preperiodic, then $\hat{H}_f(a)=1$.
So we only need to consider the product of $\hat{H}_f(a)$ when $f$ is nonexceptional and $a$ is not $f$-preperiodic.

\medskip

Let  $\sD_d$ be the set of nonexceptional polynomials in $\overline{\Q}[z]$ of degree $d.$
Let $\sT_d$ be the set of $(f,a)\in \sD_d\times \overline{\Q}$, for which there is an embedding $\tau: \overline{\Q}\hookrightarrow \C$ such that $|\tau(f^n(a_i))|, n\geq 0$ is unbounded. 
If $(f,a)\in \sT_d$, then $a$ is not $f$-preperiodic. 
\begin{pro}\label{prohproductalgintro}For $(f_i,a_i),\dots, (f_r,a_r)\in \sD_d\times \overline{\Q}$, $n_1,\dots,n_r\in \Z$,
set $T:=\{i=1,\dots, r|\,\, (f_i,a_i)\in \sT_d\}.$
Let $T=\sqcup_{j=1}^l J_j$ be the partition associated to $\sim_w.$
Then $\prod_{i=1}^r\hat{H}_{f_i}(a_i)^{n_i}$ is algebraic if and only if for every $j=1,\dots,l$, $$Q((f_s,a_s),s\in J_j)\in \{\sum_{s\in J_j} n_sz_s=0\}\subseteq \P^{|J_j|}(\Q)\setminus (\cup_{s\in J_j}\{z_s=0\}).$$
\end{pro}
See Section \ref{subsectionalgcanheight} for the definition of the equivalence relation $\sim_w$ and the point $Q((f_s,a_s),s\in J_j)\in \P^{|J_j|}(\Q)\setminus (\cup_{s\in J_j}\{z_s=0\}).$
They are geometric invariants up to Galois conjugations.
 In Corollary \ref{corlinearrelationheight}, we further discuss the case when the product equals to one.

Our Proposition \ref{prohproductalgintro} is strongly motivated by \cite[Corollary 1.6]{Nguyen}, which is the $r=1$ case of our result.
Our Corollary \ref{corlinearrelationheight} is also motivated by the disscusion in \cite[Section 4.3]{Nguyen}.
\begin{rem}Nguyen's result \cite[Corollary 1.6]{Nguyen} (which is the $r=1$ case of  Proposition \ref{prohproductalgintro}) answers a question of Silverman \cite{Silverman2013} for polynomials of degree at least $2.$
\end{rem}

\begin{rem}During the preparation of this paper, we learned from Bell and Nguyen that they got the same result as Theorem \ref{thmtranbottintro} and Proposition \ref{prohproductalgintro} in characteristic $0$ independently. Their proofs are based on a refinement of Nguyen's original proof of \cite[Theorem 1.6]{Nguyen}.
Moreover, their proof work of Theorem \ref{thmtranbottintro}  walks in positive characteristic once we have Proposition 1.8.
\end{rem}

\subsection{Invariant germs curves at infinity}
Assume that $K$ is of characteristic zero. 
In Section \ref{sectioninvargerms}, we study germs curves at infinity which are invariant under polynomial endomorphisms of $\A^N.$
For most of those germs of curves, Theorem \ref{thmgeneralalgebraicity} can be applied in the way we showed in Section \ref{subsectiongenstra}.
We introduce some classes of polynomial endomorphisms. We get more precise results in smaller classes.

Let $\bk$ be a field of characteristic zero.
For $N\geq 2, d\geq 2$, denote by $\sP(N,d)$ the space of dominant endomorphisms $f: \A^N\to \A^N$ taking form 
$$f:(x_1,\dots,x_N)\mapsto (f_1(x_1,\dots, x_N),\dots, f_N(x_1,\dots,x_N))$$ of algebraic degree $\deg_1(f)=d$.
Recall that the algebraic degree of $f$ is $$\deg_1(f):=\max\{\deg(f_1),\dots, \deg(f_N)\}.$$
It is an irreducible quasi-affine variety of dimension $N\binom{d+N+1}{N}$ (c.f. Section \ref{subsectionpolyspace}). 

\medskip

Via the standard embedding $\A^N\hookrightarrow \P^N,$
every $f\in \sP(N,d)$ extends to a rational self-map on $\P^N.$ 
Denote by $\sP^*(N,d)$ the space of $f\in \sP(N,d)$ whose extension is an endomorphism.
This space in a Zariski dense open subset of $\sP(N,d)$ (c.f. Section \ref{subsectionpolyspace}). 

Let $H_{\infty}:=\P^N\setminus \A^N$ be the hyperplane at infinity. For $f\in \sP^*(N,d)$, denote by $\overline{f}: H_{\infty}\to H_{\infty}$ the restriction of $f$ to $H_{\infty}.$
Denote by $\sP^{NS}(N,d)$ the space of $f\in \sP^*(N,d)$ such that for every $n\geq 1$ and fixed point $x$ of $\overline{f}^n$, $d\overline{f}|_x$ is invertible. 

The relations among them are as follows:
$$\sP^{NS}(N,d)\subseteq_{\text{as a dense $G_{\delta}$-set}} \sP^*(N,d)\subseteq_{\text{as a dense open subset}} \sP(N,d).$$ 
Recall that a $G_{\delta}$-set is a countable intersection of open subsets.
When $\bk$ is uncountable, this means that a very general $f\in \sP^*(N,d)$ is contained in $\sP^{NS}(N,d).$
When $\bk$ is countable, the notion ``very general" does make much sense:
even a property is satisfied by a very general point, it may not be satisfied by any $\overline{\bk}-$point.
An alternative is to use the notion of ``adelic general", which was introduced in \cite[Section 3.1]{Xie2019}.
We show, in Section \ref{subsectionadelidenseness}, that an adelic general $f\in \sP^*(N,d)$ is contained in $\sP^{NS}(N,d).$

When $N=2$, we introduce another subspace $\sP^{NP}(2,d):$ 
for $f\in \sP^*(2,d)$, it is contained in $\sP^{NP}(2,d)$ if $\overline{f}$ is not of \emph{polynomial type} i.e. for every $x\in H_{\infty}$, the backward orbits $\cup_{n\geq 0}\overline{f}^{-n}(x)$ is infinite. We have 
$$\sP^{NS}(2,d)\subseteq_{\text{as a dense $G_{\delta}$-set}} \sP^{NP}(2,d)\subseteq_{\text{as a dense open subset}} \sP^*(2,d)\subseteq \sP(2,d).$$

\medskip

Now we may present a transcendence result.
\begin{pro}\label{prospndkextendedattrintro}
Let $f\in \sP^*(N,d)$, $x\in \A^N(K)$ and $v\in M_K$ such that $f^n(x)\to \infty$ as $n\to \infty$ in the $v$-adic topology.
Let $C_v$  be an irreducible $v$-analytic curve at a point $o_v$ at infinity defined over $K$ which is $f$-invariant. 
If $x\in C_v$,
then $C_v$ is algebraic over $K$ at $o_v.$ 
\end{pro}
We also have a result for $f\in \sP(N,d)$ (c.f. Proposition \ref{prospndkattr}). In that result, we need not one but finitely many analytic curves $C_v$ for maybe different places $v$. 

\medskip

One may ask: Are there many such $C_v$ satisfying the assumptions in Proposition \ref{prospndkextendedattrintro}?
This question is answered by Lemma \ref{lemfixpfix}, Lemma \ref{lemformalcurve} and Lemma \ref{lemconvergence}.
Basically, we showed that for every point $o$ in a subset of density $100\%$ in the set of $\overline{f}$-periodic points, there is a unique irreducible formal curve $\widehat{C}$ at $o$ which is $f$-periodic. For every archimedean place $v$, such a formal curve converges to an  irreducible $v$-analytic curve. Hence Proposition \ref{prospndkextendedattrintro} applies for $C_v$. 
\begin{rem}
We suspect that such a statement also holds for all nonarchimedean places. 
However, our proof relies on a version of Hadamard-Perron theorem (c.f. \cite[Theorem 3.1.4]{Abate2001}) on stable/unstable manifolds. But the nonarchimedean version of such a Theorem does not exists in the literature. 
\end{rem}

Once we have such a $C_v$, if we can show that $C_v$ is not algebraic, then by Proposition \ref{prospndkextendedattrintro}, all points in $C_v$ are transcendental. 
Usually it is hard to determine whether a germ of curve $C_v$ is algebraic, but one can determine whether $C_v$ is contained in a curve of certain degree. 
In Corollary \ref{corbounddegree}, we proved that when $f\in \sP^{NS}(N,d)$, all $f$-periodic curves are of degree at most $2$. 
We note that such a degree bound is not true for general $f\in \sP^{*}(N,d)$. 
\begin{exe}
Let $g$ be a polynomial of degree $d$.  
Set $f:=g\times g.$ Then $f\in \sP^*(2,d).$
Let $C_n, n\geq 0$ be the Zariski closure of the graph of $g^n: \A^1\to \A^1$ in $\P^2.$
One may check that $C_n$ is $f$-invariant. But  
$\deg C_n=d^n, n\geq 2$ is unbounded. 
\end{exe}

\subsection{Endomorphisms of $\A^2$}
We now focus in the case $N=2.$
Let $f: \A^2\to \A^2$ be an endomorphism in $\sP^{*}(2,d)$ for some $d\geq 2.$ 

\medskip

We say that $f$ is \emph{homogenous}, if it takes form $f: (x,y)\mapsto (F(x,y), G(x,y))$ where $F,G$ are homogenous polynomials of degree $d.$
We say that $f$ is \emph{homogenous at $o=(o_1,o_2)\in \A^2(\bk)$} if $f$ is homogenous in the new coordinates $x'=x-o_1, y'=y-o_2.$
Easy to see that if $f$ is homogenous at $o$, then such $o$ is unique. 

If $f$ is homogenous at $o$, then for every $\overline{f}$-periodic point 
$x\in H_{\infty}$, the line $L_{o,x}$ passing through $o$ and $x$ is $f$-periodic.
Hence $f$ has infinitely many periodic curves.
We prove that its inverse also holds when $f\in \sP^{NP}(2,d)$.
\begin{thm}\label{thminfinitepercurveintro}For $f\in \sP^{NP}(2,d)$ with $d\geq 2.$ If there are infinitely many $f$-periodic curves, then $f$ is homogenous at some point $o\in \A^2(\bk)$.
Moreover, all but finitely many $f$-periodic curves are lines passing through $o$.
\end{thm}
If $f\not\in \sP^{NP}(2,d)$, such a result does not true (c.f. Example \ref{exenotnp}).
In the proof of Theorem \ref{thminfinitepercurve}, we used the theory of valuative tree introduced by Favre and Jonsson in \cite{Favre2004}.

\subsection*{Acknowledgement}
I would like to thank Xinyi Yuan for his help for the proof of Lemma \ref{lemmaauxiliarysection}.
I thank Charles Favre for telling me Lemma \ref{lemgoodendo}.
I also thank Zhiyu Tian for helpful discussions. 
I thank Thomas Scanlon for explaining to me the notions ``almost orthogonal" and ``orthogonal"  in model theory. 
I thank Jason Bell and Khoa D. Nguyen for kindly telling me their independent work in the same results as Theorem \ref{thmtranbottintro} and Proposition \ref{prohproductalgintro}.
I also thank Khoa D. Nguyen for reminding me that my notion ``independent" relates to the notions ``almost orthogonal" and ``orthogonal" in model theory.

\section{Algebraicity criteria}
The aim of this section is to prove a generalization of Theorem \ref{thmgeneralalgebraicity}, which allow us to replace the point $o$ by a closed subscheme. To state and prove this generalization, we need to introduce some notions and prove a lemma for constructing auxiliary sections. 

\medskip

Let $X$ be a projective variety over $K.$ 
Let $Z$ be a proper closed subscheme of $X$.
Denote by $\widehat{X}_Z$ the formal completion of $X$ along $Z.$

\subsection{Analytic subvarities}
Fix a place $v\in M_K$.
A \emph{local analytic subspace (along $Z$)} is a closed analytic subspace $\sV$ of some open subset $\sW$ of $X_v$ (satisfying $Z_v\subseteq \sV$). 
In particular, every closed subscheme of $X_v$ is a local analytic subspace.

A local analytic subspace $\sV$ along $Z$, is \emph{defined over $K$} if its completion $\widehat{\sV_{Z_v}}\subseteq \widehat{(X_v)_{Z_v}}=\widehat{(X_{K_v})_{Z}}$ along $Z$
is induced by some closed formal  subscheme $\widehat{V_Z}$ of $\widehat{X_{Z}}.$
We say that $\sV$ is \emph{$K$-algberaic along $Z$} if there is a subscheme $Y$ of $X$ such that $\dim Y=\sV$ and $\sV\cap Y_v$ contains a neighborhood of $Z$ in $\sV.$
It is clear that if $\sV$ is $K$-algberaic then it is defined over some finite extension of $K$.

Let $x_n,n\geq 0$ be a sequence of points of $X(K_v).$ Let $B$ be any compact subset of $X_v.$ We write $\lim\limits_{n\to \infty}x_n\subseteq B$ if for every open subset $\sW$ of $X$ containing $B$, $\{n\geq 0|\,\, x_n\not\in \sW\}$ is finite.


\subsection{Green functions for subschemes}
Let $\sV$ be a local analytic subspace of $X_v$.
Let $\sY$ be a closed subscheme of $X_v$ which is contained in $\sV$. We define a function 
$$g_{\sY/\sV,v}: \sV(K_v)\to (-\infty,+\infty]$$ as follows:
Let $\sU_i, i\in I$ be a finite open cover of $\sV$ such that 
for each $i\in I$, there is a finite set $S_i$ of generators of the ideal $I_i$ of $O_{\sV}(\sU_i)$ associated to the closed subspace $\sY\cap \sU_i$ of $\sU_i.$
Because $X_v$ is projective and $\sY$ is a closed subscheme of $X_v$, such cover exists. 
For every $x\in \sV(K_v),$ define 
\begin{equation}
  g_{\sY/\sV,v,i}(x):=
   \begin{cases}
   0&\mbox{if $x\not\in \sU_i(K_v)$}\\
   -\log \max\{0, |h(x)|, h\in S_i\}&\mbox{if $x\in \sU_i(K_v)$}
   \end{cases}
  \end{equation}
for $i\in I$ and 
$$g_{\sY/\sV,v}(x):=\max\{g_{\sY/\sV,v,i}(x), i\in I\}.$$

Observe that, up to a bounded function, $g_{\sY/\sV,v}$ does not depend on the choice of the open cover $\sU_i, i\in I$ and the set generators $S_i, i\in I.$
We still denote by $g_{\sY/\sV,v}$ a function on $\sV(K_v)$ which equals to the above construction up to a bounded function and call it a \emph{Green function} for $\sY/\sV$ and $v.$

\medskip

\begin{facts}\label{facts}
The following facts hold:
\begin{points}
\item The function $g_{\sY/\sV,v}$ is bounded from below.
\item For every $r\in \R$, $g_{\sY/\sV,v}^{-1}((r,+\infty])$ is a neighborhood of $\sY$ in $\sV(K).$
In particular, for $x_n\in \sV(K_v)$, $\lim\limits_{n\to \infty}g_{\sY/\sV,v}(x_n)=+\infty$ if and only if  $\lim\limits_{n\to 0} x_n\subseteq \sY$.
\item For every relatively compact open subset $\sW \subset\subset \sU_i$, $g_{\sY/\sV,v,i}-g_{\sY/\sV,v}$ is bounded on $\sW$.
\item Let $\sV'$ be a local analytic subspace of $X_v$ containing $\sV$. Then $g_{\sY/\sV,v}-g_{\sY/\sV',v}|_{\sV(K_v)}$ is bounded on $\sV$.
\item Let $\sY'$ be a closed subscheme of $\sY$. Then $g_{\sY/\sV,v}-g_{\sY'/\sV,v}$ is bounded from below.
Moreover, if the support of $\sY'$ and $\sY$ are the same, then there is $C>1$ such that $C^{-1}g_{\sY'/\sV,v}-C\leq g_{\sY/\sV,v}\leq Cg_{\sY'/\sV,v}+C.$
\item Let $\sH$ be a closed subsccheme of $X_v$, then $g_{\sY\cap \sH/ \sV\cap \sH,v}-g_{\sY/\sV,v}|_{\sV\cap \sH(K_v)}$ is bounded.
\item When $\sY$ is a point $o\in X_v(K_v)$, then $g_{o/\sV,v}(\cdot)-(-\log d_v(o,\cdot))$ is a bounded function on $\sV(K_v)$.
\item When $\sY$ is a Cartier divisor of $X_v$, $g_{\sY/X_v,v}(\cdot)|_{X(K_v)\setminus \sY}$ equals to a (hence every) Green function of $\sY$ up to a bounded function.
\end{points}
\end{facts}

\subsection{Auxiliary sections}
Let $L$ be an ample line bundle on $X$. We denote by $L_v$ its analytification on $X_v$.
Let $\sV$ be a local analytic subspace along $Z$ which is defined over $K$.
Let $\sI_{Z_v/\sV}$ be the ideal sheaf of $\sO_{\sV_v}$ associated to $Z_v$.
For $m,n\geq 0$, $H^0(\sV, \sI_{Z_v/\sV}^{\otimes m}\otimes L_v^n|_{\sV})$ can be naturally viewed as a space of 
$H^0(\sV, L_v^n|_{\sV}).$ The restriction gives a morphism $\phi_n: H^0(X, L^n)\to H^0(\sV, L_v^n|_{\sV}).$
Because $\sV$ is defined over $K$, $\phi_n$ factor through $\psi_n: H^0(X, L^n)\to H^0(\widehat{V_Z}, L^n|_{\widehat{V_Z}}).$

\begin{lem}\label{lemmaauxiliarysection}
Assume that $\dim \sV<\dim X.$
For every integer $C>0$, there is $N>0$ such that for every $l\geq N$,  there is $s\in H^0(X,L^n)\setminus \{0\}$ such that 
$\phi_n(s)\in H^0(\sV, (\sI_{Z_v/\sV})^{Cn}\otimes L_v^n|_{\sV}).$
\end{lem}

\proof[Proof of Lemma \ref{lemmaauxiliarysection}]
Let $\sI_{Z/X}$ the ideal sheaf of $\sO_{X}$ associated to $Z.$
Denote by $\sI_{Z/\widehat{X_Z}}$ the ideal sheaf of $\sO_{\widehat{X_Z}}$ associated to $Z.$
and $\sI_{Z/\widehat{V_Z}}$ the ideal sheaf of $\sO_{\widehat{V_Z}}$ associated to $Z.$
Denote by $\sI_{\widehat{V_Z}/\widehat{X_Z}}$ the ideal sheaf of $\sO_{\widehat{X_Z}}$ associated to $\widehat{V_Z}.$
Let $\pi: \sO_{\widehat{X_Z}}\to \sO_{\widehat{V_Z}}$ be the quotient morphism, then for $m\geq 0,$ we have
$\sI_{Z/\widehat{X_Z}}^m \subseteq \pi^{-1}(\sI_{Z/\widehat{V}}^m).$
Observe that
$$\sO_{\widehat{V_Z}}/\sI_{Z/\widehat{V}}^m=\sO_{\widehat{X_Z}}/(\sI_{\widehat{V_Z}/\widehat{X_Z}}\pi^{-1}(\sI_{Z/\widehat{V}}^m))$$
$$=(\sO_{\widehat{X_Z}}/\sI_{Z/\widehat{X_Z}}^m)/(\sI_{\widehat{V_Z}/\widehat{X_Z}}\pi^{-1}(\sI_{Z/\widehat{V}}^m)/\sI_{Z/\widehat{X_Z}}^m)$$
$$=(\sO_{X}/\sI_{Z/X}^m)/(\sI_{\widehat{V_Z}/\widehat{X_Z}}\pi^{-1}(\sI_{Z/\widehat{V}}^m)/\sI_{Z/\widehat{X_Z}}^m),$$
which is a quotient of the coherent sheaf $\sO_{X}/\sI_{Z/X}^m$.
Hence $\sO_{\widehat{V_Z}}/\sI_{Z/\widehat{V}}^m$ is a coherent sheaf of $X$. It implies that 
$\sI_{Z/\widehat{V}}^{m-1}/\sI_{Z/\widehat{V}}^m$ is a coherent sheaf of $Z$.

\medskip

For every $n,m\geq 0$ and $s\in H^0(X,L^n)\setminus \{0\}$, $\phi_n(s)\in H^0(\sV, \sI_{Z_v/\sV}^{\otimes m}\otimes L_v^n|_{\sV})$
if and only if $\phi_n(s)$ is contained in the kernal of the morphism 
$$\kappa_{n,m}:H^0(X, L^n)\to H^0(X, L^n\otimes \sO_{\widehat{V_Z}}/\sI_{Z/\widehat{V}}^m).$$
Because $L$ is ample, there is $c_1>0$ such that $\dim H^0(X, L^n)\geq c_1n^{\dim X}$ for $n\geq 0.$
Now we want to bound $\dim H^0(X, L^n\otimes \sO_{\widehat{V_Z}}/\sI_{Z/\widehat{V}}^m)$ from above. 
Observe that 
$$\dim H^0(X, L^n\otimes \sO_{\widehat{V_Z}}/\sI_{Z/\widehat{V}}^m)\leq \sum_{i=1}^m\dim H^0(X, L^n\otimes \sI_{Z/\widehat{V}}^{i-1}/\sI_{Z/\widehat{V}}^i).$$

Consider the projective morphism $\beta: E:=\Proj(\oplus_{i\geq 1}\sI_{Z/\widehat{V}}^{i-1}/\sI_{Z/\widehat{V}}^i)\to Z.$
By GAGA, the analytification of $E\otimes_K K_v$ is exactly the exceptional divisor of 
the blowup of $\sV$ along $Z_v$.  
So $\dim E=\dim \sV-1.$ There is $B\geq 1$ such that for $i\geq B$,
$$\beta_*\sO_E(i-1)=\sI_{Z/\widehat{V}}^{i-1}/\sI_{Z/\widehat{V}}^i.$$
Pick an ample line bundle $M$ on $E$ such that $M\otimes \sO_{E}(-1)$ is effective. 
Then for every $i\geq B$, we have 
$$\dim H^0(X, L^n\otimes \sI_{Z/\widehat{V}}^{i-1}/\sI_{Z/\widehat{V}}^i)=\dim H^0(X, \beta_*(\beta^*(L^n)\otimes \sO_{E}(i-1)))$$
$$=\dim H^0(E, \beta^*(L^n)\otimes \sO_{E}(i-1))\leq \dim H^0(W, \beta^*(L^n)\otimes M^{i-1}).$$
Because $\beta^*L$ is nef and $M$ is ample, by Fujita's vanishing theorem \cite[Theorem 1.4.35]{Lazarsfeld}, there is $B_1\geq B$, such that for every $i\geq B_1$,
$j\geq 1$, we have $H^j(E, \beta^*(L^n)\otimes M^{i-1})=0.$ By Riemann-Roch theorem, there is a polynomial $Q(x,y)$ of degree $\dim E$ such that 
$$\dim H^0(E, \beta^*(L^n)\otimes M^{i-1})=\chi(W, \beta^*(L^n)\otimes M^{i-1})=Q(n,i-1).$$
So there is $c_2>0$ such that $$Q(x,y)\leq c_2(x^{\dim E}+y^{\dim E}+1)$$ for $x,y\geq 0.$
Because $L$ is ample and $\sI_{Z/\widehat{V}}^{i-1}/\sI_{Z/\widehat{V}}^i$ is supported on $Z$,  there is $c_3> 0$ such that 
$$\sum_{i=1}^{B_1}\dim H^0(X, L^n\otimes \sI_{Z/\widehat{V}}^{i-1}/\sI_{Z/\widehat{V}}^i)\leq c_3(n^{\dim Z}+1)$$ for $n\geq 0.$ 
Then for $m\geq B_1$, 
$$\dim H^0(X, L^n\otimes \sO_{\widehat{V_Z}}/\sI_{Z/\widehat{V}}^m)\leq \sum_{i=1}^{B_1}\dim H^0(X, L^n\otimes \sI_{Z/\widehat{V}}^{i-1}/\sI_{Z/\widehat{V}}^i)+\sum_{i= B_1+1}^mQ(n,i-1)$$
$$\leq c_3(n^{\dim Z}+1)+c_2\sum_{i= B_1+1}^m(n^{\dim E}+(i-1)^{\dim E}+1)$$
$$\leq c_4(mn^{\dim E}+m^{\dim E+1}+1)=c_4(mn^{\dim \sV-1}+m^{\dim \sV}+1)$$
for some $c_4>0.$
For $n\geq 1$, take $m=Cn$, we get
$$\dim H^0(X, L^n\otimes \sO_{\widehat{V_Z}}/\sI_{Z/\widehat{V}}^{Cn})\leq c_4(C+C^{\dim \sV})n^{\dim \sV}+c_4\leq c_5n^{\dim \sV}$$
for some $c_5>0.$
Because $\dim \sV< \dim X$, there is $N\geq 1$ such that 
$$\dim H^0(X, L^n)\geq c_1n^{\dim X}\geq c_5n^{\dim \sV}\geq \dim H^0(X, L^n\otimes \sO_{\widehat{V_Z}}/\sI_{Z/\widehat{V}}^{Cn}).$$
Hence $\ker \kappa_{Cn,n}\neq \emptyset$ for $n\geq N$, which concludes the proof. 
\endproof

\subsection{Algebraicity criteria}
The following result generalizes Theorem \ref{thmgeneralalgebraicity}.
\begin{thm}\label{thmgengeneralalgebraicity}
Let $Z$ be a closed subscheme of $X$ over $K.$
Denote by $h: X(K)\to \R$ any Weil height. 
Let $v\in M_K$ and let $\sV$ be a local analytic subspace of $X_v$ along $Z$ defined over $K$.
Let $x_n, n\geq 0$ be a set of $K$-points in $(X\setminus Z)(K)$ such that:
\begin{points}
\item $x_n\in \sV$ for $n\geq 0;$
\item $\lim\limits_{n\to \infty}x_n\subseteq  Z_v$ in $X(K_v)$ as $n\to \infty;$
\item $\liminf_{n\to \infty} \frac{g_{Z_v/\sV, v}(x_n)}{h(x_n)}>0.$
\end{points}
Set $H:=\cap_{m\geq 0}\overline{\{x_n, n\geq m\}}\subseteq X.$ 
Then $H_v \cap \sV$ is $K$-algebraic along $H_v \cap Z$.
In other words, $\dim H=\dim H_v \cap \sV.$
\end{thm}
\rem
Assume that $Z$ is a point $o$, $\dim \sV=1$ and $\widehat{\sV_{o}}$ is irreducible. If $H_v \cap \sV$ is $K$-algebraic at $o$ implies that $\sV$ is $K$-algebraic at $o$.
This is not true in general. The reason is that, it is possible that $\sV$ is not algebraic, but it contains an analytic subspace $\sY$ which is $K$-algebraic along $Z$ and contains all $x_n, n\geq 0$.
\endrem

\proof[Proof of Theorem \ref{thmgengeneralalgebraicity}]
There is $m_0\geq 0$, such that $H=\overline{\{x_n, n\geq m_0\}}.$ After replacing $x_n, n\geq 0$ by $x_{n+m_0}, n\geq 0,$ we may assume that 
$H=\overline{\{x_n, n\geq 0\}}.$ Let $H_1,\dots,H_s$ be all irreducible components of $H$. 
Since $\{x_n, n\geq 0\}\cap Z=\emptyset,$ and $\{x_n, n\geq 0\}\cap H_i$ is Zariski dense in $H_i$, $H_i\not\subseteq Z$ for $i=1,\dots,s.$

\medskip

Write $|Z|$ for the support of $Z.$
There is an effective and ample divisor $D$ of $X$ which contains $|Z|$ and does not contain any $H_i, i=1,\dots,s.$
After replacing $x_n, n\geq 0$ by the subsequence consisting of points not in $D$, we may assume that $x_n\not\in D$ for $n\geq 0.$

\medskip

We may assume that $h$ is a Weil height associated to the divisor $D.$
For $w\in M_K$, let $g_w: X(K_v)\to (-\infty, +\infty]$ be the local heigh functions for $h$ associated to $D$.
In particular, $g_w$ is a Green function on $X_w$ for the divisor $D.$
For every $x\in (X\setminus D)(K)$, we have $$h(x)=\frac{1}{[K:\F]}\sum_{w\in M_K}n_wg_w(x).$$
where $n_w:=[K_w:\F_p]$ where $p$ is the restriction of $w$ on $\F.$
After modifying $g_w$ at finitely many places $w$, we may assume that $g_w(x)\geq 0$ for every $w\in M_K$ and $x\in X(K_w).$
Because $|Z|\subseteq D$, Assumption (ii) implies that $g_v(x_n)\to \infty$. Then $h(x_n)\to \infty$ as $n\to \infty.$
Now we may assume that for every $n\geq 0,$ $h(x_n)\geq \frac{n_v}{[K:\F]}g_v(x_n)>0$.
By Assumption (iii), we may assume that $g_{Z_v/\sV, v}(x_n)\geq ch(x_n)$ for every $n\geq 0.$

\medskip

Let $x^i_n, n\geq 0$ be the subsequence of $x_n, n\geq 0$ consisting of the points contained in $H_i.$
Now we check that the assumptions (i),(ii) and (iii) are satisfied when we replace $x_n, n\geq 0, X, Z,\sV$ by $x^i_n, n\geq 0, H_i, Z\cap H_i,\sV\cap H_i$.
It is clear that Assumptions (i) and (ii) hold.
By (vi) of Facts \ref{facts}, Assumption (iii) also holds.
Now we may assume that $X=H$. In other words, $\{x_n|\,\, n\geq 0\}$ is Zariski dense in $X.$
We only need to show that $\dim \sV=\dim X$. 

Otherwise, assume that $\dim \sV<\dim X.$ We want to get a contradiction.
Pick an integer $C\geq 2(cn_v)^{-1}.$
By Lemma \ref{lemmaauxiliarysection},  there is  $l\geq 1$ and $s\in H^0(X,\sO(lD))\setminus\{0\}$ whose restriction $\phi_l(s)$ on $\sV$
is in $H^0(\sV, (\sI_{Z_v/\sV})^{Cl}\otimes \sO_{X_v}(lD_v)|_{\sV}).$ 
View $s$ as a function in $\sO_X(X\setminus D).$
For every $x\in \sV(K_v)$, there is $C_1>0$ such that 
$$-\log |s(x)|_v\geq Clg_{Z_v/\sV, v}(x)-lg_v(x)-C_1.$$
For every $w\in M_K\setminus \{v\}$ and $x\in (X\setminus D)(K_w)$, there is $C_2>0$ such that 
$$-\log |s(x)|_w\geq -lg_w(x)-C_2$$
Moreover, there is a finite subset $F$ of $M_K\setminus \{v\}$, such that for every $w\in M_K\setminus (F\cup\{v\})$ and $x\in (X\setminus D)(K_w),$
$$-\log |s(x)|_w\geq -lg_w(x).$$
We claim that $s(x_n)=0$ for all but finitely many $n\geq 0.$ 
Assuming this claim, $x_n, n\geq 0$ is not Zariski dense in $X$, which contradicts to our assumption.

We now prove the claim.
Assume that there is an infinite subsequence $x_{n_m}, m\geq 0$ such that $s(x_{n_m})\neq 0.$
By product formula, we have 
$$0=\sum_{w\in M_K}-n_w\log |s(x_{n_m})|_w$$
$$\geq n_vClg_{Z_v/\sV, v}(x_{n_m})-l\sum_{w\in M_K}n_wg_w(x_{n_m})-n_vC_1-|F|[K:\F]C_2$$
$$=n_vClg_{Z_v/\sV, v}(x_{n_m})-lh(x_{n_m})-C_1-|F|[K:\F]C_2$$
$$\geq n_vClch(x_{n_m})-lh(x_{n_m})-C_1-|F|[K:\F]C_2$$
$$=l(Ccn_v-1)h(x_{n_m})-C_1-|F|[K:\F]C_2\geq lh(x_{n_m})-C_1-|F|[K:\F]C_2.$$
Because $h(x_{n_m})\to \infty$ when $m\to \infty$, we get a contradiction. 
\endproof

\subsection{Algebericity of analytic curves}
Applying the same method of the proof of Theorem \ref{thmgengeneralalgebraicity}, one can easily generalizes the algebericity criteria for adelic branches of curves \cite[Theorem 1.5]{Xie2015ring} to any dimension.  Because we want to avoid the definition of adelic branches of curves, we will not do the generalization in this paper.
As an alternative, we prove a result for finite sets of analytic curves. This is sufficient for our applications in the rest of the paper. 

\medskip

Let $H_{\infty}=\P^N\setminus \A^N$ be the hyperplane at infinity.  A \emph{branch of curve at infinity over $K$} is a triple $(v, o, \sC_v)$, where $v\in M_K$, $o\in H_{\infty}(K)$ and $\sC_v$ is an analytic curve in $\P^N_v$ containing $o$, which is defined over $K$ at $o.$ We say $\sC_{v}$ is irreducible at $o$ if its completion $\widehat{(C_v)}_o$ at $o$ is irreducible. 

For $v\in M_K$, denote by $g_v: \A^N(K_v)\to [0,+\infty)$ the naive local height for $H_{\infty}$ i.e. 
$$g_v((x_1,\dots,x_N))=\log\max\{1, |x_i|_v, i=1,\dots, N\},$$
for $(x_1,\dots,x_N)\in K_v^N.$
Let $h: \A^N(K_v)\to [0,+\infty)$ be the naive height i.e. 
$$h(x)=[K:\F]^{-1}\sum_{v\in M_{K}}n_vg_v$$
for $x\in K^N.$

\begin{thm}\label{thmcurvesalge}Let $(v_i,o_i, \sC_{v_i}), i\in I$ be a finite set of branches of curve at infinity over $K$. 
For every $x\in \A^N(K)$, set $I(x):=\{i\in I |\,\, x\in \sC_{v_i}\}.$
Let $x_n, n\geq 0$ be a set of $K$-points in $\A^N(K)$ such that:
\begin{points}
\item $\lim\limits_{n\to \infty}h(x_n)=+\infty$;
\item $\liminf_{n\to \infty} \frac{\max_{i\in I(x_n)}g_{v_i}(x_n)}{h(x_n)}>0.$
\end{points}
Then the Zariski closure $\{x_n, n\geq 0\}$ is of dimension one. 
In particular, for $i\in I$, if $C_{v_i}$ is irreducible and contains infinitely many $x_n, n\geq 0$, then $C_{v_i}$ is algebraic over $K$ at $o$.
\end{thm}

\begin{rem}If $K$ has the Northcott property \cite{Northcott1949, Bombieri2006} e.g. when $\F=\Q$ or a finite field, then the assumption (i) is satisfied when $x_n, n\geq 0$ are pairwise distinct. 
\end{rem}

\proof[Proof of Theorem \ref{thmcurvesalge}]
The last sentence is trivial if we know that  the Zariski closure $\{x_n, n\geq 0\}$ is of dimension one.

We now prove that  the Zariski closure $\{x_n, n\geq 0\}$ is of dimension one by induction on the cardinality $|I|.$
The assumption (ii) shows that $|I|\geq 1.$
We note that we can always remove finitely many $x_n$. 
When $|I|=1$, by (v), (vi)  of Facts \ref{facts}, the assumptions of Theorem \ref{thmgengeneralalgebraicity} are satisfied. 
This concludes the proof. 
Now assume the $|I|=l+1\geq 2$ and Theorem \ref{thmcurvesalge} holds when $1\leq |I|\leq l$. Pick $i_0\in I$.
Consider the subsequence $x_{n_j}, j\geq 0$ of $x_{n}, n\geq 0$ of those $x_n$ with $i_0\in I(x_n)$ and $g_{v_{i_0}}(x_n)= \max_{i\in I(x_n)}g_{v_i}(x_n).$
If $x_{n_j}, j\geq 0$ is finite, then we may remove them and remove $i_0$ from $I$. We conclude the proof by the induction hypothesis.
If $x_{n_j}, j\geq 0$ is infinite, we apply our theorem to the the case $I=\{i_0\}$ and for the sequence $x_{n_j}, j\geq 0$. This shows that the Zariski closure of $x_{n_j}, j\geq 0$
is of dimension $1$.  Let $x_{n}', n\geq 0$ be the subsequence obtained from $x_n, n\geq 0$ by removing those $x_{n_j}, j\geq 0.$
If $x_{n}', n\geq 0$ is finite, then we done. If $x_{n}', n\geq 0$ is infinite,
 we apply the induction hypothesis for the set $I\setminus \{i_0\}$ and $x_{n}', n\geq 0$ to conclude the proof.
\endproof

\section{Invariant subvarieties for products of endomorphisms}
In this section, $\bk$ is an algebraically closed field.

\subsection{Semi-conjugacies by dominant finite morphisms}
Recall that every dominant endomorphism on a projective variety is finite \cite[Lemma 5.6]{fa}. 

\medskip

Let $f: X\to X,g: Y\to Y$ be dominant endomorphisms of projective varieties. Write $f\gtrsim g$ if there is a dominant finite morphism $\pi: X\to Y$ such that $\pi\circ f=g\circ \pi.$
We also write $f\gtrsim_{\pi} g$ when we want to emphasize  $\pi$.  If $f\gtrsim g$ and $g\gtrsim h$, then $f \gtrsim h.$

\medskip

Easy to see that the following statements are equivalent:
\begin{points}
\item There is a dominant endomorphism $h: Z\to Z$ of a projective variety such that $h\gtrsim f^n$ and $h\gtrsim g^n$ for some $n\geq 1.$
\item There is a closed irreducible $(f\times g)$-periodic subvariety in $X\times Y$ whose projection to each factor is finite and dominant. 
\end{points}
Write $f\sim g$ and say $f$ \emph{semi-equivalent} to $g$ if the above equivalent conditions holds. Note that $f\sim g$ if and only if $f^n\sim g^n$ for some $n\geq 1.$

\rem\label{remsemieqdeg} If $f\sim g$, then $\dim X=\dim Y$, $\deg f=\deg g$ and $\deg_{\rm sep} f=\deg_{\rm sep} g.$
\endrem
\begin{lem}\label{lemrelation}The relation $\sim$ is an equivalence relation.
\end{lem}
\proof[Proof of Lemma \ref{lemrelation}]
We only need to show the transitivity.
Assume that $f_1\sim f_2$ and $f_2\sim f_3$, we want to show that $f_1\sim f_3.$
After replacing $f_1,f_2,f_3$ by some common positive iterate, we may assume that there are dominant endomorphisms $g_1, g_2$ of a projective varieties such that 
$g_1\gtrsim_{\pi_1}f_1$, $g_1\gtrsim_{\pi_2}f_2$ and $g_2\gtrsim_{\phi_1}f_2$, $g_2\gtrsim_{\phi_2}f_3.$ 
Then $g_1\times g_2\gtrsim_{\pi_1\times \phi_1}f_1\times f_2.$   Let $\Gamma_{g_1}$ be the image of $\pi_1\times \pi_2$, which is $(f_1\times f_2)$-invariant.
Let $W$ be an irreducible component of $(\pi_1\times \phi_1)^{-1}(\Gamma_{g_1})$. Then $W$ is $(g_1\times g_2)$-periodic.
Then there is a dominant endomorphism $h: Z\to Z$ of a projective variety such that $h\gtrsim g_1^n, h\gtrsim g_2^n$ for some $n\geq 1.$
Hence $h\gtrsim f_1^n, h\gtrsim f_3^n$, which concludes the proof.
\endproof

\subsection{Amplified endomorphisms}
An endomorphism $f: X\to X$ of a projective variety is \emph{amplified} \cite{Krieger2017}, if it is dominant and  there exists a line bundle $L$ on $X$ such that $f^*L\otimes L^{-1}$ is ample. 

%

We have the following properties.
\begin{points}
\item The identity endomorphism $\id: X\to X$ is amplified if and only if $X$ is a point.
\item  A dominant endomorphism $f$ on a curve is amplified if and only if $\deg f\geq 2.$
\item Let $n$ be a positive integer. Then $f$ is amplified if and only if $f^n$ is amplified \cite[Lemma 5.1]{Xie2019}.
\item If $f$ is amplified and $V$ is an $f$-invariant closed subvariety of $X$, then $f|_V$ is amplified.
\item Let $f: X\to X, g: Y\to Y$ be endomorphisms of projective varieties.  Assume that $f\sim g$, then $f$ is amplified if and only if $g$ is amplified. 
\item Let $f: X\to X, g: Y\to Y$ be amplified endomorphisms of projective varieties. Then $(X,f)\times (Y,g):=(X\times Y, f\times g)$ is amplified. 
\end{points}

\medskip

For an endomorphism $f$, denote by $\Fix(f)$ its set of fixed points and $\Per(f)$ its set of periodic points.
The proof of  \cite[Theorem 5.1]{fa} shows the following important property of  amplified endomorphisms.
\begin{pro}\label{proamplifiedfix}
If $f$ is amplified,  then $\Per(f)$ is Zariski dense and for all $n\geq 1$, $\Fix(f^n)$ is finite.
\end{pro}

Denote by $k(X)^f$ the field of $f$-invariant rational functions on $X.$
\begin{pro}\label{pronoinvfun}If $f$ is amplified, then $\bk(X)^f=\bk.$
\end{pro}
\proof[Proof of Proposition \ref{pronoinvfun}] Assume that there is $\phi\in \bk(X)^f\setminus \bk$.
We get a rational map $\phi: X\dashrightarrow \P^1.$ Let $Y\subseteq X\times \P^1$ be the graph of $\phi$ and $\pi_X: Y\to X$, $\psi: Y\to \P^1$ be the projections. Then $\pi$ is birational and $g:=\pi^{-1}\circ f\circ\pi$ is an endomorphism on $Y$ and $\psi\circ g=\psi.$  For every $c\in \P^1$, set $Y_c:=\psi^{-1}(c)$ and $X_c:=\pi_X(Y_c)$.  We note that 
$X_c$ is $f$-invariant. Let $\eta$ be the geometric generic point of $\P^1$. Then $f$ induces an endomorphism $f_{\eta}$ on $X_{\eta}$ which is also amplified. By Proposition \ref{proamplifiedfix}, $\Per(f_\eta)$ is Zariski dense in $X_{\eta}.$ After replacing $f$ by a suitable iterates, there is an $f$-fixed point $o\in X_{\eta}$ such that $\pi_X^{-1}$ is well defined at $o$. Since $\psi(\pi_X^{-1}(o))=\eta$, the Zariski closure $C$ of $o$ in $X$ is a curve. Because $o$ is $f_{\eta}$-fixed, $f|_C=\id.$ So $f|_C$ is not amplified which contradicts to Property (iii) above.
\endproof

\subsection{Independence}\label{subsectionind}

\begin{defi}For amplified endomorphisms $f: X\to X$ and $g: Y\to Y$, we say $(X,f)$ and $(Y,g)$ are \emph{independent} and write $(X,f)\perp (Y,g)$ (or $f\perp g$) if 
for every closed irreducible $(f\times g)$-periodic subvariety $Z\subseteq X\times Y$, $Z$ takes form $Z=Z_1\times Z_2$ where 
$Z_1,Z_2$ are closed periodic subvarieties for $f$ and $g$ respectively. 
\end{defi}

\begin{rem}
We may replace periodic subvarieties by preperiodic subvarieties to get the same definition. 
\end{rem}

It is clear that $f\perp g$ if and only if there is some $n\geq 1$ such that $f^n\perp g^n.$
\begin{rem}\label{remsemiind}
Let $f: X\to X$ and $g: Y\to Y$ be amplified endomorphisms. Let $\pi: X\to Y$ be a morphism satisfying $g\circ \pi=\pi\circ f.$
Let $h: Z\to Z$ be an amplified endomorphism. One may check the following properties.
\begin{points}
\item If $g\perp h$ and $\pi$ is finite, then $f\perp h$.
\item If $f\perp h$ and $\pi$ is surjective, then $g\perp h$.
\end{points}
\end{rem}

\begin{pro}\label{proindprod}Let $f: X\to X$, $g: Y\to Y$ and $h: Z\to Z$ be amplified endomorphisms. Assume that $f\perp h$ and $g\perp h$, then $(f\times g)\perp h$.
\end{pro}

We need the following two lemmas. 

\begin{lem}\label{lemprodirr}Let $X,Y$ be two normal projective varieties of dimensions $d_X,d_Y$ respectively. Let $W$ be a closed irreducible subvariety of $X\times Y$ of dimension $d_W$.
Let $\pi_X,\pi_Y$ be the projections from $X\times Y$ to $X$ and $Y$ respectively. Assume that $\pi_X(W)=X$ and $\pi_Y(W)=Y.$
Then there is a numerical class $\alpha\in N_{d_W-d_Y}(X)$ such that for a general point $y\in Y(\bk)$, $\dim(\pi_Y^{-1}(y)\cap W)=d_W-d_Y$ and for every irreducible component $R$ of $\pi_Y^{-1}(y)\cap W$, $[R]=(\pi_X|_{\pi_Y^{-1}(y)})^*\alpha\in N_{d_W-d_Y}(\pi_Y^{-1}(y)).$
\end{lem}

This lemma is a generalization of \cite[Lemma 9.3]{Xie2019} with a similar proof.
\proof[Proof of Lemma \ref{lemprodirr}]
Let $\phi: B\to Y$ be the normalization of $Y$ in the field $\bk(W)$.
Consider $\psi:=\id\times_{Y}\phi:(X\times Y)\times_{Y}B=X\times B\to X\times Y$.
Denote by $\pi_B: X\times B\to B$ the projection to $B.$
The Stein factorization of $\pi_Y|_W: W\to Y$ induces a closed immersion $\iota:W\hookrightarrow X\times B$.
We denote by $V$ its image.  In particular, the composition $\psi\circ\iota$ is exactly the inclusion $W\hookrightarrow X\times Y.$
We get the following commutative diagram. 

\[
 \xymatrix{
    V\ar@{^{(}->}[d]\ar[r]^{\simeq}& W\ar@{^{(}->}[d] \\
    X\times B \ar[r]^{\psi}\ar[d]^{\pi_B} & X\times Y \ar[d]^{\pi_Y} \\
     B\ar[r]^{\phi}& Y
    }
\]
There is a nonempty open subset $U\subseteq Y$ such that $\pi_Y|_W$ is flat above $U$ and for every $b\in \phi^{-1}(U)$, the fiber $V_b$ of $V$ above $b$ is reduced and irreducible. 
After shrinking $U$, we may assume that $\phi|_{\phi^{-1}(U)}: \phi^{-1}(U)\to U$ is finite.
Then $\psi(V\cap (X\times \phi^{-1}(U)))=W\cap \pi_Y^{-1}(U)$.
Moreover, 
for every $b\in \phi^{-1}(U)$, $\psi|_{X\times b}: X\times b\to \pi_Y^{-1}(\phi(b))$ is an isomorphism and sends $V_b$ to an irreducible component of 
$(\pi_Y|_W)^{-1}(\phi(b)).$ For every $b\in \phi^{-1}(U)$, the numerical class of $\pi_{X}(\psi(V_b))$ in $N_{d_W-d_Y}(X)$ is the same class $\alpha$. 
For every $y\in U(\bk)$ and every irreducible component $R$ of $(\pi_Y|_W)^{-1}(y)$, there is $b\in \phi^{-1}(y)$, such that $R=\psi(V_b)$. So we get 
$[R]=(\psi|_{X\times b})_*([V_b])=(\psi|_{X\times b})_*(\pi_{X}\circ \psi|_{X\times b})^*\alpha=(\pi_X|_{\pi_Y^{-1}(y)})^*\alpha.$
\endproof

\begin{lem}\label{lemcarre}Let $X,Y$ be two normal projective varieties of dimensions $d_X,d_Y$ respectively. Let $W$ be a closed subvariety of $X\times Y$ of pure dimension $d_W$.Let $\pi_X,\pi_Y$ be the projections from $X\times Y$ to $X$ and $Y$. Let $\beta_X$
be an ample numerical divisor class on $X$. Then the following two statements are equivalent:
\begin{points}
\item $W=V\times Y$ for some closed subvariety $V$ of $X$;
\item 
$W\cdot \pi_Y^*\beta_X^{d_W-d_Y+1}= 0$.
\end{points}

In particular, if $W$ numerically equivalent to $c V\times Y$ for some closed subvariety $V$ of $X$ and some $c>0$, then 
$W=V'\times Y$ for some closed subvariety $V'$ of $X$.
\end{lem}

\proof[Proof of Lemma \ref{lemcarre}]
It is clear that (i) implies (ii). Now we prove that (ii) implies (i).
After replacing $\beta_X$ by a suitable multiple, we may assume that it is represented by a very ample divisor $H$ on $X$.
Let $H_1,\dots, H_{d_W-d_Y+1}$ be general elements in the linear system $|H|$. 
Then $W$ intersects $\pi_X^*(H_1)\cdots \pi_X^*(H_{d_W-d_Y+1})$ properly.
Because $W\cdot \pi_Y^*\beta_X^{d_W-d_Y+1}= 0$, $W\cap (\pi_X^{-1}(H_1)\cap \dots \cap\pi_X^{-1}(H_{d_W-d_Y+1}))=\emptyset.$
Hence $\pi_X(W)\cap H_1\cap \dots \cap H_{d_W-d_Y+1}=\emptyset.$ Then $\dim \pi_X(W)\leq d_W-d_Y.$
For every point $x\in \pi_X(W)$, $\dim (\pi_X^{-1}(x)\cap W)=d_W-\dim \pi_X(W)\geq d_Y.$ Hence $\dim (\pi_X^{-1}(x)\cap W)=d_Y,$ which implies that 
$\pi_X^{-1}(x)\cap W=x\times Y.$ It follows that $W=\pi_X(W)\times Y$, which concludes the proof.
\endproof

\proof[Proof of Proposition \ref{proindprod}]
We do the proof by induction on $\dim X+\dim Y+\dim Z$.
It is clear that Proposition \ref{proindprod} holds when one of $\dim X,\dim Y,\dim Z$ is $0.$

Let $W$ be an irreducible closed $(f\times g\times h)$-periodic subvariety of $X\times Y\times Z.$ After replacing $f,g,h$ by $f^m,g^m,h^m$ for some $m\geq 1$, we may assume that $W$ is $(f\times g\times h)$-invariant. We want to show that $W=W_1\times W_2$ where $W_1, W_2$ are closed invariant  varieties for $f\times g$ and $h$ respectively.

Let $\pi_X,\pi_Y,\pi_Z$ be the projection from $X\times Y\times Z$ to $X,Y, Z$ respectively.
If $\pi_Y(W)\neq Y$, then $\pi_Y(W)$ is $g$-invariant. By Remark \ref{remsemiind}, $g|_{\pi_Y(W)}\perp h.$
Because $\dim \pi_Y(W)<\dim Y$, we conclude the proof by the induction hypothesis.
By similar argument, we may assume that $\pi_X(W)=X, \pi_Y(W)=Y$ and $\pi_Z(W)=Z$.

Let $\pi_{X\times Z}:X\times Y\times Z\to X\times Z$ be the projection.
By Lemma \ref{lemprodirr}, there is a numerical class $\alpha\in N_{d_W-d_Y}(X\times Z)$ and a Zariski dense open subset $U$ of $Y$ such that for every  $y\in U(\bk)$, $\dim(\pi_Y^{-1}(y)\cap W)=d_W-d_Y$ and for every irreducible component $R$ of $\pi_Y^{-1}(y)\cap W$, $$[R]=(\pi_X|_{\pi_Y^{-1}(y)})^*\alpha\in N_{d_W-d_Y}(\pi_Y^{-1}(y)).$$
In particular $[\pi_Y^{-1}(y)\cap W]=c(\pi_{X\times Z}|_{\pi_Y^{-1}(y)})^*\alpha \in N_{d_W-d_Y}(\pi_Y^{-1}(y))$ for some constant $c>0.$

By Proposition \ref{proamplifiedfix}, after replacing $f,g,h$ by $f^m,g^m,h^m$ for some $m\geq 1$, we may assume that there is a $g$-fixed point $b\in U$ which is smooth.
Every irreducible component $T$ of $\pi_Y^{-1}(b)\cap W$ is $f\times g\times h$-periodic, hence takes form $T_1\times T_2$ where $T_1$, $T_2$ are closed periodic subvarieties of $f$ and $h$ respectively.  

First assume that $T_2=Z$.  By Lemma \ref{lemcarre}, $\pi_Y^{-1}(b)\cap W= T_b\times Z$ for some closed subvariety $T_b$ of $X$.
By Lemma \ref{lemcarre} again, we get that for every $y\in U$, $\pi_Y^{-1}(y)\cap W= T_y\times Z$ for some closed subvariety $T_y$ of $X$.
Let $\pi_{X\times Y}:X\times Y\times Z\to X\times Y$ be the projection. We have that, for every $t\in \pi_{X\times Y}(W)$, $\pi_{X\times Y}^{-1}(t)=t\times Z,$ hence $W=\pi_{X\times Y}(W)\times Z$ which concludes the proof.

Now assume that $T_2\neq X.$ Pick a general point $z\in Z(\bk)$, we have $\pi_Z^*(z)\cap T=\emptyset$, hence $\pi_Z^*(z)\cdot (\pi_{X\times Z}|_{\pi_Y^{-1}(b)})^*\alpha=0.$
Then $\pi_Z^*(z)\cdot (\pi_Y^{-1}(b)\cap W)=0.$ Because $z$ is general, the intersection of $\pi_Z^*(z)$ and $(\pi_Y^{-1}(b)\cap W)$ is proper.  So $\pi_Z^{-1}(z)\cap (\pi_Y^{-1}(b)\cap W)=\emptyset.$ Let $\pi_{Y\times Z}:X\times Y\times Z\to Y\times Z$ be the projection.
Then $(b,z)\not\in \pi_{Y\times Z}(W).$ Since $\pi_{Y\times Z}(W)$ is $g\times h$-invariant, it takes form $Y_1\times Z_1$ where $Y_1,Z_1$ are invariant under $g, h$ respectively. Moreover, we have $\dim Y_1+\dim Z_1\leq \dim Y+\dim Z-1.$
We have $W\subseteq X\times Y_1\times Z_1.$ By Remark \ref{remsemiind}, $f\perp h|_{Z_1}$ and $g|_{Y_1}\perp h|_{Z_1}.$ We conclude the proof by the induction hypothesis. 
\endproof

\medskip

\begin{cor}\label{corprodindZDO}Let $f_i: X_i\to X_i, i=1,\dots, m$ be amplified endomorphisms. Assume that $f_i\perp f_j$ for $i\neq j$. 
For every point $p:=(x_1,\dots, x_m)\in X_1\times\dots\times X_m(\bk)$, its orbit under $f:=f_1\times\dots\times f_m$ is Zariski dense in $X_1\times\dots\times X_m$ if and only if 
for every $i=1,\dots,m,$ the orbit of $x_i$ under $f_i$ is Zariski dense in $X_i.$
\end{cor}

\begin{rem}
Under the assumption of Corollary \ref{corprodindZDO}.
Proposition \ref{pronoinvfun} and Corollary \ref{corprodindZDO} implies that the Zariski dense orbit conjecture \cite[Conjecture 5.10]{Medvdevv1} and its adelic version \cite[Conjecture 1.10]{Xie2019} for $f_1\times\dots\times f_m$ can be reduced to the same conjectures for each $f_i.$
\end{rem}

\proof[Proof of Corollary \ref{corprodindZDO}]
If the $f$-orbit of $p$ is Zariski dense in $X_1\times\dots\times X_m$, then for every $i=1,\dots,m,$ the orbit of $x_i$ under $f_i$ is Zariski dense in $X_i.$

Now assume that for every $i=1,\dots,m,$ the orbit of $x_i$ under $f_i$ is Zariski dense in $X_i.$
Let $Z$ be the Zariski closure of the orbit of $p.$ 
After replacing $f$ by a suitable iterate and $p$ by $f^l(p)$, we may assume that $Z$ is irreducible and $f$-invariant.
Then we have $Z=Z_1\times\dots\times Z_m$ for some closed $f_i$-invariant subvarieties $Z_i$ of $X_i.$
Our assumptions shows that $Z_i=X_i$ for $i=1,\dots,m$, which concludes the proof.
\endproof

\subsubsection{Examples of independent endomorphisms}

\begin{exe}\label{execurve}
Let $f_i: C_i\to C_i, i=1,2$ be endomorphisms of projective curves of degree at least $2$. Then $f_1\not\perp f_2$ if and only if $f_1\sim f_2$
\end{exe}

\begin{pro}\label{propolarizedendomorphism}
Let $f_i: X_i\to X_i, i=1,2$ be two polarized endomorphisms. If $\la_1(f_1)\neq \la_1(f_2)$, then $f_1\perp f_2$.
\end{pro}

\proof[Proof of Proposition \ref{propolarizedendomorphism}]
Let $V$ be an irreducible closed $(f_1\times f_2)$-periodic subvariety of dimension $d$.
After replacing $f_1,f_2$ by a suitable common positive iterate, we may assume that  $V$ is $(f_1\times f_2)$-invariant.
Let $\pi_i: X_1\times X_2\to X_i$ be the projection to the $i$-th coordinate. 
Set $V_i:=\pi_i(V)$ and $d_i:=\dim V_i$. We have $d_1+d_2\geq d$ and the equality holds if and only if $V=V_1\times V_2.$
So we may assume that $d_1+d_2> d.$

Let $L_i$ be an ample line on $X_i$ such that $f_i^*L_i=L_i^{\otimes \la_1(f_i)}.$
Denote by  $\alpha_i$ the numerical class associated to $\pi_i^*L_i.$
Observe that for $j\in \{0,\dots,d\},$ $\alpha_1^j\cdot \alpha_2^{d-j}\cdot V\geq 0$ and it $>0$ if $j= d_1$ or $d-j= d_2.$

Set $f:=f_1\times f_2$. For every real numbers $u_1, u_2\in \R$, we have 
$$\deg(f|_V)(V\cdot (u_1\alpha_1+u_2\alpha_2)^d)=(f_*V\cdot (u_1\alpha_1+u_2\alpha_2)^d)$$
$$=(V\cdot (u_1f^*\alpha_1+u_2f^*\alpha_2)^d)=(V\cdot (u_1\la_1(f_1)\alpha_1+u_2\la_1(f_2)\alpha_2)^d)$$
Compare the coefficients of the two different terms $u_1^{d_1}u_2^{d-d_1}$ and $u_1^{d-d_2}u_2^{d_2}$, we get 
$$\deg(f|_V)=\la_1(f_1)^{d_1}\la_1(f_2)^{d-d_1}=\la_1(f_1)^{d-d_2}\la_1(f_2)^{d_2}.$$
It implies that $(\la_1(f_1)/\la_1(f_2))^{d_1+d_2-d}=1$, hence $\la_1(f_1)=\la_1(f_2)$. This contradicts our assumption.
\endproof

\begin{pro}\label{procurvesurface}
Let $f:X\to X$ be an amplified automorphism of a projective surface and $g:C\to C$ be an endomorphism of a projective curve of degree at least $2$.
Then $f\perp g.$
\end{pro}
\proof[Proof of Proposition \ref{procurvesurface}]
Let $\pi_X: X\times C\to X, \pi_C:X\times C\to C$ be the projections. 
Let $V$ be a closed irreducible $(f\times g)$-periodic subvariety. After replacing $f,g$ by a suitable common positive iterate, we may assume that 
$V$ is $(f\times g)$-invariant. 
We may assume that $\dim V\leq 2$, $\dim \pi_X(V)\geq 1$ and $\dim \pi_C(V)= 1.$
Because any automorphism of a curve is not amplified, there is no $f$-periodic curve. 
So $\dim \pi_X(V)=2.$ 
After replacing $f,g$ by a suitable common positive iterate, we may assume that there is a fixed point $o\in C(\bk)$.
Then $\pi_X(\pi_C^{-1}(o))$ is an $f$-invariant curve, which is a contradiction. 
\endproof

\subsection{Products of separable endomorphisms of curves}\label{sectionsependo}
Let $C$ be a smooth projective curve.
Let $g: C\to C$ be a separable endomorphism with $\deg g\geq 2$. Then $C$ is either $\P^1$ or an elliptic curve.

\medskip

We say that $g$ is \emph{of Latt\`es type}, if it semi-equivalents to an endomorphism of an elliptic curve i.e.  there exists an endomorphism of an elliptic curve $h: E\to E$ and a finite morphism
$\pi: E\to C$ such that $f\circ\pi=\pi\circ h.$

We say that $g$ is \emph{of monomial type}, if it semi-equivalents to an endomorphism of a monomial map i.e.  there exists a monomial endomorphism $h: \P^1\to \P^1$ taking form $x\mapsto x^{\pm d}, d\geq 2$ and a finite morphism
$\pi: \P^1\to C$ such that $f\circ\pi=\pi\circ h.$ We note that in this case $C\simeq \P^1.$

We say that $g$ is \emph{exceptional} if it is of  Latt\'es type or monomial type.
Otherwise, it is said to be \emph{nonexceptional}.

\medskip

For every separable endomorphism $g: \P^1\to \P^1$ of $\deg g\geq 2$, it has
exactly one type in Latt\'es, monomial, and nonexceptional.  Moreover,
the types of $g^n,n\geq 1$ are the same.

The following facts are well known.
\begin{points}
\item If two endomorphisms of curves are semi-equivalent, then they have the same type.
In particular, for endomorphisms $f,g$ of curves, if they have different types, then $f\perp g.$
\item If there is a nonzero rational differential form $\omega$, such that $g^*\omega=\mu \omega$ for some $\mu\in \bk^*$, then $g$ is exceptional.
\item If $f: C\to C$ is nonexceptional and $C$, then $C=\P^1.$
\end{points}

\subsubsection{Nonexceptional endomorphisms}
When $\Char\, \bk=0$, the following result on the invariant subvarieties was obtained in \cite{Medvdev} using model theory and in \cite[Proposition 9.2]{Xie2019} using purely geometric method.
When $\bk=\overline{\Q}$, it was also obtained in \cite[Theorem 1.2]{Ghioca2018c}, as a consequence of their solution of the Dynamical 
Manin-Mumford Conjecture in this case. 
The proof of \cite[Proposition 9.2]{Xie2019} can be easily generalized in positive characteristic when all $f_i$ are separable. Hence we get the following result. 
\begin{pro}\label{proinvsubvspl}Assume that $N\geq 2$, $\deg f_i\geq 2, i=1,\dots,N$ and all $f_i,i=1,\dots,N$ are separable and nonexceptional. Let $V$ be a proper irreducible closed subvariety of $(\P^1)^N$ which is invariant under  $f$.
Then there exists $1\leq i< j\leq  N$ such that $V\subseteq \pi_{i,j}^{-1}(C)$ where $\pi_{i,j}: (\P^1)^N\to (\P^1)^2$ is the projection to the $i,j$-th coordinates and $C$ is an $(f_i\times f_j)$-invariant curve in $(\P^1)^2.$ Moreover, the normalization $\widetilde{C}$ of $C$ is $\P^1$ and the endomorphism on $\widetilde{C}$ induced by $(f_i\times f_j)|_C$ is nonexceptional.
\end{pro}

Easy to see that in Proposition \ref{proinvsubvspl}, we may replace the invariant subvarieties by the periodic or preperiodic subvarieties.

\subsubsection{Invariant subvarieties}
By Remark \ref{remsemieqdeg}, for semi-equivalent endomorphisms $f,g$ of projective curves, they have the same degree; they are separable at the same times; and they are nonexceptional at the same times.

\medskip

Let $\sS$ be the set of equivalent classes under $\sim$ of all separable endomorphisms of smooth projective curve of degree at least $2.$
Let $\sS'$ be the subset of $\sS$ of those equivalent classes of nonexceptional endomorphisms. 
The following result was obtained in \cite{Medvdev} in characteristic zero for polynomial endomorphisms.
\begin{pro}\label{prosplitendocurve}Let $s_1,\dots, s_m\in \sS.$
For $i=1,\dots, m$, let $I_i$ be a finite subset of $s_i.$ Then we have the following statement:
\begin{points}
\item Every irreducible closed $\prod_{i=1}^m\prod_{j\in I_i}f_j$-invariant subvariety takes form 
$\prod_{i=1}^mV_i$ where $V_i$ is a $\prod_{j\in I_i}f_j$-invariant subvariety.
\item Assume that $s_i\in \sS',$ and $V$ is a $\prod_{j\in I_i}f_j$-invariant subvariety. Then there is a partition 
$I_i=J_0\sqcup(\sqcup_{j=1}^l J_j),$ fixed points $o_s$ of $f_s$ for $s\in J_0$, $(\prod_{s\in J_j}f_s), j=1,\dots,l$-invariant curves 
$C_j\subseteq (\P^1)^{|J_j|},$ such that $V=\prod_{j=1}^lC_j$. Here we identified $(\P^1)^{|I_i|}$ with $(\prod_{j=1}^l\P^1)^{|J_j|}.$
\end{points}
\end{pro}

In the second part of Proposition \ref{prosplitendocurve}, $C_j$ may be singular.  
One may reformulate it as follows to avoid the singularities.
\begin{pro}[Reformulate of (ii) of Proposition \ref{prosplitendocurve}]\label{proreformulationii}
Assume that $s_i\in \sS',$ and $V$ is a $\prod_{j\in I_i}f_j$-invariant subvariety. Then there is a partition 
$I_i=J_0\sqcup(\sqcup_{j=1}^l J_j),$ fixed points $o_s$ of $f_s$ for $s\in J_0$,
endomorphisms $g_j: \P^1\to \P^1, j=1,\dots, l$ with $g_j\gtrsim_{\pi_{s/j}} f_s$ for every $s\in J_j$, such that 
for every $j=1,\dots,l$, $$\Phi_j:=(\pi_{s/j})_{s\in J_j}: \P^1\to (\P^1)^{|J_j|}$$ is birational to its image and $V=(\prod_{s\in J_0}o_s)\times (\prod_{j=1}^l\Phi_j(\P^1)).$
Here we identified $(\P^1)^{|I_i|}$ with $(\prod_{j=1}^l\P^1)^{|J_j|}.$
The partition $I_i=\sqcup_{j=0}^l J_j$, the points $o_s, s\in I_0$ are unique and the $(g_j; \pi_{s/j}, s\in J_j)$ is unique up to changing the coordinate on $\P^1$ i.e. changing $(g_j; \pi_{s/j}, s\in J_j)$ to 
$(h^{-1}\circ g_j\circ h; \pi_{s/j}\circ h, s\in J_j)$ where $h$ is an automorphism of $\P^1.$
\end{pro}


\begin{rem}\label{remjzeroempty}In (ii) of Proposition \ref{prosplitendocurve},  if the projection of $V$ on each factor is dominant, then $J_0=\emptyset.$
\end{rem}
\begin{rem}\label{reminvperpreper}It is easy to check that, in Proposition \ref{prosplitendocurve}, we may replace invariant subvarieties by periodic or preperiodic subvarieties.
\end{rem}
\begin{rem}\label{rempoly}
For $s\in \sS'$, if there is one $f\in s$ is a polynomial endomorphism, then every element in $s$ conjugates to a polynomial endomorphism.
Because a nonexceptional polynomial endomorphism has exactly one point with finite backward orbit, every conjugation between two nonexceptional polynomial endomorphisms is an affine automorphism.  So, in (ii) of Proposition \ref{prosplitendocurve}, if $f_j, j\in I_i$ are polynomials, one may ask that all $g_j$ and all $\pi_{s/j}$ are polynomials. Moreover, $(g_j; \pi_{s/j}, s\in J_j)$ is unique up to changing the affine coordinate on $\A^1\subseteq \P^1.$
\end{rem}

\proof[Proof of Proposition \ref{prosplitendocurve}](i) is directly implied by Corollary \ref{corprodindind}.
(ii) is implied by Proposition \ref{proinvsubvspl} by induction on $|I_i|.$
\endproof

\section{Transcendence of B\"ottcher coordinates}
\subsection{B\"ottcher coordinates}
In this section, we recall the definition and some basic properties of B\"ottcher coordinates.

Let $$f(z)=a_dz^d+\dots+ a_0\in K[z]$$ be a polynomial of degree $d\geq 2.$ 
 To avoid confusion, in this section, {\bf we use $f^n$ to denote the $n-$ power of $f$ and $f^{(n)}$ for the $n$-th iterate of $f.$}

A B\"ottcher coordinate of $f$ is a Laurent series $\phi_f(z)\in \overline{K}((z^{-1}))$
satisfying $$(\phi_f\circ f)(z)=\phi_f(z)^d$$ and
of order $-1$ at $\infty$ i.e.  
it takes form 
$$\phi_f(z)=b_1z+b_0+b_{-1}/z+b_{-1}/z^2\dots.$$

\medskip

\begin{pro}\label{probottcherexists}When $\Char\, K\not| d$, the B\"ottcher coordinate $\phi_f$ exists and have the following properties:
\begin{points}
\item $b_1^{d-1}=a_d$;
\item the coefficients $b_i, i\leq 1$ are indeed contained in $K(b_1).$
\end{points}
\end{pro}
The proof of Proposition \ref{probottcherexists} can be found in \cite[Section 2.4]{FavreGauthier} when $\Char\, K=0$. The same proof holds when $\Char\, K\not| d$.

\begin{rem}
When $\Char\, K | d$, B\"ottcher coordinate may not exist.  One may check that when $p=\Char\, K>0$, the B\"ottcher coordinate does not exist for $f(z)=z^p+z^{p-1}$.
\end{rem}

We now assume that $\Char\, K\not| d$.
We view $f$ as an endomorphism of $\A^1$ which extends to an endomorphism of $\P^1.$
Fix a place $v\in M_K.$ 
The following property is \cite[Proposition 2.13]{FavreGauthier} when $\Char\, K=0$,
the same proof works when $\Char\, K\not| d$. 
\begin{pro}\label{proconvergencebottcher}There is $B_v>0$, such that $\phi_f(z)$ converges in the neighborhood of infinity $\Omega_v(f):=\{x\in \P^1_v|\,\, |z(x)|_v> B_v\}$ and for every $z\in \Omega_v$, $f^{(n)}(z)\to \infty$.
\end{pro}

\medskip

The following lemma implies that the B\"ottcher coordinate for $f$ is unique, up to multiplying by a $(d-1)$-th root of unity. 
\begin{lem}\label{lemuniquemultibottcher}Let $s\geq 1$ be a positive integer. Let $\psi\in \overline{K}((z^{-1}))$ be a Laurent series of order $-s$ at $\infty$ which
satisfies $$(\psi\circ f)(z)=\psi(z)^d.$$
Then we have $\psi(z)=\mu\phi_f(z)^s$ for some $(d-1)$-th root of unity $\mu$.
\end{lem}

\begin{rem}\label{remiteratebott}For every $n\geq 1$, the B\"ottcher coordinates of $f^n$ are Laurent series taking forms $\mu\phi_f(z)$ where $\mu$-is a $(d^n-1)$-th root of unity.
\end{rem}

\proof[Proof of Lemma \ref{lemuniquemultibottcher}]
View $\overline{K}((z^{-1}))$ as a non-archimedean field with the $z^{-1}$-adic norm $|\cdot |$.
It is clear that $f^{(n)}(z)\to \infty$.

Both $\psi(y)$ and $\phi_f(y)^s$ are of order $-s$ at $\infty$ and have coefficients in $\overline{K}$.
So there is $c\in \overline{K}^*$ such that $$(\psi/\phi_f^s)(y)\to c$$ when $y\in \overline{K}((z^{-1}))$ tends to $\infty.$
Because  $(\phi_f\circ f)^s(z)=(\phi_f(z)^s)^d$ and $(\psi\circ f)(z)=\psi(z)^d$, we get 
\begin{equation}\label{equationpsiphi}(\psi/\phi_f^s)(f^{(n)}(z))=(\psi/\phi_f^s(z))^{d^n},
\end{equation}
for every $n\geq 0.$
Then we get 
\begin{equation}\label{equationpsiphilim}\lim_{n\to \infty}(\psi/\phi_f^s(z))^{d^n}=\lim_{n\to \infty}(\psi/\phi_f^s)(f^{(n)}(z))=c.
\end{equation}
In particular, we get 
$$1=\lim_{n\to \infty}(\psi/\phi_f^s(z))^{d^n}/(\psi/\phi_f^s(z))^{d^{n-1}}=(\lim_{n\to \infty}(\psi/\phi_f^s(z))^{d^{n-1}})^{d-1}=c^{d-1}.$$
Hence $c$ is a $(d-1)$-th root of unity. Set $g:=\psi/(c\phi_f^s(z))\in \overline{K}((z^{-1})).$  By Equality \ref{equationpsiphilim}, 
$\lim\limits_{n\to \infty}g^{d^n}=1.$ 
\begin{lem}\label{lemgmeq1}There is $m\geq 1$ such that $g^{d^m}=1.$
\end{lem}
In particular, we get $g\in \overline{K}.$ By Equality \ref{equationpsiphi} with $n=1$, we get $g^{d-1}=1$. Hence 
$\psi(z)/\phi_f(z)^s=cg$ is a $(d-1)$-th roots of unity. 
\endproof
\proof[Proof of Lemma \ref{lemgmeq1}]
There is $m\geq 1$ such that $|g^{d^m}-1|<1$. 
Set $\delta:=g^{d^m}-1$ and assume that $\delta\neq 0.$
Then for $n\geq m$, we get 
$$g^{d^n}=(1+\delta)^{d^{n-m}}=1+d^{n-m}\delta+\epsilon$$
where $|\epsilon|<|\delta|$.
Then $|g^{d^n}-1|=|d^{n-m}\delta|=|\delta|>0$ which is a contradiction.
\endproof

\begin{cor}\label{corsemibott}If $g\gtrsim_{\pi}f$, where $f,g$ are polynomials of degree $d$, then $\phi_f\circ \pi=\mu\phi_g^{\deg(\pi)}$,  where $\mu$-is a $(d-1)$-th root of unity.
\end{cor}
\proof[Proof of Corollary \ref{corsemibott}]Because $$(\phi_f\circ\pi)(g(z))=(\phi_f\circ f)(\pi(z))=(\phi_f\circ\pi)(z)^d,$$ we conclude the proof by Lemma \ref{lemuniquemultibottcher}.
\endproof

Since $\id$ is a B\"ottcher coordinate for the map $z\mapsto z^d$, we get the following result.
\begin{cor}\label{cormono}When $f$ is of monomial type, then every B\"ottcher coordinate of $f$ is algebraic.  
\end{cor}

\subsection{Orbit of a point}\label{subsectionorbit}
Let $X$ a variety over a field $\bk$. Let $f: X\to X$ be an endomorphism and $x\in X$. We denote by $O_f(x)$ the $f$-orbit of $x$.
Denote by $Z_f(x)$ the Zariski closure of $O_f(x)$. It is clear that $f(Z_f(x))\subseteq Z_f(x).$
In this section, we study the structure of $Z_f(x)$ and the action of $f$ on it.

\medskip

Write $Z_f(x)= T_f(x) \sqcup P_f(x)$ where $T_f(x)$ is the union irreducible components of $Z_f(x)$ of dimension $0$ and $P_f(x)$ is the union irreducible components of $Z_f(x)$ of dimension $\geq 1.$ Note that $P_f(x)=\emptyset$ if and only if $x$ is preperiodic for $f.$
If $x\in X(\bk)$, then every irreducible component contains a Zariski dense set of $\bk$ points, hence geometrically irreducible by \cite[Lemma 1]{Abramovich1992}.

It is clear that $f(P_f(x))\subseteq P_f(x).$ Set $t_f(x):=|T_f(x)|\geq 0$. Then $T_f(x)=\{f^n(x)|\,\, 0\leq n\leq t_f(x)-1\}.$
Becasue $O_f(x)$ is contained in $\{f^n(x)|\,\, 0\leq n\leq t_f(x)\}\cup f(P_f(x))$ and every irreducible component of $P_f(x)$ has positive dimension, 
$P_f(x)=\overline{f(P_f(x))}.$ Hence every irreducible component of $P_f(x)$ is periodic.
Let $Z$ be an irreducible component of $Z_f(x)$ containing $f^{t_f(x)}(x).$ There is a minimal $m\geq 1$ such that $f^m(Z)\subseteq Z.$
Then $P_f(x)=\cup_{i=1}^{m-1} \overline{f^i(Z)}.$ 

\begin{lem}\label{lemuniqueirrper}
For every $n\geq t_f(x)$, there is a unique $i\in \{0,\dots,m-1\}$ such that $f^n(x)\in \overline{f^i(Z)}.$
\end{lem}
\proof[Proof of Lemma \ref{lemuniqueirrper}]Assume that there are $i\neq j\in\{1,\dots,m-1\}$ such that $f^n(x)\in \overline{f^i(Z)}\cap \overline{f^j(Z)}.$ Because  $f^m(\overline{f^i(Z)}\cap \overline{f^j(Z)})\subseteq \overline{f^i(Z)}\cap \overline{f^j(Z)}$, $O_f(f^{n}(x))\subseteq \cup_{j=1}^{m-1} (\overline{f^{i+j}(Z)}\cap \overline{f^j(Z)})$ which is not dense in $P_f(x).$ We get a contradiction.
\endproof

\medskip

Now we summarize what we get. 
For $a\in \Z$ and $b\in \Z^+$, denote by $e(a/b)$ the unique element in $\{0,\dots, b-1\}$ such that $a=e(a/b) \mod b.$ 
For any endomorphism $f: X\to X$ and $x\in X$, we may associated to them the following datas: Integers $t_f(x)\geq 0$ and $p_f(x)\geq 1$, a closed irreducible subvariety $P^+_f(x)$ such that $$Z_f(x)=\{f^n(x)|\,\, 0\leq n\leq t_f(x)-1\}\sqcup (\cup_{i=1}^{p_f(x)}\overline{f^i(P^+_f(x))}),$$ and for every $i\geq t_f(x)$,  the unique $j\in \{0,\dots, p_f(x)\}$ satisfying 
$\overline{f^{j}(P^+_f(x))})$ is $e(i/p_f(x))$. In particular $\overline{f^{p_f(x)}(P^+_f(x))}=P^+_f(x).$

\subsection{Canonical heights}\label{subsectioncanheight}
In this section, we recall some basic facts on canonical height for polynomial dynamics.
All these are well known and can be found in \cite{Silverman2007} and \cite{FavreGauthier}.

\medskip

Let $f\in K[x]$ be a polynomial of degree at least $2$.
The canonical Green functions and canonical height of $a\in K$ for $f$ are
$$g_{f,v}(a):=\lim_{n\to \infty}\deg(f)^{-n}\log\max\{1,|f^n(a)|_v\}$$
and 
$$\hat{h}_{f}(a):=\sum_{v\in M_K}n_vg_{f,v}(a).$$
By Tate's limiting argument, these limits exists and non-negative.  
Moreover $g_{f,v}(a)=0$ for all but finitely many $v\in M_K$, hence $\hat{h}_{f}(a)$ is well defined and non-negative.
We also concern the multiplicative Green functions and canonical height 
$$G_{f,v}(a):=e^{g_{f,v}(a)},$$
$$\hat{H}_f(a):=e^{\hat{h}_{f}(a)}=(\prod_{v\in M_K}G_{f,v}(a)^{n_v})^{1/[K:\F]}$$
for $a\in K.$

Recall that the naive height of $a$ is 
$$h(a)=[K:\F]^{-1}\sum_{v\in M_K}n_v\log\max\{1,|a|_v\}.$$
Then we have  
$$\hat{h}_f(a):=\lim_{n\to \infty}\deg(f)^{-n}h(f^n(a)).$$

If we change $K$ by any finite extension of $K$, the values of naive and canonical heights do not change. So $h$ and $\hat{h}_f$ are indeed functions from 
$\overline{K}$ to $[0,+\infty).$
We have the following basic properties:
\begin{points}
\item $\hat{h}_f(\cdot)-h(\cdot)$ is bounded on $\overline{K};$
\item $\hat{h}_f(f(a))=\deg(f)\hat{h}_f(a)$ for every $a\in \overline{K};$
\item for $n\geq 1$ and $a\in \overline{K},$ $\hat{h}_f(a)=\hat{h}_{f^n}(a).$
\end{points}

\medskip

There is a way to compute the canonical Greens functions using B\"ottcher coordinates:
Assume that $\Char\, K\not| d.$ For every $a\in K,$ we have 
\begin{equation}\label{equgreenbott}
  g_{f,v}(a):=
   \begin{cases}
   0&\mbox{if $\{|f^n(x)|_v\}_{n\geq 0}$ is bounded};\\
   |\phi_f(f^m(a))|^{1/d^m}_v&\mbox{if $f^m(x)\in \Omega_{v}(f)$ for some $m\geq 0$}.
   \end{cases}
  \end{equation}
  Recall $\Omega_{v}(f)$ is defined in Proposition \ref{proconvergencebottcher}.
Note that if $\{|f^n(x)|_v\}_{n\geq 0}$ is not bounded, then there is $m\geq 0$ such that $f^m(x)\in \Omega_{v}(f)$, moreover $|\phi_f(f^m(a))|^{1/d^m}_v$ does not depend on the choice of $m$.

\subsection{Equivalent dynamical pairs}\label{subsectionequidp}
For $d\geq 2$ and $\Char\,K\not| d$,  let $\sD_d$ be the set of polynomials $f\in \overline{K}[z]$ which is nonexceptional as an endomorphism of $\P^1$.
When $\Char\,K=0$, such polynomials are called non-integrable in \cite{FavreGauthier} and \cite{Nguyen}.
\begin{rem}A polynomial of degree $d$ is not in $\sD_d$ if and only if it is of monomial type. 
\end{rem}
Denote by $\sP_d$ the set of $(f,a)\in \sD_d\times \overline{\F}$ such that $a$ is $f$-preperiodic. 
\begin{rem}
When $\F=\Q$, by the Northcott property, $(f,a)\in \sP_d$ if and only if $\hat{h}_f(a)=0.$
\end{rem}

Denote by  $\sD\sP_d:=(\sD_d\times \overline{\F})\setminus \sP_d$. An element $(f,a)\in \sD\sP_d$ is called a \emph{dynamical pair}.
For $(f,a), (g,b)\in \sD\sP_d$,   we say that $(f,a)$ and $(g,b)$ are equivalent and write   $(f,a)\sim(g,b)$ if
the $(f\times g)$-orbit of $(a,b)\in (\P^1)^2$ is not Zariski dense. In this case every irreducible component $C$ of $P_{f\times g}((a,b))$ is a curve.
Let $\pi_i:( \P^1)^2\to \P^1, i=1,2$ be the projection to the $i$-th coordinate. Then $\pi_i(C)=\P^1$ and the ratio $$d((f,a)/(g,b)):=\deg\pi_1|_C/ \deg\pi_2|_C\in \Q^*$$ does not depend on the choice of irreducible component $C.$  
 If $(f,a)\sim (g,b)$ and $(g,b)\sim(h,c)$, then by Proposition \ref{prosplitendocurve},  $(f,a)\sim (h,c)$. It shows that $\sim$ is an equivalent relation. 
Moreover, we have $$d((f,a)/(h,c))=d((f,a)/(g,b))d((g,b)/(h,c)).$$
Hence, for a finite set of equivalent dynamical pairs $\alpha_i, i\in I$, it defines a point 
$d(\alpha_i, i\in I)\in \P^{|I|-1}\setminus (\cup_{i\in I}\{z_i=0\})$ such that $$(z_i/z_j)(d(\alpha_i, i\in I))=d(\alpha_i/\alpha_j)$$ for every $i,j\in I.$
\begin{rem}\label{remeqpairpoly}
If $(f,a)\sim(g,b)$, then $f\sim g$.
\end{rem}

Now assume that $\Char\,\F\not| d.$
Let $(f,a),(g,b)\in \sD\sP_d$ be two equivalent dynamical pairs. Let $K$ be a finite extension of $\F$, which contains $a,b$, all coefficients and all $(d-1)$-roots of leading coefficients of $f$ and $g$. For $v\in M_K$,  we have $\lim\limits_{n\to \infty}|f^n(a)|_v\to \infty$ if and only if $\lim\limits_{n\to \infty}|g^n(b)|_v\to \infty$. By Corollary \ref{corsemibott}, if $a\in \Omega_v(f)$ and $b\in \Omega_v(g)$, then $\phi_f(a)^n=\phi_g(b)^m$ for some $m,n\in \Z_{\geq 1}$ satisfying $m/n=d((f,a)/(g,b)).$
Then by Equality \ref{equgreenbott}, we get the following result.
\begin{pro}\label{proequpairheight} Assume that $\Char\, \F\not| d.$
Let $(f,a),(g,b)\in \sD\sP_d$ be equivalent dynamical pairs defined over a number field $K.$
Then the following holds:
\begin{points}
\item For every $v\in M_K$, $g_{f,v}(a)=0$ if and only if $g_{g,v}(b)=0$. Moreover, if $g_{f,v}(a)\neq 0$, then $g_{f,v}(a)/g_{g,v}(b)=d((f,a)/(g,b)).$ 
\item The canonical heights $\hat{h}_{f}(a)$ and $\hat{h}_{g}(b)$ are zero at the same time. Moreover, if they are not zero,
$\hat{h}_{f}(a)/\hat{h}_{g}(b)=d((f,a)/(g,b)).$
\end{points}
\end{pro}

\subsection{Transcendence of B\"ottcher coordinates}\label{subsectionbottcoortran}
Let $f_1,\dots, f_r$ be polynomials of degree $d$. Assume that $K$ contains all coefficients of every $f_i$ and the $(d-1)$-th roots of its leading coefficient.
For $v\in M_K$, let $\Omega_v(f_i)$ be defined as in Proposition \ref{proconvergencebottcher}. 
We view the B\"ottcher coordinates $\phi_{f_i}$ as functions on $\Omega_v(f_i)$.

For $i=1,\dots,r$, let $a_i$ be a $K$-point in $\Omega_v(f_i)$.  We will show that the algebraicity of $\prod_{i=1}^r\phi^{n_i}_{f_i}(a_i)$ is completely determined by some 
geometric datas associated to $(f_1,\dots, f_r, a_1,\dots, a_r).$ 
\begin{rem}\label{remmono}
By Corollary \ref{cormono}, if $f_i$ is of monomial type, then $\phi_{f_i}(a_i)$ is algebraic.
\end{rem}
Now we assume that  $f_1,\dots, f_r\in \sD_d.$
Our assumption shows that $(f_i,a_i)\in \sD\sP_d$ for every $i=1,\dots,r.$

\subsubsection{Geometric datas}
In this section we will associated to $(f_1,\dots, f_r, a_1,\dots,a_r)$ a partition $\{1,\dots,r\}=\sqcup_{j=1}^l J_j$ and a group of positive integers $d_{s/j}\geq 1, s\in J_j, j=1,\dots,l.$

\medskip

Consider the endomorphism $F:=\prod_{i=1}^rf_i: \A^r\to \A^r$ which extends to an endomorphism of $(\P^1)^r.$
Set $a:=(a_1,\dots, a_r)\in (\P^1)^r.$ Following Section \ref{subsectionorbit}, we may associated to $F, a$ the data $(t_F(a), p_F(a), P_F^+(a)).$
In particular, $P_F^+(a)$ is $F$ periodic with the minimal positive period $p_F(a).$
Because $\lim\limits_{n\to \infty} |f^n(a_i)|_v\to \infty, i=1,\dots, r$, the projection of $P_F^+(a)$ to every factor of $(\P^1)^r$ is dominant.

\medskip

By Proposition \ref{prosplitendocurve} and Proposition \ref{proreformulationii}, we may associated to $p_F(a), P_F^+(a)$ the following data: 
there is a partition $\{1,\dots,r\}=\sqcup_{j=1}^lJ_j$, a polynomial $g_j: \P^1\to \P^1, j=1,\dots,l$ such that $g_j\gtrsim_{\pi_{s/j}} f_s^{p_F(a)}, s\in J_j$ where $\pi_{s/j}$ is a polynomial of degree $d_{s/j}$  such that 
for every $j=1,\dots,l$, $$\Phi_j:=(\pi_{s/j})_{s\in J_j}: \P^1\to (\P^1)^{|J_j|}$$ is birational to its image and $P_F^+(a)=(\prod_{j=1}^l\Phi_j(\P^1)).$ Here we identified $(\P^1)^{|I_i|}$ with $(\prod_{j=1}^l\P^1)^{|J_j|}.$
The partition $\{1,\dots,r\}=\sqcup_{j=1}^l J_j$ is unique and the $(g_j; \pi_{s/j}, s\in J_j)$ is unique up to changing the affine coordinate on $\A^1\subseteq \P^1.$
In particular, the degrees $d_{s/j}\geq 1, s\in J_j, j=1,\dots,l$ does not depend on the choice of $(g_j; \pi_{s/j}, s\in J_j), s\in J_j, j=1,\dots,l.$

\begin{rem}\label{remgeomdata}
Another discerption of the datas $\{0,\dots,r\}=\sqcup_{j=1}^l J_j$ and $d_{s/j}\geq 1, s\in J_j, j=1,\dots,l$  is as follows:
\begin{points}
\item[(a)]There is a unique partition $\{0,\dots,r\}=\sqcup_{j=1}^l J_j$, such that 
$P_F^+(a)=\prod_{j=1}^lC_j$ after identifying $(\P^1)^{|I_i|}$ with $(\prod_{j=1}^l\P^1)^{|J_j|},$ where $C_j, j=1,\dots,l$ is a $(\prod_{s\in J_j}f_s)^{p_f(a)}$-invariant curve of $(\P^1)^{|J_j|}.$ Indeed, $C_j=\Phi_j(\P^1)$ for $j=1,\dots,l$.
\item[(b)]For $s\in J_j, j=1,\dots,l$, $d_{s/j}=\deg (\pi_s|_{C_j})$ where $\pi_s: (\P^1)^{l}\to \P^1$ is the projection to the $s$-th coordinate.
\end{points}

\begin{rem}\label{rempartition}
 A simple way to view the partition $\{0,\dots,r\}=\sqcup_{j=1}^l J_j$ is as follows:
 For  $s,t\in \{0,\dots,r\}$, they are contained in the same $J_j$ for some $j=1,\dots,l$ if and only if 
 $(f_s,a_s)\sim (f_t,a_t).$ 
  Moreover for $s,t\in J_i$, we have $d_{s/j}/d_{t/j}=d((f_s,a_s)/(f_t,a_t)).$
\end{rem}

\end{rem}
\begin{rem}\label{remiteratedataorbit}Easy to see from Remark \ref{remgeomdata} that,
if we replace $P_F^+(a)$ by $F^n(P_F^+(a))$ for any $n\geq 0$, we get the same datas $\{0,\dots,r\}=\sqcup_{j=1}^l J_j$ and $d_{s/j}\geq 1, s\in J_j, j=1,\dots,l$.
Hence, if one replace $(f_1,\dots, f_r, a_1,\dots,a_r)$ by $(f_1^m,\dots, f_r^m, f_1^n(a_1),\dots,f_r^n(a_r))$ for any $m\geq 1, n\geq 0$, we get the same datas $\{0,\dots,r\}=\sqcup_{j=1}^l J_j$ and $d_{s/j}\geq 1, s\in J_j, j=1,\dots,l$.
\end{rem}

\subsubsection{Transcendence of B\"ottcher coordinates}
\begin{thm}[=Theorem \ref{thmtranbottintro}]\label{thmtranbott}
The following statements hold.
\begin{points}
\item For integers $n_1,\dots, n_r$, $\prod_{i=1}^r\phi^{n_i}_{f_i}(a_i)$ is either transcendence over $K$ or a root of unity.
\item The product $\prod_{i=1}^r\phi^{n_i}_{f_i}(a_i)$ is a root of unity, if and only if, for every $j=1,\dots,l,$
$\sum_{s\in J_j}n_sd_{s/j}=0.$
\end{points}
\end{thm}
Part (i) of Theorem \ref{thmtranbott} was proved in \cite[Theorem 1.4]{Nguyen} in the number field case.
Part (ii) of Theorem \ref{thmtranbott} answers the first question proposed in \cite[Section 4.3]{Nguyen}.

\proof[Proof of Theorem \ref{thmtranbott}]
Assume that $b:=\prod_{i=1}^r\phi^{n_i}_{f_i}(a_i)$ is algebraic over $K.$ After replacing $K$ by a finite extension, we may assume that $b\in K.$
Let $X:=(\P^1)^{r+1}$ and let $G: X\to X$ be the endomorphism $$G: (z_1,\dots, z_l, y)\mapsto (f_1(z_1),\dots, f_l(z_l), y^d).$$ Set $x:=(a_1,\dots,a_r, b)\in X(K).$
For $n\geq 0$, set $$x_n:=G^n(x)=(f_1^n(a_1),\dots, f_l^n(a_l), b^{d^n})\in X(K).$$
Let $\pi: (\P^1)^{r+1}\to (\P^1)^{r}$ be the projection to the first $r$-coordinates. We have $\pi\circ G=F\circ \pi.$
Set $H_{\infty}:=(\P^1)^{r+1}\setminus (\A^1)^{r+1}.$ Set $H:=\cap_{m\geq 0}\overline{\{x_n, n\geq m\}}\subseteq X.$ 
Because $\{x_n, n\geq m\}\cap H_{\infty}=\emptyset$, no irreducible component of $H$ is contained in $H_{\infty}.$
It is clear that $\pi(H)=P_F(a).$
By Remark \ref{remiteratedataorbit}, we may replace $(f_1,\dots, f_r, a_1,\dots,a_r)$ by $(f_1^m,\dots, f_r^m, f_1^n(a_1),\dots,f_r^n(a_r))$ for some suitble $m\geq 1, n\geq 0$ to assume that $H=\overline{\{x_n, n\geq 0\}}$ and is irreducible. Hence the $Z_F(a)=P_F(a)=P_F^+(a)$ is irreducible. In this case, $t_F(a)=0$ and $p_F(a)=1.$

\medskip

We first prove the ``if" part of (ii).
Because $Z_F(a)=P_F^+(a)$, there is $(u_1,\dots,u_l)\in (\prod_{j=1}^l\Phi_j)^{-1}((a_1,\dots,a_r))$. 
By Corollary \ref{corsemibott}, we have 
$$\prod_{j=1}^l\phi^{\sum_{s\in J_j}n_sd_{s/j}}_{g_j}(u_j)=\prod_{j=1}^l\prod_{s\in J_j}\phi^{n_sd_{s/j}}_{g_j}(u_j)=\mu\prod_{i=1}^r\phi^{n_i}_{f_i}(a_i)$$
for some $(d-1)$-th root of unity. If $\sum_{s\in J_j}n_sd_{s/j}=0$ for every $j=1,\dots,l,$ the left hand side of the above equality is $1$, so $\prod_{i=1}^r\phi^{n_i}_{f_i}(a_i)$ is a root of unity.

\medskip

Now we only need to prove (i) and the ``only if" part of (ii).
Consider the analytic subvariety $\sV'\subseteq \prod_{i=1}^r(\Omega_v(f_i)\setminus \{\infty\})\times \A^1_v$ defined by $y=\prod\phi_{f_i}^{n_i}(z_i).$
Because $\phi_{f_i}, i=1,\dots, r$ is meromorphic at $\infty$, $\sV'$ extends to a closed analytic subvariety $\sV\subseteq \prod_{i=1}^r(\Omega_v(f_i))\times \P^1_v.$
Set $Z_v:=\sV\cap H_{\infty}$. Because all coefficients of $\phi_{f_i}, i=1\dots,r$ are contained in $K$, $Z_v$ is the analytification of a $K$-subvariety of $H_{\infty}$ and 
$\sV$ is defined over $K$ along $Z.$ 

Because $|f_i^n(a_i)|_v\to \infty$ for every $i=1,\dots,r$, $\lim\limits_{n\to\infty}x_n\subseteq H_{\infty,v}$ in the $v$-adic topology.
Because $x_n\in \sV$ for $n\geq 0$,  $\lim\limits_{n\to\infty}x_n\subseteq Z_v.$
Let $h$ be the naive height on $\P^1(K)$ i.e. for every $z\in \P^1(K)$, $h_w(z)=\log \max\{1,|z|_w\}, w\in M_K$ and $h(z)=[K:\F]^{-1}\sum_{w\in M_K}n_wh_w(z).$
There is $C_1>1$ such that $C_1^{-1}d^n\leq h_v(x_n)\leq h(x_n)\leq C_1 d^n$ for every $n\geq 0.$ By (vi) of Facts \ref{facts}, 
$h_v(z)=g_{Z_v/\sV,v}(z)+O(1)$ for every $z\in \sV(K_v).$ Then all assumptions in Theorem \ref{thmgengeneralalgebraicity} are satisfied. 
By Theorem \ref{thmgengeneralalgebraicity}, there is an open neighborhood $\Omega$ of $(\infty,\dots,\infty)\in (\P^1)^r$ in $\prod_{i=1}^r(\Omega_v(f_i)$ such that $F(\Omega)\subseteq \Omega$, $\sV\cap H_v\cap \pi^{-1}(\Omega)\subseteq H\cap \pi^{-1}(\Omega)$ and $\dim H=\dim \sV\cap H_v.$  
In particular, $\sV\cap H_v\cap \pi^{-1}(\Omega')$ contains a non empty open subset $\sW$ of $H_v\cap \pi^{-1}(\Omega'),$ where $\Omega':=\Omega\setminus \{(\infty,\dots,\infty)\}.$
Because $G(x_n)=x_{n+1},$ $H$ is $G$-invariant. 
Because $\sV\setminus H_{\infty, v}=\sV'$ is a section of $\pi$ on $\prod_{i=1}^r(\Omega_v(f_i)\setminus \{\infty\})$, $H\setminus H_{\infty}\neq\emptyset$ and 
$\dim H=\dim \sV\cap H_v$, $H\neq \pi(H)\times \P^1.$
By (i) of Proposition \ref{prosplitendocurve}, $H$ takes form $\pi(H)\times c$ where $c$ is a $(d-1)$-th root of unity.
Because $x_0\in H$, we get $b=c$ is a $(d-1)$-th root of unity, which implies (i).  
Then $\pi(\sW)$ is open in $\pi(H)$ and $\pi^{-1}(\pi(\sW))\cap \sV=\sW.$

Note that $\pi(H)=P_F^+(a)$ is an irreducible closed $F$-invariant subvariety.
Moreover, we get
 $$\{(z_1,\dots, z_r, b)|\,\, (z_1,\dots, z_r)\in \pi(\sW)\}=\sW$$
 $$=\{(z_1,\dots, z_r, \prod_{i=1}^r\phi^{n_i}_{f_i}(z_i))|\,\, (z_1,\dots, z_r)\in \pi(\sW)\}.$$
Hence for every $(z_1,\dots, z_r)\in \pi(\sW)$, one has $\prod_{i=1}^r\phi^{n_i}_{f_i}(z_i)=b.$
By Corollary \ref{corsemibott}, for every $(w_1,\dots,w_l)\in (\prod_{j=1}^l\Phi_j)^{-1}(\pi(\sW))$, we get 
\begin{equation}\label{equationproduclp}\prod_{j=1}^l\phi^{\sum_{s\in J_j}n_sd_{s/j}}_{g_j}(w_j)=\prod_{j=1}^l\prod_{s\in J_j}\phi^{n_sd_{s/j}}_{g_j}(w_j)=b\mu
\end{equation}
for some $(d-1)$-root of unity $\mu.$
We note that $\prod_{j=1}^l\phi^{\sum_{s\in J_j}n_sd_{s/j}}_{g_j}(w_j)$ is a meromorphic function on $(\prod_{j=1}^l\Phi_j)^{-1}(\Omega)$ and 
$(\prod_{j=1}^l\Phi_j)^{-1}(\pi(\sW))$ is a nonempty open subset of $(\prod_{j=1}^l\Phi_j)^{-1}(\Omega),$
Equality \ref{equationproduclp} holds on $(\prod_{j=1}^l\Phi_j)^{-1}(\Omega).$
We note that $(\prod_{j=1}^l\Phi_j)^{-1}(\Omega)$ is a neighborhood of $(\infty,\dots,\infty)\in (\P^1)^l$.
Because for every $j=1,\dots,l$, $\lim\limits_{w\to\infty}|\phi_{g_j}(w)/w|_v= c_j$ for some $c_j\in (0,+\infty)$, there is $B>1$, such that for $(w_1,\dots,w_l)\in (\prod_{j=1}^l\Phi_j)^{-1}(\Omega(K_v))$,  $$B^{-1}\leq |\prod_{j=1}^lw_j^{\sum_{s\in J_j}n_sd_{s/j}}|_v\leq B.$$
Hence for every $j=1,\dots,l,$
$\sum_{s\in J_j}n_sd_{s/j}=0.$ This concludes the proof.
\endproof
\subsection{Algebraicity and $\Q$-linear relations of canonical heights}\label{subsectionalgcanheight}
In this section, $K$ is a number field. Denote by $M_K^{\infty}$ the archimedean places of $K$ and $M_K^{f}$ the nonarchimedean places of $K$.

\begin{pro}\label{pronguyen}For $f\in \sD_d$ with coefficients in $K$ and $a\in K$, $\hat{H}_f(a)$ is algebraic if and only if for every $v\in M_K^{\infty}$, $|f^n(a)|_v, n\geq 0$ is bounded.
\end{pro}
This result is \cite[Corollary 1.6]{Nguyen}. It is directly implied by (i) of Theorem \ref{thmtranbott} and 
the following well known fact. 
\begin{lem}\cite[Lemma 2.1]{Nguyen}\label{lemnonarchialgebraic} If $v\in M_K^f$, then $G_{f,v}(a)=p^c$ for some $c\in \Q$.
\end{lem}
\begin{rem}\label{remcomputable}The proof of Lemma \ref{lemnonarchialgebraic} shows that $G_{f,v}(a)$ is indeed computable.
\end{rem}

In \cite[Section 4.3]{Nguyen}, Nguyen suggested to study the $\Q$-linear relations of the canonical heights.
He also proved a partial result \cite[Corollary 1.8]{Nguyen} in this direction.
He suspected that some result like (ii) of Theorem \ref{thmtranbott} may be helpful for this problem.
In this section, we follows Nguyen's suggestion to get some applications of Theorem \ref{thmtranbott}.

\medskip

%

Recall that we have defined an equivalence relation $\sim$ on $\sD\sP_d$ in Section \ref{subsectionequidp}.
We now introduce another equivalence relation:
Set $\sG:=\Gal(\overline{\Q}/\Q)$. For $(f,a), (g,b)\in \sD\sP_d$, we write $(f,a)\sim_w(g,b)$ if there is $\sigma\in \sG$ such that 
$(f,a)\sim \sigma(g,b):=(\sigma(g),\sigma(b)).$ 
We say such $(f,a), (g,b)$ are \emph{weakly equivalent}.
Because $\sim$ is an equivalence relation and $\sG$ is a group, $\sim_w$ is an equivalence relation.
It is clear that $(f,a)\sim(g,b)$ implies $(f,a)\sim_w(g,b).$ 

\begin{rem}
Weak equivalence does not implies equivalence.
For example, let $f_1=f_2=z(z+1/2)$, $a_1=3-2\sqrt{2}$ and $a_2=3+2\sqrt{2}$. 
Because $(f_1,a_1)$ and $(f_2,a_2)$ are conjugate by some $\sigma\in \sG,$
$(f_1,a_1)\sim_w(f_2,a_2)$. On the other hand, when $n\to \infty$, $|f_1^n(a_1)|\to 0$ and $|f_2^n(a_2)|\to \infty.$
By Proposition \ref{proequpairheight}  $(f_1,a_1)\not\sim(f_2,a_2).$
\end{rem}

\begin{rem}\label{remweakeqheight}
For $(f,a)\in \sP\sD_d$ and $\sigma\in \sG$, $\hat{h}_f(a)=\hat{h}_{\sigma(f)}(\sigma(a)).$
Hence, by Corollary \ref{proequpairheight}, for $(f,a)\sim_w(g,b)$, 
$\hat{h}_f(a)/\hat{h}_f(b)\in \Q^*.$
\end{rem}

For a finite set of weakly equivalent elements $(f_i,a_i)\in \sD\sP_d, i=1,\dots,r$, define
$$Q((f_i,a_i), i=1,\dots,r):=[\hat{h}_{f_1}(a_1):\dots : \hat{h}_{f_r}(a_r)]\in \P^{r-1}(\Q)\setminus (\cup_{i=1}^r\{z_i=0\}).$$ 
It can be computed as follows:
there are $\sigma_i\in \sG, i=1,\dots,r$, such that $\sigma_i(f_i, a_i), i=1,\dots,r$ are equivalent.
Then we have $$Q((f_i,a_i), i=1,\dots,r)=d(\sigma_i(f_i,a_i), i=1,\dots,r).$$
\begin{rem}\label{remcomputequiheight}
Thought the definition of $Q((f_i,a_i), i=1,\dots,r)$ involves the canonical heights, the above argument shows that, up to some Galois conjugations, it is indeed  a geometric invariant.
\end{rem}

\medskip

Define $\sT_d$ to be the set of $(f,a)\in \sD_d\times \overline{\Q}$ for which there is an embedding $\tau: \overline{\Q}\hookrightarrow \C$ such that $|\tau(f^n(a_i))|, n\geq 0$ is unbounded. We have $\sT_d\subseteq \sD\sP_d$. Moreover,
by Proposition \ref{pronguyen}, $(f,a)\in\sT_d$ if and only if $\hat{H}_f(a)$ is not algebraic. 
By Lemma \ref{remweakeqheight}, for $(f,a), (g,b)\in \sD\sP_d$, if $(f,a)\sim_w(g,b)$, they are contained in $\sT_d$ in the same time. 

\medskip

The following result generalizes Proposition \ref{pronguyen}.
\begin{pro}[=Proposition \ref{prohproductalgintro}]\label{prohproductalg}For $(f_i,a_i),\dots, (f_r,a_r)\in \sD_d\times \overline{\Q}$, $n_1,\dots,n_r\in \Z$,
set $T:=\{i=1,\dots, r|\,\, (f_i,a_i)\in \sT_d\}.$
Let $T=\sqcup_{j=1}^l J_j$ be the partition associated to $\sim_w.$
Then $\prod_{i=1}^r\hat{H}_{f_i}(a_i)^{n_i}$ is algebraic if and only if for every $j=1,\dots,l$, $$Q((f_s,a_s),s\in J_j)\in \{\sum_{s\in J_j} n_sz_s=0\}.$$
\end{pro}

\begin{rem}\label{remequarchone}
By Proposition \ref{proequpairheight}, if $Q((f_s,a_s),s\in J_j)\in \{\sum_{s\in J_j} n_sz_s=0\},$ 
$\prod_{s\in J_j}\hat{H}_{f_i}(a_i)^{n_i}=1.$ Hence we get that $\prod_{i=1}^r\hat{H}_{f_i}(a_i)^{n_i}$ is algebraic if and only if 
$\prod_{i\in T}\hat{H}_{f_i}(a_i)^{n_i}=1$.
\end{rem}

\proof[Proof of Proposition \ref{prohproductalg}]
The ``if" part is obvious. Now we prove the ``only if" part. 

By Proposition \ref{pronguyen}, we may assume that $(f_i,a_i),\dots, (f_r,a_r)\in \sT_d.$ 
Let $K$ be a number field which contains $a_1,\dots,a_r$, all coefficients and all $(d-1)$-th roots of leading coefficients of $f_i, i=1,\dots,r$.
After enlarging $K$, we may assume that $K/\Q$ is Galois and $K$ can not be embedded in  $\R.$
After replacing $(f_i,a_i)$ by $(\sigma_i(f_i),\sigma_i(a_i))$ for some $\sigma_i\in \sG_K:=\Gal(K/\Q)$ for each $i$, we may assume that for $i,j=1,\dots, r$, 
$(f_i,a_i)\sim (f_j,a_j)$ if and only if $(f_i,a_i)\sim_w(f_j,a_j).$ Hence, for every $\sigma\in \sG_K$, $(\sigma(f_i),\sigma(a_i))\sim (\sigma(f_j),\sigma(a_j))$ if and only if $(\sigma(f_i),\sigma(a_i))\sim_w(\sigma(f_j),\sigma(a_j)).$
Now we only need to show that 
for every $j=1,\dots,l$, $\sum_{s\in J_j} n_s\hat{h}_{f_s}(a_s)=0.$

\medskip

Write
$$\prod_{i=1}^r\hat{H}_{f_i}(a_i)^{n_i}=\prod_{i=1}^r\prod_{v\in M_K}\hat{G}_{f_i,v}(a_i)^{n_in_v}$$
$$=(\prod_{i=1}^r\prod_{v\in M_K^{\infty}}\hat{G}_{f_i,v}(a_i)^{n_in_v})\times (\prod_{i=1}^r\prod_{v\in M_K^{f}}\hat{G}_{f_i,v}(a_i)^{n_in_v})$$
By Lemma \ref{lemnonarchialgebraic}, $\prod_{i=1}^r\prod_{v\in M_K^{\infty}}\hat{G}_{f_i,v}(a_i)^{n_in_v}$ is algebraic.

Fix an embedding $K\hookrightarrow \C$ to view $K$ as a subfield of $\C.$ We will omit the place $v_0$ associated to this embedding when it appears in any notation.
For example, the Green function $\hat{G}_{f_i,v_0}$ will be written as $\hat{G}_{f_i}.$
Because $K$ has no real embedding, for each $v\in \sM^{\infty}_K$, $n_v=2$ and it corresponds to exactly a pair of conjugate elements
$\sigma, \bar{\sigma}\in \sG_{K}$.  
We note that, over $\C$, $\phi_{\bar{\sigma}(f_i)}(\bar{\sigma}(x))=\overline{\phi_{\sigma(f_i)}(\sigma(x))}$ 
for every $i=1,\dots,r$ and $x\in K\cap \Omega(\sigma(f_i))=\Omega(\bar{\sigma}(f_i)).$
For every $i, \sigma$, we denote by $G(\sigma, i)$ the Green function $\hat{G}_{f_i,v}(a_i)=\hat{G}_{\sigma(f_i)}(\sigma(a_i))$ where $v$ is the place associated to $\sigma$ and by $h(i)$ the canonical height of $\hat{h}_{f_i}(a_i).$
For pairs $(i,\sigma),(j,\delta)\in \{1,\dots,r\}\times \sG_K$, write $(i,\sigma)\sim (j,\delta)$ if $(\sigma(f_i),\sigma(a_i))\sim (\delta(f_j),\delta(a_j))$ and $(i,\sigma)\sim_w (j,\delta)$ if $(\sigma(f_i),\sigma(a_i))\sim_w (\delta(f_j),\delta(a_j)).$

Let $A$ be the set of $(\sigma, i)$ such that $G(\sigma, i)>1.$ Note that $(\sigma, i)\in A$ if and only if $(\overline{\sigma}, i)\in A.$ Moreover if $(i,\sigma)\sim(j,\delta)$ then they are contained in $A$ at the same time. After replacing $a_i, i=1,\dots, r$ by $f_i^n(a_i),i=1,\dots, r$ for some large enough $n\geq 0,$ we may assume that 
$\sigma(a_i)\in \Omega(\sigma(f_i))$ for every $(\sigma, i)\in A.$ For every $(\sigma, i)\in A$, set $\Phi(\sigma,i):=\phi_{\sigma(f_i)}(a_i).$
For $u=(\sigma,i)\in A$, set $\sigma_u:=\sigma$ and $i_u:=i.$
Then $$\prod_{i=1}^r\prod_{v\in M_K^{\infty}}\hat{G}_{f_i,v}(a_i)^{n_in_v}=\prod_{i=1}^r\prod_{\sigma\in \sG_K}G(\sigma, i)^{n_i}.$$
$$=\prod_{u\in A}G(u)^{n_{i_u}}=\prod_{u\in A}|\Phi(u)|^{n_{i_u}}=\prod_{u\in A}\Phi(u)^{n_{i_u}}$$
is algebraic. By Theorem \ref{thmtranbott} and Proposition \ref{proequpairheight}, $\prod_{u\in A}|\Phi(u)|^{n_{i_u}}=1$ and for every equivalence class $\alpha\in A/\sim$, we have 
\begin{equation}\label{equationaequsum}
\sum_{u\in \alpha}n_{i_u}h(i_u)=0.
\end{equation}

\medskip

For every $i=1,\dots,r$, define $\sG_i:=\{\sigma\in \sG_K|\,\, (\sigma,i)\in A\}.$ Because $T=\{1,\dots,r\},$
$\sG_i\neq\emptyset.$
If $(f_i,a_i)\sim (f_j,a_j)$, then $\sG_i=\sG_j.$
Hence for every $j=1,\dots,l$ the set $\sG_s$ does not depend on $s\in J_j$. We denote it by $\sG^j.$ 
Then 
$$|\sG^j|\sum_{s\in J_j} n_s\hat{h}_{f_s}(a_s)=\sum_{u\in A, i_u\in J_j}n_{i_u}h(i_u).$$
By Equality \ref{equationaequsum}, we only need to show that that for
$(\sigma, s)\sim (\delta,t)$, if $(\sigma, s)\in \{u\in A|\,\, i_u\in J_j\}$, then $(\delta,t)\in  \{u\in A|\,\, i_u\in J_j\}.$
Because $(\sigma, s)\in \{u\in A|\,\, i_u\in J_j\}$, $s\in J_i$. Because $(\sigma, s)\sim (\delta,t)$, $(\id, s)\sim_w(\id, t).$
By our assumption, we have $(\id, s)\sim (\id, t),$ which shows that $t\in J_i$. This concludes the proof.
\endproof

\medskip

Proposition \ref{proequpairheight} and Proposition \ref{prohproductalg} directly imply the following result on $\Q$-linear relations of canonical heights.
\begin{cor}\label{corlinearrelationheight}
Let $(f_i,a_i),\dots, (f_r,a_r)\in (\sD_d\times \overline{\Q})$ be dynamical pairs.
Set $T:=\{i=1,\dots, r|\,\, (f_i,a_i)\in \sT_d\}$ and $U:=\{i=1,\dots, r|\,\, (f_i,a_i)\in \sD\sP_d\setminus \sT_d\}$
Let $T=\sqcup_{j=1}^l J_j, U=\sqcup_{m=1}^n M_m$ be the partition associated to $\sim_w.$
Let $n_1,\dots, n_r\in \Z.$
Then the following holds:
\begin{points}
\item 
If for every $j=1,\dots,l, m=1,\dots,n,$
$$Q((f_s,a_s),s\in J_j)\in \{\sum_{s\in J_j} n_sz_s=0\} \text{ and } Q((f_s,a_s),s\in M_m)\in \{\sum_{m\in M_m} n_mz_m=0\},$$
then $\sum_{i=1}^r\hat{h}_{f_i}(a_i)^{n_i}=0.$
\item
If $U=\emptyset$, and $\sum_{i=1}^r\hat{h}_{f_i}(a_i)^{n_i}=0,$
then for every $j=1,\dots,l, m=1,\dots,n,$
$$Q((f_s,a_s),s\in J_j)\in \{\sum_{s\in J_j} n_sz_s=0\}.$$
\end{points}

In particular, if $(f_i,a_i),\dots, (f_r,a_r)\in \sT$, $\hat{h}_{f_i}(a_i), i=1,\dots,r$ are linearly independent over $\Q$ if and only if $(f_i,a_i)\not\sim (f_j,a_j)$ for every $i\neq j.$
\end{cor}

Example \ref{exenongeorelheight} shows that the condition in (i) is not necessary.

\subsubsection{Discussion}
Now we discuss the linear relations between canonical heights of dynamical pairs.
For every pair $\alpha=(f,a)\in \sD$, set $h(\alpha):=\hat{h}_f(a).$
Write $\sD_d=\sP_d\sqcup \sT_d\sqcup (\sD\sP_d\setminus \sT_d).$
Because $h(\sP_d)=0$, the pairs in $\sP_d$ are not interesting for our discussion. 

\medskip

Let $\langle h(\sT_d) \rangle$ and $\langle h((\sD\sP_d\setminus \sT_d)) \rangle$ be the $\Q$-spaces generated by $\sT_d$ and $\sD\sP_d\setminus \sT_d$ respectively.
By Corollary \ref{corlinearrelationheight}, these two space are $\Q$-linearly independent. So we only need to study the linear relations in $h(\sT_d)$ and in $h((\sD\sP_d\setminus \sT_d))$ separately.

\medskip

Corollary \ref{corlinearrelationheight} completely classified the $\Q$-linear relations between $h(\sT_d)$. After Galois conjugations, all the relations are come from geometry.

\medskip

For $h((\sD\sP_d\setminus \sT_d))$, the situation is different. The condition in (i) of Corollary \ref{corlinearrelationheight} is not necessary. In other words, the $\Q$-linear relations between $h(\sT_d)$ are not all come from geometry even up to Galois conjugations.
\begin{exe}\label{exenongeorelheight}
Let $f=z(z+1/2)$, $g=z(z+3/8)$, $a=b=1/16.$ Easy to check that $(f,a), (g,b)\in \sD\sP_d\setminus \sT_d.$
We first view them as complex dynamical systems via the unique embedding $\Q\hookrightarrow \C.$
Because the unique critical point of $f$ (resp. $g$) in $\C$ is contained in the attracting basin of the attracting fixed point $0$ of $f$ (resp. $g$),
$0$ is the unique attracting periodic point of $f$ (resp $g$). Because the multiplers $1/2$ and $3/8$ for $f$ and $g$ are multiplicatively independent,
$f\not\sim g.$ Hence $(f,a)\not\sim (g,b).$ Because $(f,a)$ and $(g,b)$ are defined over $\Q$, $(f,a)\not\sim_w (g,b).$
Easy to check that $\hat{g}_{f,v}(a)=\hat{g}_{g,v}(b)=0$ for all $v\in M_{\Q}\setminus \{2\}$ and $\hat{g}_{f,2}(a)=\hat{g}_{g,2}(b)=4\log 2.$
Hence $\hat{h}_f(a)=\hat{h}_g(b)=4\log 2$.
\end{exe}
However, as said in Remark \ref{remcomputable}, the nonarchimedean Green functions are computable for dynamical pairs. 
For pairs in $\sD\sP_d\setminus \sT_d$, their archimedean green functions are $0$ and their nonarchimedean green functions are computable, hence their canonical heights are computable.  By Lemma \ref{lemnonarchialgebraic}, those canonical heights take form
\begin{equation}\label{equsumplogp} \sum_{p \text{ primes }} c_p\log p
\end{equation} 
where $c_p$ are rational numbers which are $0$ except finitely many. Because $\log p, p \text{ primes }$ are $\Q$-linearly independent,
it is easy to find all $\Q$-linear relations between finitely many real numbers having form (\ref{equsumplogp}).

\section{Invariant germs of curves at infinity}\label{sectioninvargerms}
In the section, $K$ is of characteristic zero. 

\subsection{Some classes of polynomial endomorphisms}\label{subsectionpolyspace}
Let $\bk$ be a field of characteristic zero.
For $N\geq 2, d\geq 2$, denote by $\sP(N,d)$ the space of dominant endomorphisms $f: \A^N\to \A^N$ taking form 
$$f:(x_1,\dots,x_N)\mapsto (f_1(x_1,\dots, x_N),\dots, f_N(x_1,\dots,x_N))$$ of algebraic degree $\deg_1(f)=d$.
Recall that the algebraic degree of $f$ is $$\deg_1(f):=\max\{\deg(f_1),\dots, \deg(f_N)\}.$$
Write $f_i=\sum_{|I|\leq d}a_{i,I}x^I, i=1,\dots, N$
 where $I=(i_1,\dots,i_N)\in \Z_{\geq 0}^N, x^I=x_1^{i_1}\dots x_N^{i_N}$ 
 and $|I|=i_1+\dots+i_N.$

The Jacobian $J_{N,d}$ of $f$ is a polynomial in variables $x_1,\dots, x_N$ and $a_{i, I}, i=1,\dots, N, |I|\leq d$ of degree at most $d^N$ in $x_1,\dots, x_N$.
So we may write 
$$J_{N,d}=\sum_{|J|\leq d^N}b_{J}(a_{i,I}, i=1,\dots, N, |I|\leq d)x^J,$$
where $b_{J}$ are polynomials in $a_{i,I}, i=1,\dots, N, |I|\leq d$.  
Then we get
$$\sP(N,d)=\Spec (\bk[a_{i,I}, 1\leq i\leq N, |I|\leq d])\setminus ((\cup_{|J|\leq d^N}\{b_J=0\})\cup(\cap_{1\leq i\leq N, |J|=d}\{a_{i, J}=0\})),$$
which is an irreducible quasi-affine variety of dimension $N\binom{d+N+1}{N}$. 

\medskip

For general $N\geq 2$, we will introduce two spaces $\sP^*(N,d)$ and $\sP^{NS}(N,d)$ of $\sP(N,d)$. The relations among them are as follows:
$$\sP^{NS}(N,d)\subseteq_{\text{as a dense $G_{\delta}$-set}} \sP^*(N,d)\subseteq_{\text{as a dense open subset}} \sP(N,d).$$ 
When $N=2$, we introduce another subspace $\sP^{NP}(2,d).$ We have 
$$\sP^{NS}(2,d)\subseteq_{\text{as a dense $G_{\delta}$-set}} \sP^{NP}(2,d)\subseteq_{\text{as a dense open subset}} \sP^*(2,d)\subseteq \sP(2,d).$$
In Section \ref{sectioninvargerms}, we will prove algebricity results for invariant analytic germs of curves for  
$f\in\sP(N,d)$.  Moreover, we will show that there are a lot of invariant analytic germs of curves when $f\in\sP^*(N,d)$.
For $f\in \sP^{NS}(N,d)$, we will show that the $f$-periodic curves are of degree at most 2. This leads to a criteria for the transcendence of them. 
In Section \ref{sectionendomorphismplan}, we will focus on those $f\in\sP^{NP}(2,d)$ and classify such $f$ having infinitely many invariant curves.

\medskip
\subsubsection{Extendable endomorphisms}
Via the standard embedding $\A^N\hookrightarrow \P^N,$
every $f\in \sP(N,d)$ extends to a rational self-map on $\P^N.$ 
Denote by $\sP^*(N,d)$ the space of $f\in \sP(N,d)$ whose extension is an endomorphism.
This space is a Zariski dense open subset of $\sP(N,d)$. 
Indeed, $\sP^*(N,d)$ can be viewed as follows:
Let $H_{\infty}:=\P^N\setminus \A^N$ be the hyperplane at infinity.
Then $f$ induces a rational self-map $\overline{f}:=f|_{H_{\infty}}$ on $H_{\infty}$.
Then $f\in \sP^*(N,d)$ if and only if $f|_{H_{\infty}}$ is everywhere well-defined.
For each $i=1,\dots,N$, let $\overline{f_i}=\sum_{|I|=d} a_{i,I}x^I$ be the leading terms of $f_i.$
Then $\overline{f}$ is the endomorphism $H_{\infty}=\P^{N-1}\to \P^{N-1}$ sending 
$[x_1:\dots: x_N]$ to $[\overline{f_1}:\dots: \overline{f_N}].$
Denote by $\Rad(N-1,d)$ the space of rational self-maps on $\P^{N-1}$ of degree $d.$
It is clear that $\Rad(N-1,d)=\Proj\, \overline{K}[a_{i,J}, 1\leq N, |J|=d]=\P^{N\binom{d+N}{N-1}}.$
The morphism $\Phi:\sP(N,d)\to \Rad(N-1,d)$ 
defined by $f\mapsto \overline{f}$ is surjective.
Observe that $\overline{f}$ is everywhere well defined if and only if the zero set 
$$I_{\overline{f}}:=\{\overline{f_1}=\dots=\overline{f_N}=0\}\subseteq \P^{N-1}$$ is empty. 
This happens if and only if $\overline{f_i}\neq 0$ for $i=1,\dots, N$ and the intersections of the divisors 
$\{\overline{f_i}=0\}, i=1,\dots,N$ on $\P^{N-1}$ is a proper intersection.
Hence the space of endomorphisms of $\P^{N-1}$ of degree $d$
$$\End(N-1,d):=\{\overline{f}\in \Rad(N-1,d)|\,\, I_{\overline{f}}\neq\emptyset\}$$ is Zariski dense and open in $\Rad(N-1,d).$
Then $$\sP^*(N,d)=\Phi^{-1}(\End(N-1,d))$$ is a Zariski dense open subset of $\sP(N,d)$. 
\subsubsection{Endomorphisms without supperattracting periodic points} 
For $\overline{f}\in \End(N-1,d)$, we say that $\overline{f}$ has  \emph{NS property} if for every $n\geq 1$, every fixed point $x$ of $\overline{f}^n$, $d\overline{f}|_x$ is invertible. 
For every $f\in \sP^*(N,d)$, we say that $f$ has \emph{NS property at infinity} if $\overline{f}$ has \emph{NS property}, and define  
$\sP^{NS}(N,d)$ to be the set of $f\in\sP^*(N,d)$ of has NS property at infinity.
For every $f\in \sP^*(N,d)$, by Proposition \ref{proamplifiedfix}, $\Fix(f^n)$ is finite for every $n\geq 0$. 
Hence for every $n\geq 1$, one get a function $\theta_n$ on $\sP^*(N,d)$ sending $f$ to  $\prod_{x\in \Fix(f^n)}\det(df^n|_x).$
Here we count the fixed points with multiplicities. Because $\theta_n$ does not vanish at a general point of $\sP^*(N,d)$,
its zero set $Z_n:=\{\theta_n=0\}$ is a proper closed subset of $\sP^*(N,d)$.  Then we have 
$$\sP^{NS}(N,d)=\sP^*(N,d)\setminus (\cup_{n\geq 1}Z_n).$$
As a consequence, a very general $f\in \sP^*(N,d)$ satisfies the NS property at infinity. 
When $\bk$ is uncountable, this implies that most of $f\in \sP^*(N,d)$ satisfies the NS property at infinity. 
When $\bk$ is countable, the notion ``very general" does make much sense:
even a property is satisfied by a very general point, it may not be satisfied by any $\overline{\bk}-$point.
An alternative is to use the notion of ``adelic general", which was introduced in \cite[Section 3.1]{Xie2019}.
We will show that the NS property at infinity is satisfied by an adelic general point in Section \ref{subsectionadelidenseness}.

\subsubsection{Endomorphisms of NP type}
In this section, we focus the case $N=2$.

For an endomorphism $\overline{f}: \P^1\to \P^1$ of degree $d\geq 2$. A point $x\in \P^1(\overline{\bk})$ is called \emph{exceptional}, if its inverse orbit $\cup_{n\geq 0}\overline{f}^{-n}(x)$ is finite. It is well known that for every $\overline{f}$ there are at most $2$ exceptional points. 
We say that $\overline{f}$ is \emph{of NP-type} if it does not have any exceptional point.  Easy to see that $\overline{f}$ is of NP type if and only if $\overline{f}^n$ is of NP type for some $n\geq 1.$ Moreover, for two endomorphisms $\overline{f}, \overline{g}$ of $\P^1$ of degree $d\geq 2.$ Then $\overline{f}$ is of NP type if and only if $g$ is of NP-type.

Let $\Rad^{NP}(1,d)$ be the space of $f\in \Rad(1,d)$ which is of NP-type.  It is a Zariski dense open subset of $\Rad(1,d)$. Indeed, for $f\in \Rad(1,d)$, it is not of NP-type, if and only if there is a fixed point $x$ of $f^2$ satisfying $f^{-2}(x)=x.$ This a closed condition. So $\Rad^{NP}(1,d)$ is Zariski open in $\Rad(1,d)$. Easy to check that the map $z\mapsto (z^2+0.001)/(0.001 z^2+1)$ is of NP-type, hence $\Rad^{NP}(1,d)$ is non-empty.

For $f\in \sP^*(2,d)$, we say that $f$ is \emph{NP at infinity}, if $\overline{f}$ is of NP type. Then $\sP^{NP}(2,d):=\Phi^{-1}(\Rad^{NP}(1,d))$ is the space of $f\in \sP^*(2,d)$ which are NP at infinity, Because $\Phi|_{\sP^*(2,d)}$ is dominant, $\sP^{NP}(2,d)$ is a Zariski dense open subset of $\sP^*(2,d)$.

\subsection{Adelic denseness of NS endomorphisms}\label{subsectionadelidenseness}
In this section, we will show that the NS property at infinity is satisfied by an adelic general point.

\medskip

We first gather some of the properties satisfied by the {\bf{adelic topology}}. See \cite[Section~3]{Xie2019} for its precise definition and more properties.
Assume that $\bk$ is an algebraically closed extension of $\overline{\Q}$ of finite transcendence degree.
Let $X$ be a variety over $\bk$. The adelic topology is a topology on $X(\bk)$, which is defined by considering all 
possible embeddings of $\bk$ in $\C$ and $\C_p$, for all primes $p$.

An impotent example of adelic open subsets is as follows:  Let $L$ be a subfield of $\bk$ such that its algebraic closure ${\overline{L}}$ is equal to~$\bk$,$L$  is finitely generated over $\Q$, and $X$ is defined over $L$ i.e. $X=X_L\otimes_L \bk$ for some variety $X_L$ over $L$. Fix any embedding $\tau: L\hookrightarrow \C_p$(resp. $\C$).
Then,  given any open subset $U$ 
of $X_L(\C_p)$ for the $p$-adic (resp. euclidean) topology, the union $X_L(\tau, U):=\cup_\iota \Phi_\iota^{-1}(U)$ for all embeddings $\iota\colon \bk \to \C_p$ extending $\tau$ is, by definition, an open subset of $X(\bk)$ in the adelic topology. Moreover $X_L(\tau, U)$ is empty if and only if $U=\emptyset.$

The adelic topology has the following basic properties.
\begin{enumerate}
\item It is stronger than the Zariski topology. If $\dim(X)\geq 1$, there are non-empty, adelic, open subsets $\mathcal{U}$ and $\mathcal{U}'$ of $X(\bk)$
such that $\mathcal{U}\setminus \mathcal{U}'$ is Zariski dense in $X$.
\item It is $T_1$, that is for every pair of distinct points $x, y\in X(\bk)$ there are adelic open subsets $\mathcal{U}$, $\mathcal{V}$ of $X(\bk)$ such that 
$x\in \mathcal{U}, y\not\in \mathcal{U}$ and $y\in \mathcal{V}, x\not\in \mathcal{V}$.
\item Morphisms between algebraic varieties over $\bk$ are continuous in the adelic topology.
\item Flat morphisms are open with respect to the adelic topology.
\item The irreducible components of $X(\bk)$ in the Zariski topology coincide with the irreducible components of $X(\bk)$ in the  adelic topology.
\item Let $L$ be a subfield of $\bk$ such that (a) its algebraic closure ${\overline{L}}$ is equal to~$\bk$, (b) $L$  is finitely generated over $\Q$, and (c) $X$ is defined over $L$.
Endow the Galois group  $\Gal(\bk/L)$ with its profinite topology. Then the action 
$\Gal(\bk/L)\times X(\bk)\to X(\bk)$ is continuous for the adelic topology.
\end{enumerate}
\begin{rem}
From (5), when $X$ is irreducible, the intersection of finitely many nonempty adelic open subsets of $X(\bk)$ is nonempty.
\end{rem}
We say that a property $P$ holds for an {\bf{adelic general point}} if there exists an adelic dense open subset $\mathcal{U}$ of $X(\bk)$, 
such that $P$ holds for all points in $\mathcal{U}$.
\begin{thm}\label{thmadelicns}
Assume that $\bk$ is an algebraically closed extension of $\overline{\Q}$ of finite transcendence degree.
For $N\geq 2$, an adelic general $f\in \sP^*(N,d)$ (resp. $f\in \Rad(N-1,d)$) satisfies the NS property at infinity (resp. NS property). 
\end{thm}
\proof[Proof of Theorem \ref{thmadelicns}]
By (3), $\Phi$ is continuous w.r.t. the adelic topology. Since $\Phi(\sP^*(N,d))$ contains a Zariski dense open subset of $\Rad(N-1,d)$, by (1) and (5), we only need to show that 
an adelic general $f\in \Rad(N-1,d)$ satisfies the NS property.

Let $L$ be a subfield of $\bk$ such that its algebraic closure ${\overline{L}}$ is equal to~$\bk$, $L$  is finitely generated over $\Q$.
Because $\Rad(N-1,d)$ is defined over $\Q$, it is defined over $L$.

\begin{lem}\label{lemgoodendo}There is an endomorphism $g\in \Rad(N-1,d)$, a prime $p>0$ a field embedding $\tau: \bk\hookrightarrow \C_p$, such that 
$g\otimes_{\bk}\C_p$ has good reduction i.e. it extends to a finite endomorphism on $\P^N_{\C_p^{\circ}}$ and  for every fixed point $\widetilde{x}\in \Fix(\widetilde{g}^n), n\geq 1$ where $\widetilde{g}$ is the specialization of $g\otimes_{\bk}\C_p$,  $d\widetilde{g}^n|_{\widetilde{x}}$ is invertible.
\end{lem}

The  $g$ in the above lemma satisfies the NS property. Our idea is to should that a small perturbation of $g$, w.r.t. the topology induced by $\tau$, still satisfies the NS property.
One may write $g=[g_1:\dots: g_N]$ where $g_i, i=1,\dots, N$ are contained in $\tau^{-1}(\C_p^{\circ})\cap \bk$ and the leading coefficients of $\tau(g_i), i=1,\dots, N$ are of norm $1.$
Via $\tau$, view $g$ as an endomorphisms over $\C_p$.  
Consider the $p$-adic open neighborhood $U$ of $g$ in $\Rad(N-1,d)(\C_p)$
of endomorphisms taking form $[x_1:\dots :x_N]\to [g_1+F_1(x_1,\dots, x_N):\dots :g_N+F_N(x_1,\dots, x_N)]$ where $F_i, i=1,\dots,N$
are homogenous polynomials of degree $d$ with coefficients in $\C^{\circ\circ}=\{a\in \C_p|\,\, |a|_p<1\}.$
For every $f\in U$, we have $\widetilde{f}=\widetilde{g}.$ Hence for every fixed point $\widetilde{x}\in \Fix(\widetilde{f}^n), n\geq 1$ ,  $d\widetilde{f}^n|_{\widetilde{x}}$ is invertible.
This implies that for every $f\in \Phi_{\tau}^{-1}(U)$, $f$ satisfies the NS property. Because the the NS property is an algebraic property,  for every $f\in \Phi(\sP^*(N,d))_L(\tau|_{L}, U)$, $f$ satisfies the NS property, which concludes the proof.
\endproof

\proof[Proof of Lemma \ref{lemgoodendo}]
Because $\overline{\Q}\subseteq \bk,$ we may assume that $\bk=\overline{\Q}.$ 
I learned the following lemma from Charles Favre.
\begin{lem}\label{lemlatte}For every $d\geq 2$, there is an endomorphism $h: \P^1\to \P^1$ of Latt\`es type of degree $d$ defined over $\overline{\Q}$.
\end{lem}
An endomorphism $h$ of Latt\`es type is hyperbolic PCF i.e. all critical points of $h$ are preperiodic, but not periodic. 
We may assume that $h$ is defined over a number field $K$. There is a finite set of places $S$ of $K$ containing $M_K^{\infty}$ such that $h$ extends to a finite endomorphism $h_{O_S}$ of 
$\P^1_{O_S}$ over $O_S$ and for every place $v\in M_K^{\infty}\setminus S$, the specialization $\widetilde{h}$ of $h$ satisfies the following property: all critical points of $\widetilde{h}$ are preperiodic, but not periodic. In other words, for every fixed point $\widetilde{z}\in \Fix(\widetilde{h}^n)(\overline{\widetilde{K_v}}), n\geq 1$, $d\widetilde{h}^n|_{\widetilde{x}}\neq 0.$
Here $O_S$ is the ring of $S$-integers in $K.$
Consider the diagonal endomorphism $h\times\dots\times h: (\P^1)^{N-1}\to (\P^1)^{N-1}.$ The symmetric group $\Sigma_{N-1}$ acts on $(\P^1)^{N-1}$ by permuting the the coordinates. This action commutes with the endomorphism $h\times\dots\times h$. 
By \cite[Corollary 2.6]{Maakestad2005}, one gets $(\P^1)^{N-1}/\Sigma_{N-1}\simeq \P^{N-1}.$ Hence $h\times\dots\times h$ descends to an endomorphism $g$ on $(\P^1)^{N-1}/\Sigma_{N-1}\simeq \P^{N-1}.$ Moreover it extends to a finite endomorphism $(h\times\dots\times h)_{O_S}=h_{O_S}\times_{O_S}\dots\times_{O_S} h_{O_S}$ on $(\P^1_{O_S})^{N-1}.$
After enlarging $S$, we may assume that $(\P^1_{O_S})^{N-1}/\Sigma_{N-1}\simeq \P^{N-1}_{O_S}$ and $g$ extends to a finite endomorphism $g_{O_S}$ on $\P^{N-1}_{O_S}.$
For every closed point $y\in (\P^1_{O_S})^{N-1}$ which is fixed by $(h\times\dots\times h)_{O_S}^n$ for some $n\geq 1$, $(h\times\dots\times h)_{O_S}$ induces an \'etale endomorphism on the local ring $\sO_{(\P^1_{O_S})^{N-1},y}.$ Let $y'$ be the image of $y$ in $\P^{N-1}_{O_S}$. It is fixed by $g^n_{O_S}$.
Indeed all periodic points of $g_{O_S}$ are image of periodic points of $(h\times\dots\times h)_{O_S}$
By \cite[Proposition 41.20.6]{stacks-project}, the morphism $g_{O_S}$ induces an \'etale endomorphism on the local ring $\sO_{\P^{N-1}_{O_S},y'}.$
Pick $v\in M_K\setminus S$ and  let $p$ the prime number associated to $v$. Then $v$ induces an embedding $\tau': K\hookrightarrow \C_p,$ which extends to an embedding $\tau: \bk\hookrightarrow \C_p.$ Then for every fixed point $\widetilde{x}\in \Fix(\widetilde{g}^n), n\geq 1$, 
$d\widetilde{g}^n|_{\widetilde{x}}$ is invertible.
\endproof

\proof[Proof of Lemma \ref{lemlatte}]First assume that $d=m^2$ a prefect square. Pick any elliptic curve $E$ over $\overline{\Q}$, let $[m]: E\to E$ be the multiplication by $m.$
Then $[m]$ descents to an endomorphism $g$ on $\P^1=E/[\pm 1].$ Then $\deg g=m^2=d.$
Otherwise, $d=m^2l$ where $m\in \Z$ and $l$ is a square free integer.  Let $E$ be any elliptic curve over $\overline{\Q}$ with $\End(E)\otimes \Q=\Q(\sqrt{-l})$. Then multiply by $m\sqrt{-l}$ defines a complex multiplication on $E.$ It descents to an endomorphism $g$ on $\P^1=E/[\pm 1].$ Then $\deg g=|m\sqrt{-l}|^2=d.$
\endproof

\subsection{Algebraicity of invariant branches of curves at infinity}
Let $d\geq 2$ and $f: \A^N\to \A^N$ be an endomorphism in $\sP(N,d)(K).$ 
We embed $\A^N$ in $\P^N$ and denote by $H_{\infty}:=\P^N\setminus \A^N$ the hyperplane at infinity.
Let $g_v, v\in M_K$ and $h$ be the naive Green functions and naive height on $\A^N$ respectively.

\medskip

Let $S$ be a finite set of places $S\subseteq M_K$ containing $M_K^{\infty}$ such that $f$  extends  to an endomorphism $f_{O_S}: \A^N_{O_S}\to \A^N_{O_S}$ where $O_S$ is the ring of $S$-integers in $K.$
Let $(v,o_v, C_v)$ be a branch of curve at infinity defined over $K.$  
We note that if $C_v$ is irreducible and not contained in $H_{\infty}$, $C_v\cap (H_{\infty})_v$ is finite and it contains $o_v.$
We say that $C_v$ is $f$-invariant if $f(C_v\setminus\{o_v\})\subseteq C_v\setminus\{o_v\}.$ if $C_v$ is $f$-invariant, then $f(C_v)$ contains an open neighborhood of $o_v$ in $C_v$. 

The following proposition is directly implied by the Northcott property of $K$ and Theorem \ref{thmcurvesalge}.
\begin{pro}\label{prospndkattr}
Let $x$ be a point in $\A^N_{O_S}(O_S)$ which is not $f$-preperiodic. 
Assume that for $v\in S$, there is a branch of curve $(v,o_v, C_v)$ at infinity defined over $K$ which is $f$-invariant and containing 
$x.$ Then the Zariski closure of $O_f(x)$ has dimension $1$. Moreover, for $v\in S$, if $g_v(f^n(x)), n\geq 0$ is unbounded and 
$C_v$ is irreducible at $o_v$, then $C_v$ is algebraic over $K.$
\end{pro}

If further $f\in \sP(N,d)(K)$, we get a result which applies for a single branch of curve.

\begin{pro}[=Proposition \ref{prospndkextendedattrintro}]\label{prospndkextendedattr}
Let $x\in \A^N(K)$ and $v\in M_K$ such that $g_v(f^n(x))\to \infty$ as $n\to \infty.$
Let $(v,o_v, C_v)$  be an irreducible branch of curve at infinity defined over $K$ which is $f$-invariant and containing $x$.
Then $C_v$ is algebraic over $K$ at $o_v.$ 
\end{pro}
\proof[Proof Proposition \ref{prospndkextendedattr}]
It is well known that there is $B_v>0$ such that for every $y\in K_v^N$ with $g_v(y)\geq B_v$, we have 
$$\lim_{n\to \infty}g_v(f^n(y))/d^n\in (0, +\infty).$$
On the other hand, because $x$ is not preperiodic, its canonical height exists and positive i.e. 
$$\lim_{n\to \infty}h(f^n(x))/d^n\in (0,+\infty).$$
Then we conclude the proof by Theorem \ref{thmcurvesalge}.
\endproof

\subsection{Periodic points at infinity}
Let $d\geq 2$ and $f: \A^N\to \A^N$ be an endomorphism in $\sP^*(N,d)(\bk).$ 
For $n\geq 0$, let $F'(f^n)$ be the set of fixed points $x\in H_{\infty}$ of $\overline{f}^n: H_{\infty}\to H_{\infty}$ such that $d\overline{f}^n|_x$ is invertible. 
We note that when $f\in \sP^{NS}(N,d)$, $F'(f^n)=\Fix(\overline{f}^n)$ for every $n\geq 1.$
The following lemma shows that even in the general case, most of the fixed points of $\overline{f}^n$ are contained in $F'(f^n)$.
\begin{lem}\label{lemfixpfix}For $f\in \sP^*(N,d)(\bk),$ we have $$|F'(f^n)|/|\Fix(\overline{f^n})|\to 1 \text{ and }  |\Fix(\overline{f^n})|/d^{(N-1)n}\to 1$$ for $n\to \infty.$
Moreover $\cup_{n\geq 1}F'(f^n)$ is Zariski dense in $H_{\infty}.$
\end{lem}

\proof[Proof of Lemma \ref{lemfixpfix}]
Fix an embedding $\tau: \bk\hookrightarrow \C.$ Via $\tau$, we view $f$ as a complex dynamical system.
For $n\geq 0$, let $F''(f^n)$ be the set of repelling fixed points of $\overline{f^n}$. Then we have $$F''(f^n)\subseteq F'(f^n)\subseteq \Fix(\overline{f^n}).$$
By Lefschetz fixed point theorem, we have $$\limsup_{n\to\infty}|\Fix(\overline{f^n})|/d^{(N-1)n}\leq 1.$$
By \cite[Theorem 1.1]{Dinh2015}, $$\liminf_{n\to\infty}|F'(f^n)|/d^{(N-1)n}\geq \lim_{n\to\infty}|F''(f^n)|/d^{(N-1)n}= 1.$$
Then we get $$|F'(f^n)|/|\Fix(\overline{f^n})|\to 1 \text{ and }  |\Fix(\overline{f^n})|/d^{(N-1)n}\to 1$$ for $n\to \infty.$
We now prove that $\cup_{n\geq 1}F'(f^n)$ is Zariski dense in $H_{\infty}.$

We note that $\overline{f}(F'(f^n))=F'(f^n)$ for every $n\geq 0.$ Let $Z$ be the Zariski closure of $\cup_{n\geq 1}F'(f^n)$, then we have $\overline{f}(Z)=Z.$
After replacing $f$ by a suitable positive iterate, we may assume that every irreducible component of $Z$ is fixed by $\overline{f}$.
Let $Z_i, i=1,\dots, r$ be all irreducible component of $Z$. Set $f_i:=\overline{f}|_{Z_i}: Z_i\to Z_i$.
Then $f_i^*(O(1)|_{Z_i})=O(d)|_{Z_i}$.  
\begin{lem}\label{lemcountfixpoint}Let $g: X\to X$ be an endomorphism of a projective variety. Let $L$ be an ample line bundle on $X$ satisfying $g^*L=L^{\otimes d}$ for some $d\geq 2.$ Then there is $C>0$ such that $|\Fix(g^n)|\leq Cd^{n\dim X}$ for every $n\geq 0.$
\end{lem}
Then we have $\sum_{i=1}^r |\Fix(g_i^n)|\leq C'd^{n\dim Z}$ for some $C'>0.$ On the other hand, we have  $\sum_{i=1}^r |\Fix(g_i^n)|\geq |F'(f^n)|\geq cd^{n(N-1)}, n\geq 1$ for some $c>0.$
Then we get $\dim Z\geq N-1$, which implies that $Z=H_{\infty}.$
\endproof

\begin{rem}
Lemma \ref{lemcountfixpoint} holds in any characteristic. 
\end{rem}

\proof[Proof of Lemma \ref{lemcountfixpoint}]
After replacing $X$ by its normalization, we may assume that $X$ is normal. 
We denote by $\Delta$ the diagonal of $X\times X$ and $\Gamma_n, n\geq 1$ the graph of $g^n.$
Let $\pi_i: X\times X\to X$ the projection to the $i$-th coordinate. We have $\Fix(g^n)=\pi_1(\Gamma_n\cap \Delta), n\geq 1.$
Since $\Fix(g^n)$ is finite and $\pi_1|_{\Delta}$ is an isomorphism,  we get 
$|\Fix(g^n)|= |\Gamma_n\cap \Delta|.$ Because $X$ may be singular, the intersection of $\Gamma_n$ and $\Delta$ may be not well-defined. That is why we insist of counting fixed points without multiplicities. 
Set $M:=\pi_1^*L\otimes \pi_2^*L$. It is ample on $X\times X.$
Let $\sI_{\Delta}$ be the ideal sheaf of $\Delta$. 
There is $a>0$ such that $M^{\otimes a}\otimes \sI_{\Delta}$ is generated by global section. 
\begin{lem}\label{lemdeltagamnhyp}We have $$|\Gamma_n\cap \Delta|\leq a\Gamma_n\cdot c_1(M)^{\dim X}.$$
\end{lem}
We now compute $\Gamma_n\cdot c_1(M)^{\dim X}.$
Set $D:=c_1(L), D_i:=\pi^*_i(D), i=1,2$. Then $c_1(M)=D_1+D_2.$
We have 
$$\Gamma_n\cdot c_1(M)^{\dim X}=\Gamma_n\cdot (D_1+D_2)^{\dim X}=\sum_{i=0}^{\dim X}\binom{\dim X}{i}(\Gamma_n\cdot D_1^{\dim X-i}\cdot D_2^i)$$
$$=\sum_{i=0}^{\dim X}\binom{\dim X}{i}((g^n)^*D^i\cdot D^{\dim X-i})=\sum_{i=0}^{\dim X}\binom{\dim X}{i}d^{ni}(D^i\cdot D^{\dim X-i})$$
$$=(d^n+1)^{\dim X}(D^{\dim X}).$$
Then we get $$|\Fix(g^n)|\leq a(d^n+1)^{\dim X}(D^{\dim X}),$$
which concludes the proof.
\endproof
\proof[Proof of Lemma \ref{lemdeltagamnhyp}]
Let $s_1,\dots, s_{\dim X}$ be general elements in $H^0(M^{\otimes a}\otimes \sI_{\Delta})\subseteq H^0(M^{\otimes a}).$ 
Let $H_i$ be the zero locus of $s_i\in H^0(M^{\otimes a}).$ 
Then $H_1\cap\dots\cap H_{\dim X}\cap \Gamma_n$ is zero dimensional. Because $\Delta\subseteq H_1\cdots H_{\dim X}$, we get
$$|\Gamma_n\cap \Delta|\leq |H_1\cap\dots\cap H_{\dim X}\cap \Gamma_n|\leq \Gamma_n\cdot (H_1\cdots H_{\dim X})=a\Gamma_n\cdot c_1(M)^{\dim X},$$
which concludes the proof.
\endproof

\subsection{Existence of invariant branches of curves at infinity}
In this section, set $\bk=\overline{K}.$
Let $d\geq 2$ and $f: \A^N\to \A^N$ be an endomorphism in $\sP^*(N,d)(\bk).$ 

\begin{lem}\label{lemformalcurve}For every $o\in F'(f)$, there is a unique irreducible formal curve $\widehat{C}$ at $o$ which is $f$-invariant and not contained in $H_{\infty}$. It intersects $H_{\infty}$ transversally at $o$.
If both $f$ and $o$ are defined over $K$, then $\widehat{C}$ is defined over $K$ at $o$.  
\end{lem}

\proof[Proof of Lemma \ref{lemformalcurve}]
In some suitable coordinates, we may assume that $o=[0:\dots: 0:1]$ and $H_{\infty}=\{x_1=0\}.$

{\bf Step 1:}
We first show that such a formal curve exists after some base changing $\bk$ by a field extension $L$ over $\bk.$
Let $L:=\bk((t))$ with the $t$-adic norm. Then $f$ induces an endomorphism $f_L$ on $\P^N_L$. Let $U$ be the affinoid subdomain of $(\P^N_L)^{\an}$ defined by 
$U:=\{[x_1:\dots: x_{N-1}: 1]|\,\, |x_i|\leq |t|, i=1,\dots, N-1\}.$
Then $f_L$ induces an endomorphism $g:=f_L|_U: U\to U$ which fixes the hyperplane $Y:=(H_{\infty})_L^{\an}\cap U.$
Because $d\overline{f}|_o$ is invertible, $g|_Y: Y\to Y$ is an automorphism. Let $\widetilde{g},\widetilde{U},\widetilde{Y}$ be the reduction of $g$, $U$ and $Y$. Then $\widetilde{g}(\widetilde{U})=\widetilde{Y}.$ Then, by \cite[Theorem 8.3]{Xie2019}, there is a unique morphism $\phi: U\to Y$ such that $\phi\circ g=g|_Y\circ \phi.$ 
In particular, $\widetilde{\phi}=\widetilde{g|_Y}^{-1}\circ \widetilde{g}.$ One may check that $d\widetilde{g}$ has rank $N-1$ everywhere, so $d\widetilde{\phi}$ has rank $N-1$ everywhere. Hence, for every point $x\in Y(L)$, $\phi^{-1}(x)$ is a smooth irreducible curve which intersects $Y$ transversally.
Because $o$ is fixed by $g$, $\phi^{-1}(o)$ an analytic curve passing through $o$, which is invariant under $g.$ 
Its completion at $o$ defines a formal curve $\widehat{C_L}$ with coefficients in $L$. It is invariant under $f_L$.

\medskip

{\bf Step 2:} Show the uniqueness of such formal curves even after any base of $\bk.$ 
Let $\bk'$ be a field extension of $\bk$.
Let $\widehat{C_1}, \widehat{C_2}$ be two irreducible formal curves passing through $o$ in $\P^N_{\bk'}$ which are invariant under $f_{\bk'}$ and not contained in $H_{\infty}.$
Let $L:=\bk'((t))$ with the $t$-adic norm. Then $f$ induces an endomorphism $f_L$ on $\P^N_L$. Let $U$ be the affinoid subdomain of $(\P^N_L)^{\an}$ defined by 
$U:=\{[x_1:\dots: x_{N-1}: 1]|\,\, |x_i|\leq |t|, i=1,\dots, N-1\}.$ We get exactly the same situation as in the previous paragraph. We keep the notation in the previous paragraph.
We note that $\widehat{C_1}$ and $\widehat{C_2}$ induce closed analytic curves $C_1,C_2$ in $U$ respectively.  We claim that both $C_1, C_2$ equal to $\phi^{-1}(o).$ This implies that $\widehat{C_1}=\widehat{C_2}.$ Now we prove the claim. We only need to show it for $C_1.$ Assume that $\phi(C_1)\neq o$. Then $\phi(C_1)$ is an irreducible curve in $Y.$ 
In particular, $\phi(C_1)$ is $g|_Y$-invariant, compact and connected.  Set $C':=\cap_{n\geq 0}g|_Y^n(\phi(C_1)).$ Then $C'$ is compact and is invariant under $g|_Y$ and $g|_Y^{-1}.$ 
By the following lemma, $C'$ is connected. 
\begin{lem}\label{leminsconn}Let $Y_n, n\geq 0$ be a sequence of decreasing connected compact sets. Then $\cap_{n\geq 0}Y_n$ is connected. 
\end{lem}
We also note that $o\in C'$. 
\begin{lem}\label{lemcpotherpoint}We have $C'\neq \{o\}.$
\end{lem}
Then $C'$ contains infinitely many points. 
For every $y\in C'$ and $n\geq 0$, we have $g|_Y^{-n}(y)\in C'\in \phi(C)$.
Hence, there is $y_n\in C$ such that $g|_Y^{-n}(y)=\phi(y_n).$ By (vi) of \cite[Theorem 8.3]{Xie2019}, we get 
$g^n(y_n)\to y$ as $n\to \infty$. Because $C$ is $g$-invariant, $g^n(y_n)\in C$ for $n\geq 0$. Since $C$ is closed, $y\in C$. Then we get $C'\subseteq C.$
Because $C'\subseteq Y$ and $C'$ is infinite, we get $C\subseteq Y$ which is a contradiction. 

{\bf Step 3:} Assume that $f,o$ are defined over $K$, show that there is a such formal curve over $K.$ Note that such a curve is unique if exists by Step 2 and it intersects $H_{\infty}$  transversally by Step 1. 
Let $L$ be an algebraically closed field extension of $K$ such that there is an irreducible formal curve $C_L$ in $\P^N_L$ passing through $o$ which is $f_L$-invariant and is not contained in $(H_{\infty})_L$. We only need to show that $C_L$ is defined by an ideal $I_L$ of $\widehat{O}_{\P^N_L,o}$ which is defined over $K.$ Set $G:=\Gal(L/K).$ For $\sigma\in G$, $\sigma(C_L)$ is also an irreducible formal curve in $\P^N_L$ passing through $o$ which is $f_L$-invariant and is not contained in $(H_{\infty})_L$. By Step 2, we get that $\sigma(C_L)=C_L$. It follows that $I_L$ is $G$-invariant.  Let $m$ be the local ring of $o\in \P^N_K$ and $m_L:=m\otimes_K L$.  For every $n\geq 1$,
$(I_L+m_L^n)/m_L^n$ is an ideal of $\widehat{O}_{\P^N_L,o}/m_L^n=(\widehat{O}_{\P^N_K,o}/m^n)\otimes_KL$ which is $G$-invariant. By Galois descent, there is an ideal $J_n$ of $\widehat{O}_{\P^N_K,o}/m^n$ such that $(I_L+m_L^n)/m_L^n=J_n\otimes_KL.$ It follows that $I_L$ is defined over $K$. This concludes the proof.
\endproof
\proof[Proof of lemma \ref{leminsconn}]
Assume that $Y_{\infty}:=\cap_{n\geq 0} Y_n$ is not connected. Then there are disjoint open subsets $W_1,W_2$ of $Y_0$ such that $Y_{\infty}\subseteq W_1\sqcup W_2$ and $Y_{\infty}\cap W_i\neq \emptyset$ for $i=1,2.$ Set $W:=W_1\cup W_2$. Then $Y_0\setminus W=\cup_{n\geq 0}((Y_0\setminus W)\setminus (Y_n\setminus W)).$ Since $Y_0\setminus W$ is compact and $(Y_0\setminus W)\setminus (Y_n\setminus W), n\geq 0$ are open in $Y_0\setminus W$, there is $m\geq 0$ such that $Y_n\setminus W=\emptyset$ for $n\geq m.$
Then $Y_m$ is not connected, which is a contradiction. 
\endproof
\proof[Proof of Lemma \ref{lemcpotherpoint}]

For  $l\geq 1$, let $V_l$ be the open neighborhood of $o$ in $Y$ defined by $V_l:=\{[0: x_2:\dots: x_{N-1}:1]|\,\, |x_i|<|t^l|, i=2,\dots, N-1\}.$
Since $g|_Y$ is an automorphism of $Y$ and $g|_Y(o)=o$, we have $g(V_l)=V_l$. There is $l\geq 1$ such that $D:=\phi(C_1)\setminus V_l\neq \emptyset$.
Then $g|_Y(D)\subseteq D.$ Since $\cap_{n\geq 0}g|_Y^n(D)$ is an intersection of decreasing compact sets, it is not empty. Moreover it is contained in $C'\setminus V_l$. This concludes the proof.
\endproof

\begin{lem}\label{lemconvergence}
Keep the notations in Lemma \ref{lemformalcurve}. For every $v\in M_K^{\infty}$, there is an irreducible analytic curve $C_v$ at $o$ such that  the completion of $C_v$ at $o$ is induced by $\widehat{C}$.
\end{lem}
\proof[Proof of Lemma \ref{lemconvergence}]
By \cite[Theorem 3.1.4]{Abate2001}, for every $v\in M_K^{\infty}$, there is an irreducible analytic curve $C_v$ at $o$ which not contained in $H_{\infty}$ and is $f$-invariant. 
Then it induces an irreducible formal curve $\widehat{C_{v,o}}$ of $\P^N_{K_v}$ passing through $o$, which not contained in $H_{\infty}$ and is $f_{K_v}$-invariant. 
The uniqueness part of Lemma \ref{lemformalcurve} implies that $\widehat{C_{v,o}}$
is induced by $\widehat{C}$.
\endproof

Assume further $f\in \sP^{NS}(N,d)$, the we have the following criteria for $f$-periodic curves.
\begin{cor}\label{corbounddegree}Let $f\in \sP^{NS}(N,d)$ and let $C$ be a periodic curve of $f$. If $C\not\subseteq H_{\infty}$, then $\deg C\leq 2.$
\end{cor}
\proof[Proof of Corollary \ref{corbounddegree}]
After replacing $f$ by a suitable positive iterate, we may assume that $f(C)=C.$
Because $C\setminus H_{\infty}$ is an affine curve and  $f(C\setminus H_{\infty})=C\setminus H_{\infty}$ and $\deg f|_C=d\geq 2$, $C$ has at most two branches at infinity. 
After replacing $f$ by $f^2$, we may assume that all branches of $C$ at infinity are fixed by $f$. Because each branch of $C$ defines an irreducible $f$-invariant formal passing through 
some point in $\Fix(\overline{f})=F'(f)$ which is not contained in $H_{\infty},$ by Lemma \ref{lemformalcurve}, it intersects $H_{\infty}$ transversally. Hence $C\cdot H_{\infty}$ equal to the number of branches of $C$ at infinity, which is at most $2.$
\endproof

\section{Endomorphisms of $\A^2$}\label{sectionendomorphismplan}
In this section, $\bk$ is an algebraically closed field of characteristic $0.$
Let $f: \A^2\to \A^2$ be an endomorphism in $\sP^{*}(2,d)$ for some $d\geq 2.$ It extends to an endomorphism of $\P^2$. Let $H_{\infty}:=\P^2\setminus \A^2$ be the hyperplane at infinity. 

\medskip

We say that $f$ is \emph{homogenous}, if it takes form $f: (x,y)\mapsto (F(x,y), G(x,y))$ where $F,G$ are homogenous polynomials of degree $d.$
We say that $f$ is \emph{homogenous at $o=(o_1,o_2)\in \A^2(\bk)$} if $f$ is homogenous in the new coordinates $x'=x-o_1, y'=y-o_2.$
\begin{rem}\label{remuniqueorigin}
Easy to see that if $f$ is homogenous at $o$, then such $o$ is unique. Moreover we have $f^{-1}(o)=o.$
\end{rem}

If $f$ is homogenous at $o$, then for every point $x\in \Per(\overline{f})\subseteq H_{\infty}$, the line $L_{o,x}$ passing through $o$ and $x$ is $f$-periodic.
Hence $f$ has infinitely many periodic curves.
The aim of this section is to prove the following result which shows that its inverse also holds when $f\in \sP^{NP}(2,d)$.
\begin{thm}[=Theorem \ref{thminfinitepercurveintro}]\label{thminfinitepercurve}For $f\in \sP^{NP}(2,d)$ with $d\geq 2.$ If there are infinitely many $f$-periodic curves, then $f$ is homogenous at some point $o\in \A^2(\bk)$.
Moreover, all but finitely many $f$-periodic curves are lines passing through $o$.
\end{thm}
\begin{rem}\label{remstrongcomp}
This result strengthens  the $N=2$ case of Corollary \ref{corbounddegree} in two directions.
First, $\sP^{NP}(2,d)$ is strictly larger than $\sP^{NS}(2,d)$. Second, in Theorem \ref{thminfinitepercurve}, we get a complete classification of $f\in \sP^{NP}(2,d)$  having infinitely many $f$-periodic curves.
\end{rem}

If $f$ is not of NP-type at infinity, Theorem \ref{thminfinitepercurve} may not be true.
\begin{exe}\label{exenotnp}
Let $f$ be the endomorphism taking form $f=g\times h$ where $g$ is a polynomial of degree $d$.
Then for every $a\in \Per(g)$, the line $\{x=a\}$ is periodic for $f.$
The same, for every $b\in \Per(h)$, the line $\{y=b\}$ is periodic for $f.$
Hence $f$ has infinitely many periodic curves. But it is easy to see that $f$ is homogenous at some point $o\in \A^2$ if and only if 
$g$ and $h$ are conjugate to $z\to z^d$ by some affine transformations of $\A^1.$
\end{exe}

The proof of the Theorem \ref{thminfinitepercurve} relies on the theory of valuative tree introduced by Favre and Jonsson in \cite{Favre2004}, which was developed by Favre, Jonsson and the author in \cite{Favre2007,Favre2011,Xie2015ring, Xie2017a, Xie2017}. 
\subsection{Valuative tree at infinity}\label{sectiontvaluationtreeinfty}
 We refer to \cite{Jonsson} for details, see also \cite{Favre2004,Favre2007,Favre2011}.
%
%
Denote by $V_{\infty}$ the space of all normalized valuations centered at infinity normalized valuation centered at infinity i.e. the set of functions $v:\bk[x,y]\rightarrow \mathbb{R}\cup \{+\infty\}$ such that 
\begin{points}
\item $v(k^*)=0$;
\item $v(fg) = v(f) + v(g)$ and $v(f +g) \ge \min \{v(f), v(g) \}$ for $f,g\in \bk[x,y]$
\item  $\min\{v(x),v(y)\}=-1$.
\end{points}
The topology on $V_{\infty}$ is defined to be the weakest topology making the map $v\mapsto v(P)$ continuous for every $P\in k[x,y]$.
Under this topology, $V_{\infty}$ is compact.

\smallskip

The set $V_{\infty}$ is equipped with a {\em partial ordering} defined by $v\leq w$ if and only if
$v(P)\leq w(P)$ for all $P\in k[x,y]$. Then  $-\deg: P \mapsto -\deg (P)$ is the unique minimal element.

\smallskip

For $v\in V_\infty\setminus \{-\deg\}$, the set
$ \{ w \in V_\infty, \,| - \deg \le w \le v \}$ is isomorphic as a poset to the real segment
$[0,1]$ endowed with the standard ordering. In other words, $(V_\infty, \le)$ is a rooted tree in the sense of \cite{Favre2004,Jonsson}.
Given any two valuations $v_1,v_2\in V_{\infty}$,
there is a unique valuation in $V_{\infty}$ which is maximal in the set $\{v\in V_{\infty}|\,\,v\leq v_1 \text{ and } v\leq v_2\}.$ We denote it by  $v_1\wedge v_2$.
The segment $[v_1, v_2]$ is by definition the union of $\{w \,| \, v_1\wedge v_2 \le w \le v_1\}$
and $\{w \,| \, v_1\wedge v_2 \le w \le v_2\}$.

\smallskip

Pick any valuation  $v\in V_\infty$. We say that two points $v_1, v_2$
lie in the same direction at $v$ if the segment $[v_1, v_2]$ does not contain $v$.
A {\em direction} (or a tangent vector) at $v$ is an equivalence class for this relation.
We write $\Tan_v$ for the set of directions at $v$.

\smallskip

When $\Tan_v$ is a singleton, then $v$ is called an endpoint\index{endpoint}. In $V_\infty$, the set of endpoints is exactly the set of all maximal valuations. 
When $\Tan_v$ contains exactly two directions, then $v$ is said to be regular.
When $\Tan_v$ has more than three directions, then $v$ is a branch point.

\smallskip

Pick any $v\in V_\infty$.
For any tangent vector $\vec{v}\in \Tan_v$, denote by $U(\vec{v})$ the subset of those elements in $V_\infty$ that determine $\vec{v}$. This is an open set whose boundary is reduced to the singleton $\{v\}$. If $v\neq -\deg$, the complement of $\{w\in V_\infty, \,| w \ge v\}$ is equal to $U(\vec{v}_0)$ where $\vec{v}_0$ is the tangent vector determined by $-\deg$.

It is a fact that finite intersections of open sets of the form $U(\vec{v})$ form a basis for the topology of $V_\infty$.

There are four types of valuations in $V_{\infty}$. However, we only need one of them in this paper, which is called curve valuations.
Recall that $H_{\infty}$ is the line at infinity of $\mathbb{P}^2_k$.
For any formal curve $s$ centered at some point $q\in H_{\infty}$ which is not contained in $H_{\infty}$, denote by $v_s$ the valuation defined by $P\mapsto (s\cdot l_{\infty})^{-1}\ord_{\infty}(P|_s)$. Then we have $v_s\in V_{\infty}$ and call it a {\em curve valuation}.  We note that the branch $s$ is determined by its associated valuation $v_s.$
All curve valuations are maximal in $V_{\infty}.$

Let $C$ be a curve in $\P^2$ which is not $H_{\infty}$. Assume that $C$ has a branch $C_q$ at a point $q\in H_{\infty}$. Then $C_q$ defines a curve valuation in $V_{\infty}.$

\medskip

Recall that there is a {\em skewness} function.
$\alpha:V_{\infty}\rightarrow [-\infty,1]$ 
It has the following properties:
\begin{points}
\item $\alpha$ is strictly decreasing, and upper semicontinuous;
\item $\alpha$ is continuous on segments;
\item $\alpha(-\deg)=1;$
\item for every curve valuation $v$, $\alpha(v)=-\infty$;
\item every valuation satisfying $\alpha(v)=-\infty$ is maximal.
\end{points}
\subsubsection{Computation of local intersection numbers of curves at infinity}\label{sectionlocalinter}
Let $s_1$, $s_2$ be two different formal curves at infinity.
We denote by $(s_1\cdot s_2)$ the intersection number of these two formal curves in $\P^2$. This intersection number is always nonnegative, and it is positive if and only if $s_1$ and $s_2$ are centered at the same point. When $s_1,s_2$ are not contained in $H_{\infty}$, 
denote by $v_{s_1},v_{s_2}$ the curve valuations associated to $s_1$ and $s_2.$
By \cite[Proposition 2.2]{Xiec}, we have 
\begin{equation}\label{equationdegloc}
(s_1\cdot s_2)=(s_1\cdot H_{\infty})(s_2\cdot H_{\infty})(1-\alpha(v_{s_1}\wedge v_{s_2})).
\end{equation}
\subsection{Local valuative tree}
See \cite{Jonsson, Favre2004} for details.
For a point $q\in \P^2$, there exist local coordinates $(z,w)$ at $q$ such that $\widehat{O_{q}}=\bk[[z,w]]$ with the maximal ideal $m:=(z,w)$.

Let $V_q$ be the space of normalized valuations centered at $q$, i.e. 
the functions $v: \bk[[z,w]]\to [0,+\infty]$ such that
\begin{points}
\item $v(k^*)=0;$
\item for $f,g\in \bk[[z,w]]$, $v(fg)=v(f)+v(g)$ and $v(f+g)\geq \min\{v(f),v(g)\}$;
\item $\min\{v(z),v(w)\}=1.$
\end{points}
The topology on $V_q$ is defined to be the weakest topology making the map $v\mapsto v(P)$ continuous for every $P\in k[[z,w]]$.
Under this topology, $V_q$ is compact.
The space $V_{q}$  is
equipped with a partial ordering defined by $v\leq w$ if and only if $v(f)\leq w(f)$ for all $f\in k[[z,w]]$ for which is again a real tree (see \cite{Favre2004,Favre2007,Jonsson}).
The order of vanishing $\ord_q$ at the point $q$ is a valuation in $V_q$. Moreover, it is the unique minimal point in $V_q.$

As in the global case, for every formal curve $s$ at $q$, it defines a \emph{curve valuation} $v_s$ sending $P\in k[[z,w]]$ to $m_s^{-1}\ord_q(P|_{s})$, where $m_s$ is the multiplicity of $s$ at $q.$

Every endomorphism $g^*$ of $\bk[[z,w]]$ induces a map $g_{\d}: V_q\to V_q$ by 
$$g_{\d}: v\mapsto (P\mapsto \min\{\ord_q(g^*z),\ord_q(g^*w)\}^{-1}v(g^*P)).$$
This action on $V_q$ is continuous.  Moreover every point in $V_q$ has at most finitely many preimages.
\begin{rem}\label{remfiniteboundpull}Let $U$ be an open subset of $V_q$. Then $g_{\d}(\partial (g_{\d}^{-1}(U)))\subseteq \partial U.$
For every connected component $V$ of $\partial (g_{\d}^{-1}(U))$, we have $\partial V\subseteq \partial (g_{\d}^{-1}(U))$. Hence 
$\partial V\subseteq g_{\d}^{-1}(\partial U).$ If $\partial U$ is finite, then $\partial V$ is finite. 
\end{rem}
For a branch of curve $s$ at $q$, let $g(s)$ be its image under the map induced by $g^*$. Then $g_{\d}(v_s^q)=v^q_{g(s)}.$

\medskip

Recall that there is a {\em local skewness} function
$\alpha_q:V_{q}\rightarrow [1,+\infty]$.
It has the following properties:
\begin{points}
\item $\alpha$ is strictly increasing, and lower semicontinuous;
\item $\alpha$ is continuous on segments;
\item $\alpha(\ord_{q})=1;$
\item for every curve valuation $v$, $\alpha(v)=+\infty$;
\item every valuation satisfying $\alpha(v)=+\infty$ is maximal.
\end{points}

At last we describe the connection between the local valuative tree and the global one.
Assume that $q\in H_{\infty}.$
There exists a valuation $v^q_{H_{\infty}}$ defined by $P\mapsto \ord_q(P|_{H_{\infty}})$ for $P\in k[[z,w]].$

Denote by $U(q)$ the set of
valuations in $V_{\infty}$ whose center in $X$ is $q$ and set $\overline{U(q)}:=U(q)\cup \{-\deg\}$. For any $v\in U(q)$, there exists $r_q(v)\in \mathbb{R}^+$ such that $r_q(v)v\in V_q.$ Set $v^q:=r_q(v)v$ when $v\in U(q)$ and $v^q:=v_{H_{\infty}}^q$ when $v=-\deg.$
The map $\phi_q: \overline{U(q)}\rightarrow V_q$ defined by $v\mapsto v^q$ is a homeomorphism. Moreover for a branch of curve $s$ at $q$ not contained in $H_{\infty}$,  we have
$\phi_q(v_s)=v^q_s.$
By \cite[Lemma 11.6]{Xie2017a}, we have $$\phi_q(\{v\in \overline{U(q)}|\,\, \alpha(v)=-\infty\})=\{v\in V_q|\,\, \alpha_q(v)=+\infty\}.$$

\subsection{Proof of Theorem \ref{thminfinitepercurve}}
Let $f\in \sP^{NP}(2,d)$ with $d\geq 2.$ Let $S$ be the set of supperattracting fixed points of $\overline{f}$. It is exactly the union of orbits of periodic critical points of $\overline{f}$, hence $S$ is finite.

We denote by $PC(f)$ the set of irreducible $f$-periodic curves in $\P^2$. Write $PC(f)=PC(f)_o\sqcup PC(f)_s \sqcup \{H_{\infty}\}$, where $PC(f)_s$ is the set of $C\in PC(f)\setminus \{H_{\infty}\}$ with $C\cap H_{\infty}\neq\emptyset.$ This decomposition is $f$-invariant. Elements in $PC(f)_s$ are called special.  

\begin{lem}\label{lemboundbranches}For every $C\in PC(f)\setminus \{H_{\infty}\}$, it has at most two branches at infinity.
\end{lem}
\begin{rem}This lemma holds for every $f\in \sP(N,d)$ in any dimension $N$.
\end{rem}
\proof[Proof of Lemma \ref{lemboundbranches}]
After replacing $f$ by a suitable positive iterate, we may assume that $f(C)=C.$
Because $C\setminus H_{\infty}$ is an affine curve and  $f(C\setminus H_{\infty})=C\setminus H_{\infty}$ and $\deg f|_C=d\geq 2$, $C$ has at most two branches at infinity. 
\endproof

The following lemma bounds the special invariant curves. 
\begin{lem}\label{lemboundexccurve}The set $PC(f)_s$ is finite.
\end{lem}
\proof[Proof of Lemma \ref{lemboundexccurve}]
After replacing $f$ by a suitable positive iterate, we may assume that all points in $S$ are fixed points. 
At every fixed point of $f$, we still denote by $f$ the germ of $f$ at $q$ and $f_{\d}$ the induced endomorphism on $V_q.$

For every $q\in S,$ the branch $H_{\infty}^q$ satisfies $f^*H_{\infty}^q=dH_{\infty}^q$ and $f_*H_{\infty}^q=r_qH_{\infty}^q$ where $r_q$ is the ramification index of $\overline{f}: H_{\infty}\to H_{\infty}$ at $q.$ Since $\overline{f}$ is of NP-type, $r_q<d.$
By \cite[Lemma 5.7]{Xie2017a}, there is $w_q<v_{H_{\infty}^q}$ such that 
\begin{points}
\item $\overline{f_{\d}(W_q)}\subseteq W_q$ where $\{v\in V_q|\,\, v\vee v_{H_{\infty}^q}>w_q\}$;
\item $f_{\d}|_{W_q}^n\to v_{H_{\infty}^q}$ as $n\to \infty;$
\item $\{v\in V_q|\,\, \alpha_q(v)<+\infty\}\subseteq \cup_{n\geq 1}f_{\d}^{-n}(W_q)$.
\end{points}
We note that, for every $n\geq 0,$ the only $f_{\d}$-periodic point in $f_{\d}^{-n}(W_q)$ is $v_{H_{\infty}^q}.$
In particular, let $C$ be an $f$-periodic curve passing through $q$ and $s$ be a branch of $C$ at $q$, then 
$v_s^q\not\in \cup_{n\geq 1}f_{\d}^{-n}(W_q).$

Set $K:=\{v\in V_{\infty}|\,\, \alpha(v)\geq -3\}.$ It is a compact. 
Then $K_q:=\phi_q(K\cap \overline{U(q)})\subseteq V_q$ is compact. 
We note that $$K_q\cap \{v\in V_q|\,\, \alpha_q(v)=+\infty\}=\{v_{H_{\infty}^q}\}\subseteq W_q.$$
The compactness of $K_q$ implies that there is $N_q\geq 0$ such that 
$K_q\subseteq f_{\d}^{-N_q}(W_q).$  By Remark \ref{remfiniteboundpull}, the boundary of $f_{\d}^{-N_q}(W_q)$ is finite.

Then $K\cap \overline{U(q)}$ is contained in the unique connected component $Y_q$ of $\phi_q^{-1}(f_{\d}^{-N_q}(W_q))$ containing $-\deg.$
Since $Y_q$ has finite boundary, it takes form 
$$Y_q=U(q)\setminus (\cup_{w\in B_q}\{v\geq w\})$$ where $B_q=\partial Y_q\setminus \{-\deg\}$.
We note that for distinct points $w_1,w_2\in B_q$, $w_1\not\geq w_2.$
Then for every periodic branch of curve $s$ at $q$, $v_s\geq w(s)$ for exactly one $w(s)\in B_q.$ 
Because $K_q\subseteq Y_q$, $\alpha(w)<-3$ for every $w\in B_q.$

For $C\in PC(f)_s$, because it has at most two branches at infinity, there is one branch $s_C$ at infinity with $s_C\cdot H_{\infty}\geq \deg(C)/2.$
If $C$ has one branch $s$ at infinity which does not meet $S$, then by Lemma \ref{lemformalcurve}, $(s\cdot H_{\infty})=1$. So we may choice $s_C$ to meet $S.$
Then we get $s_C\geq w(C)$ for a unique $w(C)\in U(q)$ where $q\in S$ is the center of $s.$
We claim that for $C_1,C_2 \in PC(f)_s$, if $w(C_1)=w(C_2)$, then $C_1=C_2.$
Otherwise $C_1\neq C_2,$ then by Equality (\ref{equationdegloc}), we have
$$\deg(C_1)\deg(C_2)=(C_1\cdot C_2)\geq s_{C_1}\cdot s_{C_2}=(s_{C_1}\cdot H_{\infty})(s_{C_1}\cdot H_{\infty})(1-\alpha(v_{s_{C_1}}\wedge v_{s_{C_2}}))$$
$$\geq (s_{C_1}\cdot H_{\infty})(s_{C_1}\cdot H_{\infty})(1-\alpha(w(C_1))>4(s_{C_1}\cdot H_{\infty})(s_{C_1}\cdot H_{\infty})\geq \deg C_1\deg C_2,$$
which is a contradiction. Then we get $|PC(f)_s|\leq \sum_{q\in S}|B_q|$ is finite. 
\endproof

\medskip

\proof[Proof of Theorem \ref{thminfinitepercurve}]Since $PC(f)$ is infinite and $PC(f)_s$ is finite, $PC(f)_o$ is infinite. The proof of Corollary \ref{corbounddegree} shows that every $C\in PC(f)_o$ has degree at most two and each brach of $C$ at infinity intersects $H_{\infty}$ transversally.  So for every $C\in PC(f)_o$, $\deg C$ is exactly the number of branches of $C$ at infinity. Moreover, by the uniqueness part of Lemma \ref{lemformalcurve}, it is also $|C\cap H_{\infty}|.$
Then for every $C\in PC(f)_o$,  $\deg f(C)=\deg C.$ 
Moreover for $C_1,C_2\in PC(f)_o$, $C_1=C_2$ if and only if $(C_1\cap H_{\infty})\cap (C_2\cap H_{\infty})\neq \emptyset.$

Let $B$ be the set of $C\in PC(f)_o$ with $\deg C=2.$ 
Set $e:=1$ if $B$ is finite and $e=2$ if $B$ is infinite. 
Set $A=B$ if $B$ is infinite and $A=PC(f)_o\setminus B$ if $B$ is finite. 
Then $A$ is infinite, $f$-invariant and for every $C\in A$, $\deg C=e\in \{1,2\}.$
We note that if $e=1$, then $PC(f)\setminus A$ is finite.
Set $R:=\cup_{C\in A}(C\cap H_{\infty})\subset \Per(\overline{f}).$ Then $R$ is infinite.

\medskip

For any $n\geq 1,$ set $M_n:=\P(H^0(\P^2, O(n)))$. It is the space of effective divisors of degree $n$ in $\P^2.$
Then $f$ induces a morphism $f_*: M_n\to M_{dn}$ sending $D$ to $f_*D.$
There is another morphism $\psi : M_n\to M_{nd}$ sending $D$ to $dD.$
We note that $\psi$ is an isomorphism from $M_n$ to its image. 
Then $A$ can be viewed as a subset of $M_e$. Let $Z$ be the Zariski closure of $A$ in $M_e.$
For every $C\in A$, $f_*(C)= df(C)$ and $f(C)\in A\subseteq Z$. Hence $f_*(Z)\subseteq \psi(M_n)$ and $\psi^{-1}\circ f_*(Z)\subseteq Z.$
For every $C\in Z$, we have $\psi^{-1}\circ f_*(C)=f(C).$ Then it induces an endomorphism $g: Z\to Z$ sending $C$ to $f(C).$
Set $\Sigma:=\{(x,C)\in \P^2\times Z|\,\, x\in C\}$. It is a closed in $\P^2\times Z.$ 
Set $$\Gamma:=\Sigma\cap (H_{\infty}\times Z)=\{(x,C)\in H_{\infty}\times Z|\,\, x\in C\}.$$
Let $\pi_i, i=1,2$ be the $i$-th projection of $H_{\infty}\times Z$. 
Then $R=\pi_1(\pi^{-1}_2|_{\Gamma}(A)).$

Since $R\subseteq \pi_1(\Gamma)$ is dense in $H_{\infty}$, 
for every irreducible component $\Gamma'$ of $\Gamma$ satisfying $\pi_1(\Gamma')=H_{\infty}$,
we have
\begin{equation}\label{equationgammpd}\pi_1(\Gamma')=H_{\infty}.
\end{equation}

\begin{lem}\label{lemopoc}
For every $q\in R$, $|\pi_2^{-1}(q)|=1.$
\end{lem}
By this lemma, for every $(x,C)\in \Gamma$, $x\in R$ if and only if $C\in A.$
Let $\Gamma''$ be the union of irreducible component $\Gamma_1$ of $\Gamma$ with $\pi_1(\Gamma_1)=H_{\infty}.$
Since $R$ is dense in $H_{\infty}$, the above lemma with Equation \ref{equationgammpd} implies that $\Gamma''$ is irreducible and $\pi_1|_{\Gamma''}:\Gamma''\to H_{\infty}$ is birational. 
Because $H_{\infty}$ is of dimension $1$ and normal, $\pi_1|_{\Gamma''}$  is an isomorphism. In particular $\Gamma''\simeq \P^1$. 
By Lemma \ref{lemopoc} again $$\Gamma''\cap \pi^{-1}_1(A)=\Gamma''\cap \pi^{-1}_1(R)=\Gamma\cap \pi^{-1}_1(R)=\Gamma\cap \pi^{-1}_2(A).$$
Because $A$ is dense in $Z$, $\pi_2(\Gamma'')=Z.$
For a general $C\in Z$, $|C\cap H_{\infty}|=e$ and $C$ intersects $H_{\infty}$ transversally.  Then $\pi_2|_{\Gamma''}: \Gamma''\to Z$ is surjective and of degree $e.$
Then $h:=\pi_2|_{\Gamma''}\circ (\pi_1|_{\Gamma''})^{-1}: \P^1\to Z$ is surjective of degree $e.$
 In particular, we get $\overline{f}\gtrsim_{h}g$, hence $g$ is not of NP-type. 
 
 \medskip
 
 Let $\Delta$ be the diagonal of $Z\times Z.$
 Let  $K$ be the Zariski closure of $$K':=\{(C_1,C_2,x)\in Z\times Z\times \P^2|\,\, C_1\neq C_2, x\in C_1\cap C_2\}$$ in $Z\times Z\times \P^2$.
 Set $$W:=\{(C_1,C_2,x)\in Z\times Z\times \P^2|\,\, x\in C_1\cap C_2\}.$$ 
 We have $K'\subseteq K\subseteq W.$
 Let $p_1: (Z\times Z)\times \P^2\to Z\times Z$ and $p_2: (Z\times Z)\times \P^2\to \P^2$ be the projections. 
 For $C_1,C_2\in Z$ with $C_1\neq C_2$, $1\leq |{p_1|_K}^{-1}((C_1,C_2))|\leq e^2$.
 So $p_1|_K: K\to Z\times Z$ is surjective and generically finite of degree at most $e^2.$ In particular, $\dim K=2.$ Moreover, for every irreducible component $K_1$ of $K$,
 $p_1(K_1)=Z\times Z.$ Note that $K$ is invariant under the endomorphism $(g\times g)\times f$ on $(Z\times Z)\times \P^2.$ Since every point in $p_1|_{K'}^{-1}(A\times A)$ is $(g\times g)\times f$-periodic, $K_1$ is also periodic. 

\medskip

For every $x\in p_2(K)$, let $Z_x$ be the set of $C\in Z$ with $x\in C$. 
Then $p_2|_{W}^{-1}(x)=Z_x^2$ and $p_2|_{K'}^{-1}(x)=Z_x^2\setminus \Delta_{Z_x}$ where $\Delta_{Z_x}$ is the diagonal of $Z_x^2.$ Then we get $$2\dim Z_x=\dim p_2|_{K'}^{-1}(x)\leq \dim p_2|_{K}^{-1}(x)\leq \dim p_2|_{W}^{-1}(x)=2\dim Z_x.$$  So $p_2|_{K}^{-1}(x)=0$ or $2.$ Hence $\dim p_2(K)=0$ or $2$. 

\medskip

We now prove that $\dim p_2(K)=0.$ Otherwise $\dim p_2(K)=2$. Let $K_1$ be an irreducible component of $K$ of dimension $2$ such that $p_2(K)=\P^2.$ Then $p_1(K_1)=Z_1^2$ and $K_1$ is $(g\times g)\times f$-periodic. There are Zariski dense open subset $V_1$, $V_2$ of 
$Z_1^2$ and $\P^2$ respectively such that $p_i|_{K_1}, i=1,2$ is finite flat over $V_i.$
Then $p_i(p_1^{-1}(V_1)\cap p_2^{-1}(V_2)), i=1,2$ are nonempty open subsets.
Because there is $C_0\in A$ such that $C_0\cap p_2(p_1^{-1}(V_1)\cap p_2^{-1}(V_2))\neq\emptyset.$ Pick $y\in C_0\cap p_2(p_1^{-1}(V_1).$
Pick $z\in p_1(p_2|_{K_1}^{-1}(y))$. 
Let $C'$ be an irreducible component of $p_1(p_2^{-1}(C_0))\subseteq Z^2$ containing $z$. Then $C'$ is a curve.
Let $C''$ be an irreducible component of $p_1|_K^{-1}(C')$ satisfying $p_1(C'')=C'$. We have $p_2(C'')=C_0.$
Then $C''$ is $(g\times g\times f)$-preperiodic. Moreover, for every $n\geq 0$, both $p_1((g\times g\times f)^n(C''))=(g\times g)^n(C')$ and $p_2((g\times g\times f)^n(C''))=C_0$ are curves.
Then there is $m\geq 0$ and $l\geq 1$ such that $D:=(g\times g\times f)^m(C'')$ is $(g\times g\times f)^l$-invariant. 
Then we get $(g\times g\times f)^l|_D\sim f^l|_{C_0}$ is of NP-type. Since $(g\times g\times f)^l|_D\sim (g\times g)^l|_{p_1(D)}$ and $(g\times g)^l|_{p_1(D)}\sim g^l$, $g$ is of NP-type. This is a contradiction.

\medskip

Now set $O':=p_2(K)$ which is finite. 
Then $K\subseteq Z^2\times O'$. There is a nonempty subset $O$ of $O'$ such that 
the irreducible components of $K$ of dimension $2$ are exactly $\{Z^2\times x|\,\, x\in O\}$.
In other words, $O=\cap_{C\in Z} C.$ For two different $C_1,C_2\in A$, $C_1\cap C_2\cap H_{\infty}=\emptyset.$ Hence $O\subseteq \A^2.$
\begin{lem}\label{lemopoint}We have that $|O|=1$ and $e=1.$ 
\end{lem}
Set $O=\{o\}.$ Then $Z$ is exactly the space of line passing through $o$. After a changing of origin, we may assume that $o=(0,0)\in \A^2(\bk).$ 
Since $f$ preserves the lines passing through the origin, $f$ is homogenous. Since $A\subseteq Z$ and $PC(f)\setminus A$ is finite, all but finitely many periodic curves of $f$ are lines passing through $o$.
\endproof
\proof[Proof of Lemma \ref{lemopoc}]
Otherwise, there is $C_1\in A$ and $C_2\in Z$ such that $q\in C_1\cap C_2.$
After replacing $f$ by a suitable positive iterate, we may assume that $q$ is $f$-fixed and the branch $s_1$ of $C_1$ at $q$ is fixed. Let $s_2$ be a branch of $C_2$ at $q.$
Since $f_*C_i=dC_i, i=1,2$, $f_* s_i=ds_i.$ Since $f$ has topological degree $d$ at $q$, $f^*s_1=s_1.$ In particular, $s_1$ is $f^{-1}$-invariant.
Since $s_1\neq s_2$ and they both pass through $q$, $f^n(s_2)\neq s_1$ for all $n\geq 0.$
Hence $(s_1\cdot f^n(s_2))\geq 1$ for all $n\geq 1.$ 
Then for every $n\geq 0$, we have $$(s_1\cdot s_2)=((f^n)^*s_1\cdot s_2)=(s_1\cdot f^n_{*}s_2)=d^n(s_1\cdot f^n(s_2))\geq d^n,$$
which is a contradiction. 
\endproof

\proof[Proof of Lemma \ref{lemopoint}]
Pick $E\in A$ such that $\{E\}\times Z\not\subseteq p_1(K\setminus (Z^2\times O))$. After replacing $f$ by a suitable positive iterate, we may assume that $f(E)=E.$ 
Let $C$ be a general point in $Z$. Then $C\cap E=O.$
Since $C$ is general, $g^{-1}(C)$ is exactly $d$ general points $C_1,\dots, C_d$.
We have $\cup_{i=1}^d C_i\subseteq f^{-1}(C).$ Since 
$$de=\deg (\cup_{i=1}^d C_i)\leq \deg f^{*}(C)=de,$$
we get $\cup_{i=1}^d C_i=f^{-1}(C).$
Then we get $$(f|_E)^{-1}(O)\subseteq f^{-1}(C)\cap E=(\cup_{i=1}^d C_i)\cap E=O.$$
Since $O\subseteq E\setminus H_{\infty}$ and all points in $E\cap H_{\infty}$ are exception. We get $e+|O|\leq 2.$
Hence $e=|O|=1.$
\endproof

\newpage
\bibliography{dd}
\end{document}